# A Classification of Six Functor Formalisms via Structured Spaces

Salash Tolan Nabaala

March 2025

## Introduction

In this paper, we accomplish three things. We lay out an $\infty$-categorical interpretation of reconstruction theorems which are germane to the *symmetric monoidal* perspective[1] of noncommutative algebraic geometry, present sufficient conditions which allow for the factorization of certain six functor formalisms through animated $S$-stacks, and give a "universal" six functor formalism through which the aforementioned six functor formalisms factor through. Furthermore, and what is arguably the main feat of this article, these achievements, though in appearance arising from disparate concerns, are realized in the dissipation of a familiar thematic tension: that between space and quantity.

    We mean here by a theory of $\infty$-categories, the theory of quasicategories developed by J. Lurie in [i] and by a six functor formalism, a lax symmetric monoidal functor of $\infty$-categories $\mathcal{D} : \mathrm{Corr}(C, E) \longrightarrow \mathsf{Pr}^L$ whose image has closed symmetric monoidal $\infty$-categories as objects. Here $\mathrm{Corr}(C, E)$ is the category whose objects are objects of $C$—a $(small)$[2] $\infty$-category which admits pullbacks—and whose morphisms are spans $Y \longleftarrow Z \longrightarrow X$ with the rightward arrow belonging to $E$; $E$ being a collection of morphisms of $C$ stable under base change, composition and containing all equivalences in $C$ (we shall, by convention, call the pair $(C, E)$ a *geometric setup*). The latter formalism is due to L. Mann and is developed in A.5 in [ii] and Lecture II-IV in [iii]. On the other hand, $\mathsf{Pr}^L$ is the subcategory[3] of the $\infty$-category of *(large)* $\infty$-categories spanned by presentable $\infty$-categories and cocontinuous maps between them. In what is to follow, unless we are exempted from ambiguity, we will explicitly state when an $\infty$-category is being considered.

    There is in this work three theorems of climactic importance. The first of these theorems establishes that particular types of functors between $\infty$-categories factor through animated $S$-stacks; that is, the $\infty$-category $\mathcal{P}(\mathsf{SCR}_S^{\mathrm{op}})$ which we will label $\mathsf{Stk}_S$ (cf. [iv] 4.1.2). Specifically, if we let $\mathfrak{X}^L \subseteq \mathsf{Pr}^L$ be a particular kind of (very large) self-dual $\infty$-topos and $\mathsf{Stk}_{S\,|}^{\mathrm{op}} \subseteq \mathsf{Stk}_S^{\mathrm{op}}$ some subcategory of animated $S$-stacks which we shall define, then we can say the following:

**Theorem A.** *Let $\mathcal{O}^{op} : \mathcal{T} \longrightarrow \mathfrak{X}^L$ be a suprematic space. Then there exists $(\_)^{\mathcal{O}} : \mathcal{T}'^{op} \longrightarrow \mathsf{Stk}_{S\,|}^{op}$ where $\mathcal{T}' \subseteq \mathcal{T}$, and $\mathcal{O}_0 : \mathsf{Stk}_{S\,|}^{op} \longrightarrow \mathfrak{X}^L$, extending $\mathcal{O}|\mathcal{T}'^{op} : \mathcal{T}'^{op} \longrightarrow \mathfrak{X}^L$ as $\mathcal{O}_0 \circ (\_)^{\mathcal{O}}$. Furthermore, $\mathcal{O}_0$ admits a section $\mu_{\mathcal{O}} : \mathcal{O}(\mathcal{T}'^{op}) \longrightarrow \mathsf{Stk}_{S\,|}^{op}$. (cf. 2.1.31).*

---

[1] The view that holds the symmetric monoidal structure of the ($\infty$-) categories being considered central to the subject; for example, derived algebraic geometry as seen from the lens of homotopical algebraic geometry.

[2] The use of *small* references set theoretic considerations of size which we shall keep in mind by fixing Grothendieck universes. The operative terms are going to be *small, large* and *very large*.

[3] We will use the convention that refers any sub-$\infty$-category as simply a subcategory.

A *suprematic space* is, given a pregeometry $\mathfrak{T}$, a $\mathfrak{T}$-structure on an $\infty$-topos $\mathfrak{X}$ meeting certain criteria which we shall establish. Equivalently, it is a particular kind of $\mathcal{G}$-structure on an $\infty$-topos $\mathfrak{X}$, where $\mathcal{G}$ is the geometry enveloping $\mathfrak{T}$; see [iv] 3.1 and 3.4 for an in-depth exposition on pregeometries and their geometric envelopes. Suprematic spaces constitute the spaces that shall be the mainstay of this article. Specifically, we shall consider a certain subcategory of the $\infty$-category of $\mathfrak{T}$-structures in $\mathfrak{X}^L$, $\mathsf{Str}_\mathfrak{T}(\mathfrak{X}^L)$, which we label $\mathsf{Sup}_\mathfrak{T}^\otimes(E, \mathfrak{X}^L)$ and whose objects will be suprematic spaces. Indeed, the second of the important theorems mentioned above concerns itself with the latter $\infty$-category and the $\infty$-category of six functor formalisms $\mathrm{Corr}(\mathsf{Stk}_{S|}, E) \longrightarrow \mathfrak{X}^L$, which we label[1] $\mathsf{Fun}^{\otimes,\mathrm{lax}}(\mathrm{Corr}(\mathsf{Stk}_{S|}, E), \mathfrak{X}^L)$.

**Theorem B.** *There exists $f : \mathsf{Sup}_\mathfrak{T}^\otimes(E, \mathfrak{X}^L)^{op} \longrightarrow \mathsf{Fun}^{\otimes,\mathrm{lax}}(\mathrm{Corr}(\mathsf{Stk}_{S|}, E), \mathfrak{X}^L)$ which is a fully faithful map of $\infty$-categories. (cf. 2.2.8).*

This theorem epitomizes the paper in that it displays a parametrization of a full subcategory of $\mathsf{Fun}^{\otimes,\mathrm{lax}}(\mathrm{Corr}(\mathsf{Stk}_{S|}, E), \mathfrak{X}^L)$ by suprematic spaces: which are, in other words, structured spaces of some definite kind. Additionally, we will show that $f^{\mathrm{op}}$ can be seen as factoring through a map $\mathsf{Sup}_\mathfrak{T}^\otimes(E, \mathfrak{X}^L) \longrightarrow \mathcal{S}\mathsf{tr}_{\mathsf{Stk}_{S|}}(\mathfrak{X}^\mathcal{L}) \subseteq \mathsf{Fun}(\mathsf{Stk}_{S|}^{\mathrm{op}}, \mathsf{CAlg}(\mathfrak{X}^L))^{\mathrm{op}}$. This will allow us to give a "universal" six functor formalism through which factor all the six functor formalisms lying in the essential image of $f$. That we are able to achieve this is borne out by observing that structured spaces admit a universal structured space (cf. [iv] 1.4.2). It is the statement of this appropriation of the latter fact that gives the last of our three important theorems.

**Theorem C.** *There exists a geometric setup $(\mathsf{Stk}_{S|}, \overline{E})$ and a lax $\infty$-symmetric monoidal map $\chi : \mathrm{Corr}(\mathsf{Stk}_{S|}, \overline{E}) \longrightarrow L_{\text{ét}}(\mathsf{Stk}_{S|})$, whose image is a subcategory of $\mathsf{Pr}^L$, such that for every $\mathcal{D} \in \mathsf{Fun}^{\otimes,\mathrm{lax}}(\mathrm{Corr}(\mathsf{Stk}_{S|}, \overline{E}), \mathfrak{X}^L)$ in the image of $f : \mathsf{Sup}_\mathfrak{T}^\otimes(\overline{E}, \mathfrak{X}^L)^{op} \longrightarrow \mathsf{Fun}^{\otimes,\mathrm{lax}}(\mathrm{Corr}(\mathsf{Stk}_{S|}, \overline{E}), \mathfrak{X}^L)$, there exists a map $\widetilde{\mathcal{D}} : L_{\text{ét}}(\mathsf{Stk}_{S|}) \longrightarrow \mathfrak{X}^L$ so that $\mathcal{D} \simeq \widetilde{\mathcal{D}} \circ \chi$. (cf. 3.2.18).*

In an extremely loose manner of speaking, this theorem can be seen as an actualization of Grothendieck's motivic dream. The shortcoming is of course that, other than not adhering to precise formulation, classically, motives are defined at the level of cohomological theories and not at the level of the (triangulated) categories giving rise to said cohomologies; and it remains in this case to be seen if the theorem descends into an instantiation of motives at the cohomological level (if and) when the fecundity of the involved $\infty$-categories allows for such a descent[2]. Nonetheless, the theorem remains worthwhile for us since it takes a shape evocative of the "yoga of motives"; if not an outright incarnation of the same.

Additionally, we present in the appendix work which will form the rudiments of a subterranean through which one may weave a common thread between these three theorems. This work comes to the foundational results of *tensor triangulated geometry* using topos theoretic methods. Indeed, in this light, and keeping in mind the theorems above, this article can be seen as a single continuous movement after the following questions: What is a point? What is a space? What is a motive?

### What is a Point?

The beginning for *tensor triangulated geometry* is P. Balmer's reconstruction of a reduced noetherian scheme $X$ from the symmetric monoidal category of perfect complexes on $X$, $(D_{\mathrm{perf}}(X), \otimes^{\mathbf{L}})$, in [v]; the backbone of this approach being R.W Thomason's classification of *thick* tensor subcategories of $D_{\text{perf}}(X)$. Balmer, then proceeding in [vi], improves on this result to obtain a noetherian scheme $X$ from $(D_{\text{perf}}(X), \otimes^{\mathbf{L}})$ by introducing the category of *classifying support data* associated to a symmetric monoidal triangulated category (cf. [vi] 3.1 and 5.1). The latter approach is further recast in [vii] using the language of ideal lattices by Buan-Krause-Solberg. The retrieval of a concentrated scheme $X$ from $(D_{\text{perf}}(X), \otimes^{\mathbf{L}})$ and a noetherian scheme $X$ from the symmetric monoidal category of coherent sheaves on $X$, $(\text{Coh}, \otimes_X)$, then follows: this being a redetermination of Gabriel's original reconstruction theorem. And, as a matter of fact, due to the flexibility brought about by the theory of ideal lattices, G.Garkusha in [viii] is able to give a reconstruction of $X$ from the symmetric monoidal category of quasicoherent sheaves on $X$, $(\text{QC}(X), \otimes_X)$. This gives a brief history which, although nowhere near exhaustive, captures what is ultimately the germ of tensor triangulated geometry.

---

[1] $\mathsf{Fun}(X, Y)$ is the simplicial enrichment of $\mathsf{Hom}_{\mathsf{Set}_\Delta}(X, Y)$.
[2] For example, when the objects of $L_{\text{ét}}(\mathsf{Stk}_{S|})$ and $\mathfrak{X}^L$ are stable $\infty$-categories or when both are $\infty$-topoi.

In functorial speak, this is that tensor triangulated geometry explicitly gives a contravariant functor $\text{TriCat}^{\otimes} \longrightarrow \text{LRS'}$ which is called the *Balmer spectrum*; a functor from symmetric monoidal triangulated categories and triangulated symmetric monoidal functors to locally ringed spaces and ringed maps. When this functor is restricted to derived categories of perfect complexes on concentrated schemes, the underlying schemes are recovered. Furthermore, this functor (roughly)[1] extends to a contravariant functor from symmetric monoidal exact categories and exact symmetric monoidal functors to ringed spaces and ringed morphisms, $\text{Exact}^{\otimes} \longrightarrow \text{RS}$. This covers Gabriel's reconstruction theorem and is in part the overall content of [vii] (cf. 7.2) and [viii].

If we look at the Balmer spectrum pointwise, we notice that the topological space underlying the image is the collection of all *prime* thick tensor subcategories given a certain topology supplied by the universal *classifying support datum* (cf. [vii] 4.2 and 5.2). Moreover, each category in the domain admits finite direct sums and a 'multiplication' in the form of the tensor product which distributes over direct sums; so that we may, in fact, identify each category with the semi-ring whose objects are equivalence classes of isomorphic objects (cf. [vii] 6.3) and whose operations are the former. And from here, we can arrive at *prime* thick tensor subcategories as certain kinds of prime ideals of this semi-ring. This procedure, when we restrict to those categories of the form $(D_{\text{perf}}(X), \otimes^{\mathbf{L}})$ for $X$ concentrated, is reminiscent of the classical construction of an affine scheme. We will, at least set theoretically speaking, make this analogy precise by formulating in a functorial context the prime ideals of a commutative ring and the prime thick tensor subcategories of a subclass of the categories in question. The decisive observation is the following.

**Proposition 0.1.0.** *Let $R$ be a commutative ring considered as a semigroup under multiplication and $i : I \hookrightarrow R$ the inclusion of a proper ideal $I \subset R$ considered as a sub-semigroup inclusion in the category of semigroups and semigroup maps. Then the pullback of multiplication along the inclusion, $I \times_R (R \times R)$, is isomorphic to $(I \times R) \cup (R \times I)$, if and only if $I$ is a prime ideal of $R$.*

*Proof.* The forgetful functor $U : \text{SmGrp} \longrightarrow \text{Set}$ from the category of semigroups and semigroups maps to the category of sets is right adjoint. Thus, given the fact that a semigroup operation can be defined pointwise on the pullback of the sets underlying the semigroups, pullbacks can be computed in Set. In the one direction, it is straightforward to see that for $x, y \in R$, $xy \in I$ if and only if $x \in I$ or $y \in I$. In the other direction, if $I$ is a prime ideal, then $I \times R$ must be a retract of $I \times_R (R \times R)$. This follows because the projection $I \times_R (R \times R) \longrightarrow R \times R$ must have an image isomorphic to $(I \times R) \cup (R \times I)$ since otherwise we will have $xy \in I$ with neither $x$ nor $y$ being an element of $I$. The same argument is made for $R \times I$. But by properties of pullbacks, the projection $I \times_R (R \times R) \longrightarrow R \times R$ is a monomorphism. $\square$

---

[1] In certain instances of exact categories, such as those of quasicoherent sheaves on concentrated schemes, and in the context of the main result of [viii], we have to restrict to *flat* morphisms of schemes with the property that under their inverse image functors, the preimage of a localizing Serre subcategory of finite type is again a localizing Serre subcategory of finite type.

The situation in TriCat$^\otimes$ almost allows us to adapt this formulation of prime ideals to prime thick tensor subcategories. That is, the semigroup structure is "almost" reflected by the symmetric monoidal structure. Here "almost" alludes to the observation that the symmetric monoidal structure comes with a unit while the semigroup structure does not. However, it turns out that once we change our setup to the category whose objects are the same as those of TriCat$^\otimes$, and whose morphisms resemble triangulated symmetric monoidal functors but with the unit preservation requirement removed, we can devise a notion of primeness which is a generalization of 0.1.0: albeit which captures primes in the context of derived categories of perfect complexes of qcqs[1] schemes, under certain strong impositions. Indeed, if we let TriCat$^{\otimes'}$ be the afore-described category, then we can give the following characterization of a prime ideal in TriCat$^{\otimes'}$.

**Definition 0.1.1.** *Let $\mathcal{I}$ be a thick subcategory of $\mathcal{K} \in$ TriCat$^{\otimes'}$ and assume the pullback of the inclusion $\mathcal{I} \hookrightarrow \mathcal{K}$ along $\otimes : \mathcal{K} \times \mathcal{K} \longrightarrow \mathcal{K}$ exists in TriCat$^{\otimes'}$. Then we say $\mathcal{I}$ is prime if and only if the pullback of inclusion along the tensor product in TriCat$^{\otimes'}$ is the full subcategory of $\mathcal{K} \times \mathcal{K}$ spanned by objects of both $\mathcal{K} \times \mathcal{I}$ and $\mathcal{I} \times \mathcal{K}$.*

When we restrict ourselves to the full subcategory of TriCat$^{\otimes'}$ containing objects of the kind $D_{\text{qc}}(X)$ where $X$ is a qcqs scheme—that is, the derived categories of unbounded complexes of $\mathcal{O}_X$-modules with quasicoherent cohomology—then given an open immersion $j : U \hookrightarrow X$ where $U \subseteq X$ is quasicompact open, the right derived functor $j_* : D_{\text{qc}}(U) \longrightarrow D_{\text{qc}}(X)$ obeys the projection formula (cf. [ix] 3.9.4); this is to say, $j_*E \otimes G \simeq j_*(E \otimes j^*G)$ given $E \in D_{\text{qc}}(U)$ and $G \in D_{\text{qc}}(X)$. It is then immediate that given any $A \in D_{\text{qc}}(X)$ and $B \in j_*(D_{\text{qc}}(U))$ we have that $A \otimes B$ is in the essential image $j_*(D_{\text{qc}}(U))$. Furthermore, since the counit $j^*j_* \longrightarrow 1$ is an isomorphism, we see that $j_*(E \otimes F) \simeq j_*(E \otimes j^*j_*F) \simeq j_*E \otimes j_*F$. Consequently, after taking into account the observation that triangulated categories admit finitary biproducts, one deduces that $j_*(D_{\text{qc}}(U))$ is equivalent to a thick tensor subcategory of $D_{\text{qc}}(X)$ in TriCat$^{\otimes'}$. A similar conclusion follows when instead of open immersions like $j$, one considers the closed immersion $i : X \backslash U \longrightarrow X$. In fact, where $i$ has *finite Tor-dimension*, $i_* : D_{\text{qc}}(X \backslash U) \longrightarrow D_{\text{qc}}(X)$ preserves compact objects (cf. [xxvii] 4.4). Compact objects in turn coincide with perfect complexes when $X$ is qcqs (cf. [x] 75.16.1) and hence our conclusion descends to the setting of tensor triangulated geometry, $D_{\text{perf}}(X \backslash U) \longrightarrow D_{\text{perf}}(X)$; this being also true when $j$ is a *quasi-perfect* map (cf. [ix] 3.22). And in full generality (when $X$ is concentrated), we observe that prime thick tensor subcategories are found as kernels of left adjoints of such maps (cf. [xv] 4.1).

These conclusions, notwithstanding the strong impositions on $X$, point toward a description of thick tensor subcategories—and by extension to a description of the prime thick tensor subcategories through 0.1.1—that is *concretely* functorial; which is to say, a description as certain kinds of morphisms in TriCat$^{\otimes'}$ (from without the objects of TriCat$^{\otimes'}$) in contrast to a description that is "internal" to the objects of TriCat$^{\otimes'}$. And fleshing out this newfound understanding is an accompanying challenge: though we have been able to give a functorial illustration of primeness in TriCat$^{\otimes'}$, how are we to arrive at the Balmer spectrum only from this? In other words, how does the topology on the set[1] of these primes follow?

---
[1] Quasicompact quasiseparated/concentrated.

In the appendix, on occasion, we will limit ourselves to concentrated schemes $X$ which have $U = X \setminus \overline{\{x\}}$ quasicompact and $j : U \hookrightarrow X$ quasi-perfect whenever $x \in X$ is a closed point. This allows us to show that with thick tensor subcategories understood as certain kinds of functors, and especially prime thick tensor subcategories understood as in 0.1.1, topos theoretic methods enable us to recover the Balmer spectrum[2] (cf. Theorem A.3.4). Furthermore, here, the "point" of tensor triangulated geometry is united with its topos theoretic counterpart; bringing home a notion which until this deduction, seems a far road gone from its day-to-day (topological) occurrence.

The use of topos theoretic methods above is preceded by identifying objects in $\mathsf{TriCat}^{\otimes'}$ with pre-ordered sets (prosets) having Grothendieck topologies (prosites) and identifying the morphisms in $\mathsf{TriCat}^{\otimes'}$ with maps of sites. It is the spaces of points of the Grothendieck topoi resulting from these prosites which enable us to recover the space underlying the Balmer spectrum. In this endeavor, we will utilize the work of O. Caramello in [xi]. It is from this background that the article sets off.

We will begin by introducing $\infty$-*prosets* which will facilitate the usage of results, as well as techniques, from the theory of prosets (above) in the $\infty$-categorical world. Briefly put, an $\infty$-proset is a map of simplicial sets[3] $X : \mathsf{N}(\Delta)^{\mathrm{op}} \longrightarrow \mathsf{Cat}_\infty$ which satisfies the Segal and completeness conditions and whose essential image has prosets as objects. Equivalently, it is a simplicial proset $\Delta^{\mathrm{op}} \longrightarrow \mathsf{Proset}$ of a certain kind which when composed with the forgetful functor $U : \mathsf{Proset} \longrightarrow \mathsf{Set}$ results in an $\infty$-category. Informing this move, the $\infty$-categorical setting affords us numerous advantages if we are to insist on the line of thinking explored in the previous paragraphs. Most noticeably, the $\infty$-categorical version of $\mathsf{TriCat}^{\otimes'}$ is much better behaved under limits and colimits; for example, in the $\infty$-categorical world the limit referenced in 0.1.1 necessarily exists in the appropriate sense needed to sharpen the definition i.e as a homotopy pullback. More so, as we shall see below, when we wish to make a "big picture" account of the results obtained via topos theoretic methods regarding the aforementioned spaces, we find the $\infty$-categorical setting all the more conducive.

### What is a Space?

What we see exemplified in the appendix is a general procedure which produces for any map $D : \mathcal{C}^{\mathrm{op}} \longrightarrow \mathsf{Cat}^\otimes$, a functor $D(\mathcal{C}^{\mathrm{op}})^{\mathrm{op}} \longrightarrow \mathsf{Top}$ (cf. Remark A.2.7). Here $\mathsf{Cat}^\otimes$ is the category of (essentially) small symmetric monoidal categories and symmetric monoidal functors and Top is the category of topological spaces and continuous maps. This is possible whenever we can come up with a systematic way to identify objects of $\mathsf{Cat}^\otimes$ with prosites and maps in $(\mathsf{Cat}^\otimes)^{\mathrm{op}}$ with maps of sites. Furthermore, in the case where $D$ has an essential image rich enough to admit certain notions of localization[4], the functor above upgrades to a functor $D(\mathcal{C}^{\mathrm{op}})^{\mathrm{op}} \longrightarrow \underline{\mathsf{Stk}}_\mathbb{Z}$; where $\underline{\mathsf{Stk}}_\mathbb{Z}$ is the category of stacks over $\mathbb{Z}$.

In our case, this systematic process is achieved owing to certain "combinatorial" properties available to the functor $D_{\mathrm{qc}} : \mathsf{Sch}_{\mathrm{conc}}^{\mathrm{op}} \longrightarrow \mathsf{TriCat}^\otimes$ on the category of concentrated schemes and separated morphisms of finite type. Namely, for any $f : Y \longrightarrow X$ in $\mathsf{Sch}_{\mathrm{conc}}$ we have an adjunction $(f^* \dashv f_*) : D(Y) \longrightarrow D(X)$ in $\mathrm{TriCat}^{\otimes'}$ which comes furnished with the projection formula $f_*(x \otimes f^*y) \simeq f_*x \otimes y$, where $f^* = D_{\mathrm{qc}}(f)$. We also need the specification of coverages $\{f_i : U_i \longrightarrow X\}_{i \in I}$ for each $X \in \mathrm{Sch}_{\mathrm{conc}}$ which interact with both $f_*$ and $f^*$ in a determined way[1]. The resulting functor $D_{\mathrm{qc}}(\mathrm{Sch}^{\mathrm{op}}_{\mathrm{conc}})^{\mathrm{op}} \longrightarrow \underline{\mathrm{Stk}}_{\mathbb{Z}}$, as has been alluded to earlier, is the Balmer spectrum[2] (cf. Theorem A.1.11). And as a matter of fact, it is seen that with slight variations of the procedure giving rise to the latter map, one obtains divers functors into $\underline{\mathrm{Stk}}_{\mathbb{Z}}$ from $D_{\mathrm{qc}}(\mathrm{Sch}^{\mathrm{op}}_{\mathrm{conc}})^{\mathrm{op}}$ (cf. Remark A.2.7). However, by a result in its nascent form attributable to P. Balmer, the Balmer spectrum stands out among these functors since it is final in a certain category of the aforementioned functors (cf. [vi] 3.2).

---

[1] The objects of $\mathsf{TriCat}^{\otimes'}$ are essentially small and hence the collection of these primes up to equivalence is a set.
[2] In fact, if we are willing to relax this understanding of thick tensor subcategories slightly, we are able to recover the Balmer spectrum for any tensor triangulated category. This is shown in A.1.11.
[3] The target is the $\infty$-category of *small* $\infty$-categories.
[4] For example Verdier localization and Serre localizations.

All in all, these results intimate the viewpoint that a space is not merely a category as is suggested by the dominant gist of noncommutative algebraic geometry. Rather, it is a pair constituted by a category and "a structure"—a structure that is determined relative to all other categories in some class of categories and maps therein. Succinctly expressed, it is a functor into $\mathrm{Cat}^{\otimes}$ satisfying certain conditions. And at bare minimum, it is a class of pairs of adjoint functors between pairs of objects in a class of symmetric monoidal categories which is determined from the image of a contravariant pseudofunctor into $\mathrm{Cat}^{\otimes}$: where for each functor in its image, a right adjoint is demanded so that together the two functors observe a particular desired relationship relative to the symmetric monoidal structure. For now, we may (conventionally) identify such a pseudofunctor as a *four functor formalism/sheaf theory* satisfying extra properties (cf. [xii] 2.1).

Moreover, following from what we have seen in the previous paragraphs, we can show that given the desired kind of four functor formalism, we obtain a triad of functors which "wants" to be what we call an *inverse Tannakian formalism*. This is a triple of functors $(D, T_D, T'_D)$, where $D$ is a 4-functor formalism and $S$ some scheme, such that the following diagram commutes up to natural isomorphism

$$\begin{array}{ccc} \mathcal{C}^{\mathrm{op}} & \xrightarrow{D} & \mathrm{Cat}^{\otimes} \\ {\scriptstyle T'_D} \downarrow & \nearrow {\scriptstyle T_D} & \\ \underline{\mathrm{Stk}}^{\mathrm{op}}_S & & \end{array}$$

and the functor $T'_D$ is the composition of $D$ with a functor $D(\mathcal{C}^{\mathrm{op}}) \longrightarrow \underline{\mathrm{Stk}}^{\mathrm{op}}_S$ which is a section of $T_D$. There is a variety of concrete examples which inspire this definition. We give a few.

**Example 0.2.0.** Let $D$ be the functor $\mathrm{QC} : \mathrm{Sch}^{\mathrm{op}}_{\mathrm{geom}} \longrightarrow \mathrm{Cat}^{\otimes}$ from the category of quasi-compact separated schemes and morphisms of schemes, mapping each such scheme to its category of quasicoherent sheaves and morphisms of schemes to their inverse image functors. If we set $T'_{\mathrm{QC}} = \mathrm{Spec}^{\mathrm{op}} \circ \mathrm{QC}$ and $T_{\mathrm{QC}} = \mathrm{QC}$, then we obtain an inverse Tannakian formalism. $\mathrm{Spec}^{\mathrm{op}} : \mathrm{QC}(\mathrm{Sch}^{\mathrm{op}}_{\mathrm{conc}}) \longrightarrow \underline{\mathrm{Stk}}^{\mathrm{op}}_{\mathbb{Z}}$ is the functor $\mathrm{QC}(X) \longmapsto \mathrm{Hom}^{\otimes}(\mathrm{QC}(X), -)$; the image of this functor acts on an affine scheme $A$ by mapping it to the groupoid of all cocontinuous symmetric monoidal functors $\mathrm{QC}(X) \longrightarrow \mathrm{QC}(A)$. This is a restriction of the QC functor which in its fullness gives rise to *geometric Tannaka duality* as presented in [xiii] Theorem 5.11. When QC is considered without restriction, it acts on geometric stacks and the essential image expands in such a manner as to lie outside the definition we seek. Specifically, the essential image is equivalent to the 2-category with tame, complete and symmetric monoidal abelian categories as objects and the groupoid of symmetric monoidal functors preserving flat objects as hom-objects. It is the problem mentioned under Remark 5.12 in [xiii] that inspires the "inverse" in *inverse Tannakian formalism*. In fact, a general construction that gives a geometric stack for each object in some ambient category containing the essential image of QC would, by definition, give a solution for the problem whenever a subcategory is identified where this construction is a section of $T_{\text{QC}}$.

---

[1] Among these coverages being those maps that describe primes in the target category as per the aforementioned functorial formulation.

[2] The category LRS' can be mapped to $\underline{\mathrm{Stk}}_{\mathbb{Z}}$ via composition with the Yoneda embedding and sheafification.

**Example 0.2.1.** For an $\infty$-category $\mathcal{E}$ of nice enough algebraic stacks over a field $k$ (cf. [xiv] 5.12 and 2.3), the functors described as follows give rise to an inverse Tannakian formalism when arranged appropriately: $\mathcal{D}_{\text{qcoh}} : \mathcal{E}^{\text{op}} \longrightarrow \widehat{\mathsf{Cat}}_{\infty}^{\otimes}$ where the target is the $\infty$-category of large symmetric monoidal categories and symmetric monoidal functors, the functor $\mathcal{F} : \mathcal{D}_{\text{qcoh}}(\mathcal{E}^{\text{op}}) \longrightarrow \mathsf{N}(\underline{\mathsf{Stk}}_k)^{\text{op}}$ given pointwise as $\mathcal{D}_{\text{qcoh}}(\mathcal{X}) \longmapsto \mathcal{F}_{\mathcal{X}}$ where $\mathcal{F}_{\mathcal{X}}$ is the $\infty$-stack on the etale site of affine $k$-schemes which maps an affine $k$-scheme $S$ to $\text{Map}_k^{\otimes}(\mathcal{D}_{\text{qcoh}}(\mathcal{X}), \mathcal{D}_{\text{qcoh}}(S))$[1] and, $\mathcal{D}_{\text{qcoh}} : \mathsf{N}(\underline{\mathsf{Stk}}_k)^{\text{op}} \longrightarrow \widehat{\mathsf{Cat}}_{\infty}^{\otimes}$. This setup presents a version of *derived Tannaka duality*.

**Example 0.2.2.** Let $\mathcal{C}$ be the category of finite dimensional locally compact Hausdorff spaces and let $D : \mathcal{C}^{op} \longrightarrow \text{TriCat}^{\otimes}$ be the 4-functor formalism mapping each object $X \in \mathcal{C}$ to the derived category of its abelian sheaves. Let $(\bullet)^{\text{Sch}} : \mathcal{C} \longrightarrow \underline{\mathsf{Stk}}_{\mathbb{Z}}$ be the functor $X \longmapsto \underline{X}^{\text{Sch}}$, where $\underline{X}^{\text{Sch}}$ maps each affine scheme $S$ to the discrete groupoid $C^0(|S|, X)$. These two functors together with $D_{\text{qcoh}} : \underline{\mathsf{Stk}}_{\mathbb{Z}}^{\text{op}} \longrightarrow \text{TriCat}^{\otimes}$ present a cartoon imitating an inverse Tannakian formalism; i.e. $D(X) \simeq D_{\text{qcoh}}(\underline{X}^{\text{Sch}})$ (cf. [iii] 1.7).

This richness, insofar as capturing diverse notions of space is concerned, gives reason to consider the concept of *inverse Tannakian formalisms* a common thrust among different approaches to noncommutative algebraic geometry. Actually, a successful endeavor for this article consists in part having, to an extent, brought together through this notion two currents of noncommutative algebraic geometry. The *Tannakian type*[2] current which considers particular subcategories of commutative monoid objects of the $\infty$-category of stable $\infty$-categories and exact functors as stacks of anima on subcanonical sites of slices of $\mathsf{N}(\text{CRing})^{\text{op}}$ (as exemplified by 0.2.0 and 0.2.1), and the current of *tensor triangulated geometry* discussed before.

There are two obstructions in the way of our desired 4-functor formalisms giving rise to inverse Tannakian formalisms. In the case of $D_{\text{qc}} : \text{Sch}_{\text{conc}}^{\text{op}} \longrightarrow \text{TriCat}^{\otimes}$, even though we can construct $T'_{D_{\text{qc}}} : \text{Sch}_{\text{conc}}^{\text{op}} \longrightarrow \underline{\mathsf{Stk}}_{\mathbb{Z}}^{\text{op}}$, it is not always the case that $D_{\text{qc}}(\text{Sch}_{\text{conc}}^{\text{op}}) \longrightarrow \underline{\mathsf{Stk}}_{\mathbb{Z}}^{\text{op}}$ is a section of $D_{\text{qc}} : \underline{\mathsf{Stk}}_{\mathbb{Z}}^{\text{op}} \longrightarrow \text{TriCat}^{\otimes}$; which is the most natural choice for $T_{D_{\text{qc}}}$. This is the first impediment; the lack of a canonical way to obtain a $T_D$ that fits the sought schema given both $D$ and $T'_D$. In the event that we are able to overcome this obstacle, we are again faced with the possibility of a multitude of ways to extend $D$ through $T'_D$ i.e to construct $T_D$. In this case, we would like to know that there is a way of doing such extensions with the most "efficiency". This is the second impediment. Thence, we should hope for not only an upgrade of 4-functor formalisms to inverse Tannakian formalisms, but also an upgrade that satisfies a universal property. This seems to be the situation on occasion of certain restrictions of $D_{\text{qc}} : \text{Sch}_{\text{conc}}^{\text{op}} \longrightarrow \text{TriCat}^{\otimes}$ given the result implied by Balmer's distinquishing of Balmer spectra (cf. [vi] 3.2). However, the universality here is not as explicit as one would like, since one still has to consider the extensions (of $D_{\text{qc}}$ through $T'_{D_{\text{qc}}}$) in the context of an extraneous category; the category of *classifying support data*.

---

[1]The keen reader may get an in-depth introduction in [xiv] Section 2.2.

[2]The word "Tannakian" is chosen in the situation of this writing over the more historically natural "functor of points" to emphasize that these follow from generalizations of classical Tannaka duality.

A possible way to naturally enhance (the desired) 4-functor formalisms to inverse Tannakian formalisms is through Kan extensions. We would need only that a Kan extension of the inclusion $p : D(\mathcal{C}^{\mathrm{op}}) \hookrightarrow \mathrm{Cat}^{\otimes}$ along the constructed $\mu_D : D(\mathcal{C}^{\mathrm{op}}) \longrightarrow \underline{\mathsf{Stk}}_{\mathbb{Z}}^{\mathrm{op}}$ exist for such an enhancement to occur. For example, when the latter functor is fully faithful, both the left Kan extension ($\mathrm{Lan}_{\mu_D} p$) and right Kan extension ($\mathrm{Ran}_{\mu_D} p$) of $p$ necessarily exist and hence also the desired enhancement of $D$. On the other hand, in practice, the essential image of $D$ is usually a subcategory of $\mathrm{TriCat}^{\otimes}$. Thus, it is suitable that both $\mathrm{Ran}_{\mu_D} p$ and $\mathrm{Lan}_{\mu_D} p$ also have their essential images as subcategories of $\mathrm{TriCat}^{\otimes}$. However, when $\mu_D$ is fully faithful this possibility relies on the existence of (co)limits in $\mathrm{TriCat}^{\otimes}$; and this is already well known to be a slippery task when it is posed in the form of gluing derived categories of complexes of sheaves on schemes. Altogether, it is in keeping with these aforementioned challenges, and after the groundwork laid down by the introduction of $\infty$-prosets, that we will introduce *suprematic spaces* as artifacts which promise to resolve these difficulties.

And the deterrent that first gives away is the preceding one. In the $\infty$-categorical world where suprematic spaces live, $\mathrm{TriCat}^{\otimes}$ is replaced by the $\infty$-category of stable $\infty$-categories and exact functors and this is closed under small $\infty$-categorical limits and filtered colimits (cf. [xvi] 1.1.4.4 and 1.1.4.6 ). Furthermore, $\mathsf{Pr}^L$ and any $\infty$-topoi $\mathcal{X}$ also admit small $\infty$-categorical limits and colimits[1] (cf. [i] 5.5.3.13, 5.5.3.18 and 6.1.0.6); and we will work primarily from within either of these $\infty$-categories. This means that given any $F : \mathcal{C} \longrightarrow \mathsf{Pr}^L$, left Kan extensions of $F$ along inclusions of full $\infty$-subcategories $\mathcal{C} \hookrightarrow \mathcal{C}'$ always exist (cf. [i] 4.3.2.2 and 4.3.2.6). This turns out to be not too strong an imposition on the $\infty$-categories we will be working with. Therefore, a suprematic space $\mathcal{O}^{\mathrm{op}} : \mathcal{T} \longrightarrow \mathcal{X}^L$, coming with the datum $\mathcal{X}^L \subseteq \mathsf{Pr}^L$, is found easily enhanced to (the "derived" version of) an inverse Tannakian formalism as soon as the question of this enhancement is reduced to one pertaining to the existence of Kan extensions.

Indeed, the utility of suprematic spaces in producing inverse Tannakian formalisms is demonstrated by Theorem A. What remains to be shown is the "efficacy" of this production. Given one suprematic space, this question becomes fairly straightforward to answer, seeing that Kan extensions are put to use; these naturally enjoy a universal property. The more strenuous exercise is to determine an efficacious way to carry the factorization over a collection of suprematic spaces that essentially "look the same" from the viewpoint of animated stacks; that is, suprematic spaces that have the same *geometric content* (cf. Definition 2.2.5). And such a factorization follows when we weaken how strictly the triangle defining inverse Tannakian formalisms commutes: at this point, one may take the definition of the latter to be a lifting of the initial definition via the nerve functor $\mathsf{N} : \mathsf{Cat} \longrightarrow \mathsf{Cat}_{\infty}$. If we elect that the triangle commutes only up to natural transformation, then we will be able to show that among all the extensions, which are "coherent enough" with our needs, of suprematic spaces with the same geometric content into inverse Tannakian formalisms, there exists an initial one. This will follow from a theorem which "sits" between Theorem B and Theorem C; it is sufficient for the latter given the former and is itself implied by a mild version of the former (cf. Theorem 3.2.9). This theorem identifies each object in the image of $f$ in Theorem B with a $\mathcal{G}$-structure on $\mathcal{X}^L$, where $\mathcal{G}$ is the geometry enveloping a pregeometry on $\mathsf{Stk}_{S}^{\mathrm{op}}$ determined by $f$; this determination being in part the content of the theorem. The $\infty$-category of the structures of this pregeometry on $\mathcal{X}^L$ then happens to be equivalent with the $\infty$-category of those extensions "coherent enough" with our needs, given the suprematic space identifying the $\mathcal{G}$-structure in question. Now, by [iv] Proposition 3.4.5, this $\infty$-category contains an initial object.

---

[1]Henceforth, unless where there is danger of ambiguity, $\infty$-categorical (co/filtered) limits and (co/filtered) colimits will simply be referred to as (co/filtered) limits and (co/filtered) colimits.

Additionally, as testified by Theorem B, suprematic spaces correspond to 6-functor formalisms. We will show that in this interaction, there is naturally a confluence of the protagonists that gives exactly the kind of 4-functor formalisms we need to replicate in the $\infty$-categorical world, the project initiated in the appendix: that is, to obtain in this new setting the ability to extract through topos theoretic methods a contravariant map from a certain subcategory of $\mathsf{Pr}^L$ to $\mathsf{Stk}_S$. There is, to begin with, an innate tension between the two sides of this correspondence. Structured spaces can be thought of as $\infty$-algebras over an essentially $\infty$-algebraic theory and, 6-functor formalisms as pointers towards a geometric world that underpins the algebraic structure inherent to them (cf. [xvii] Structured $\infty$-topos 2.4 and [i] 6.2.3.20). A similar (eponymous) tension sits at the heart of noncommutative algebraic geometry: that between "algebra/quantity" and "geometric spaces". However, that there is an alignment between the two sides is the surprising fact that we (always) wish to exploit. For example, we will see that the datum of structured spaces encoding covers will allow us to pass from $\infty$-prosets to sites of $\infty$-prosets ($\infty$-prosites); at the same time, the behavior of these covers under functors of $\infty$-categories, as encoded by 6-functor formalisms, will allow us to capture maps of both $\infty$-prosets and $\infty$-prosites. And so from such and in this manner, the desired construction of animated $S$-stacks unfolds.

**Remark 0.2.3.** The map sought above is conceptually in the same neighborhood as the *smashing spectrum* functor of condensed mathematics. In the latter, one considers (large) frames of *smashing colocalizations* of cocomplete symmetric monoidal stable $\infty$-categories obtained via the identification of the latter with *coidempotent objects* of said $\infty$-categories (cf. [xviii] 2.5 and 3.17). In our context, as we shall observe in Remark 2.1.22, the information involving coverages, and which, as mentioned, is crucial to our construction, possibly entails smashing colocalizations once subjected to six functor formalisms. This is to say, at an appropriate specialization, the two constructions are coincident; in fact, we will see, when the images of the covering maps are fully faithful functors of $\infty$-categories, they necessarily produce smashing colocalizations when properties of 6-functor formalisms are brought fully to bear (cf. Remark 2.1.23).

### What is a Motive?

Six functor formalisms sprung out of Grothendieck's work in étale cohomology. For a morphism of schemes $f: X \longrightarrow Y$ (usually separated of finite type) there are relations which arise between their étale cohomologies. This is also the case for select other cohomology theories; such as when one considers finite dimensional locally compact Hausdorff spaces and their sheaf cohomologies (cf. [iii] Lecture I). These relations, it was observed, were formal consequences of a handful of occurences which involve six functors $(f^*, f_*, f_!, f^!, \otimes, \underline{\mathrm{Hom}})$ and coherence conditions which dictate how they relate to each other. Collectively, these relations have come to be known as six functor formalisms. Classically, a six functor formalism takes the shape of a functor $\mathrm{Sch}^{\mathrm{op}}_{/S} \longrightarrow \mathrm{TriCat}^{\otimes}$ which comes with extra properties which we will specify in the second section of this article. And taken at face value, this functor is a sheaf theory (cf. [xii] 2.1).

A Weil cohomology theory over an algebraically closed field $k$ is a symmetric monoidal functor $H^*: \mathrm{SmProj}/k^{\mathrm{op}} \longrightarrow \mathrm{Gr}^{\geq 0}\mathrm{Vec}_K$ which satisfies some extra properties. Here, the domain is the category of smooth projective varieties over $k$ whereas the codomain is the tensor category of graded finite dimensional $K$-vector spaces ($K$ is a field of characteristic 0) and graded linear maps. In a quest to explain analogous phenomena present across the vast array of these cohomologies, Grothendieck envisioned the theory of motives. Concretely put, the theory suggests a category (of motives), $\mathcal{M}ot$, which is abelian semisimple, symmetric monoidal, Tannakian over $K$, enriched over

$\mathrm{Vec}_K$ and most importantly, which extends any Weil cohomology as follows (cf. [xix] 4.6):

$$\begin{CD}
\mathrm{SmProj}/k^{\mathrm{op}} @>H^*>> \mathrm{Gr}^{\geq 0}\mathrm{Vec}_K \\
@VhVV @AA\tau_h^*A \\
\mathcal{M}ot
\end{CD}$$

There is, of course, a wealth of remarkable reasons why the category $\mathcal{M}ot$ is conceived in the way that it is, and why it carries immense value for the field of algebraic geometry. We find these reasons to be beyond the motivation of this article and thus to be cost-effective in the economy of this article, we will altogether ignore them. That said, we should find it worthwhile to observe that some Weil cohomologies can be formulated as sheaf theories. This has been shown to be the case for $\ell$-adic cohomology, algebraic de Rham cohomology, and Betti cohomology (cf. [xii] 2.2).

It is a curious fact that once formulated as sheaf theories, these sheaf theories also happen to be six functor formalisms. It is, however, not exactly clear that all Weil cohomology theories can be understood this way. In fact, it takes considerable technical effort to show that one is a sheaf theory. However, what is clear is that some of the comparisons that prompt motives can still be accessed at the level of sheaf theories. For example, in the case of complex algebraic varieties, the equivalence between de Rham cohomology and Betti cohomology—established by Grothendieck—is expressed as an equivalence between their respective sheaf theories (cf. [xx] Section I). If we are taken to task, and should we find ourselves contemplating the reality of enhancement to sheaf theories for the full expanse of Weil cohomologies, and taking the position that such enhancements are always six functor formalisms, then we should find presented to us this question (which somewhat turns motivic aspiration on its head): what does the register of motives look like in this new setting?

There is already a response to this call. But first, we need to transition from looking at $\mathrm{SmProj}/k$ to looking at $\mathrm{Sch}/B$ (schemes of finite type over a noetherian scheme $B$ of finite Krull dimension) and, to transmute ourselves over to the world of cocomplete symmetric monoidal stable $\infty$-categories and cocontinuous exact symmetric monoidal functors. In [xx], a *coefficient system* is introduced as a functor[1] $C : \mathrm{Sch}_B^{\mathrm{op}} \longrightarrow \mathsf{CAlg}(\mathsf{Cat}_\infty^{\mathsf{st,ex}})$ satisfying extra axioms needed to define a 6-functor formalism (some are a priori satisfied). Notice that restricting oneself to functors of such kind that map to $\mathsf{Pr}^L$, one automatically obtains sheaf theories. The domain for a coefficient system, other than being the category of finite type $B$-schemes, is also equipped with the Cartesian symmetric monoidal structure (cf. [xx] 7.2). Hereafter, a subcategory $\mathsf{CoSy}_B^{\mathsf{c}} \subseteq \mathsf{Fun}(\mathsf{Sch}_B^{\mathrm{op}}, \mathsf{CAlg}(\mathsf{Cat}_\infty^{\mathsf{st,ex}}))$ is isolated: which by definition is not full. Then it is shown in [xx] Theorem 7.14 that this $\infty$-category admits an initial object. This coefficient system is none other than the association $X \longmapsto \mathcal{SH}(X)$ where the right-hand side is Morel-Voevodsky $\mathbb{A}^1$-homotopy theory.

Generally (and loosely) speaking, Theorem C offers a different approach to the question posed. Where the latter result offers a universal map only among maps of a certain kind into $\mathsf{CAlg}(\mathsf{Cat}_\infty^{\mathsf{st,ex}})$, and by restriction into $\mathsf{Pr}^L$, Theorem C offers a 6-functor formalism through which certain 6-functor formalisms factor through. That is, one is closer in spirit to the motivic project than the other. Furthermore, there is a sense in which the construction of $L_{\mathrm{\acute{e}t}}(\mathsf{Stk}_{S_|})$ can be likened with the formulation of $\mathcal{SH}(X)$. That is, it can be said that the universal property satisfied—pointwise—by the latter, as stated in Theorem 5.2 of [xxix], is of the same ilk as the one employed to arrive at Theorem C (cf. Remark 3.2.20). However, whereas we are able to present somewhat of a refinement of B. Drew and M. Gallauer's result (relative to the question posed here), we fall short of meeting its standard in that, as yet, we are not able to guarantee the universality of this factorization. At any rate, the initial formulation of motives does not adhere to the proper functorial sense of universality, and, therefore, we are still keeping in the ethos of their original architecture.

---

[1] Here we ignore set theoretic size specifications needed to formulate the target.

The justification for Theorem C occurs in two parts. In the first, one notices that the maps in the essential image of $f$ correspond to pregeometry structures on $\mathcal{X}^L$ and, therefore, admit a universal pregeometry structure (cf. [iv] 1.4.2, 3.4.3 and 3.4.12). In the second, one begins with $p : \mathsf{Fun}(\Delta^1, \mathcal{X}) \longrightarrow \mathcal{X}$ (where $\mathcal{X}$ is any (*large*) $\infty$-topos), which is evaluation at the endpoint $\{1\}$. This is a Cartesian fibration (cf. [i] 6.1.1.1). Therefore, through straightening/unstraightening and properties of $\infty$-topoi, it is classified by a limit preserving functor $\mathcal{X}^{\mathrm{op}} \longrightarrow \mathsf{Pr}^L$ (cf. [i] 6.1.3.9). From this, one can construct a "universal self-dual" category of $\mathcal{X}^{\mathrm{op}}$ keeping in mind Remark 6.3.5.10 in [i]. In a handful of more steps, one arrives at Theorem C by recalling that given geometric envelopes, the universal structured space takes values in the $\infty$-topos determined by sheaves of anima on the site fixed by the Grothendieck topology carried by the *geometry* (cf. [iv] 1.4.2 and 3.4.3).

In the last section of this article, we will allocate our energies to the explication of these claims. In the meantime, one observes, it is altogether pleasing that this theorem comes, as a reward of sorts, at the resolution of ostensibly innocuous questions: What is a point? What is a space?

In terms of organization, this paper is given in three bullet points (excluding the appendix). In the first, we introduce $\infty$-prosets and $\infty$-prosites and adapt some fundamental results of pointless topology to the setting. In the second, our main goal will be to prove both Theorem A and Theorem B. We will introduce suprematic spaces and in the process, contextualize the results in the appendix. This will pave the way for the third act, where we will obtain results which culminate in a proof of Theorem C (as mentioned previously).

**Acknowledgements.** The author thanks P. Balmer for correcting an error originally made regarding the historicity of their work. The author thanks Julia R. Gonzalez for having agreed to read a draft of the appendix of this paper and on occasions pointed out errors and suggested more precise and compact language. The author also thanks M. Gallauer for pointing out an initial error in the formulation of Proposition 0.1.0. Lastly, the author would like to give immense gratitude to his family (foremost his parents) for the unspeakable amount of support accorded whilst writing this article; and up to the point of writing it. This article is dedicated to them.

# The Yoga of $\infty$-Prosets

We will begin by collecting theorems from pointless topology and the study of ordinary prosets (and prosites) that we find useful. We will then introduce $\infty$-prosets and $\infty$-prosites and prove results about them (relevant to us) which mimic those of pointless topology. Of special interest will be $\infty$-prosites induced by Grothendieck topologies of a certain make-up; which we will refer to as *packeted* (resp. *retro-packeted*) $\infty$-prosites. These $\infty$-prosets behave reasonably under classical Stone-type dualities; particularly those involving *spatial locales*(resp. *coherent locales*). The main theorem we wish to prove concerns the establishment of functors from $\infty$-prosites to animated $S$-stacks. We will provide two such functors. Overall, familiarity with pointless topology is assumed:

we will only reintroduce terms we find pivotal to the work intended to be done.

## Recollections

We give a primer on the state-of-the-art as construed according to the aims of this article.

Let $(\mathcal{C}, J)$ be a small[1] site.

- A *J-ideal* on $\mathcal{C}$ is a subset $I \subseteq obj(\mathcal{C})$ such that for any morphism $f : b \longrightarrow a$ in $\mathcal{C}$ if $a \in I$ then $b \in I$, and for any $J$-covering sieve $R$ on an object $c$ of $\mathcal{C}$, if $dom(f) \in I$ for every $f \in R$ then $c \in I$. For each object of $c \in \mathcal{C}$, we have the smallest $J$-ideal containing $c$ which we label as $(c) \downarrow_J$: $J$-ideals of this kind are called *principal J-ideals*. Note that when $\mathcal{C}$ is a proset, $(c) \downarrow = \{d \in \mathcal{C} : d \leq c\}$ corresponds to $(c) \downarrow_J$ when $J$ is subcanonical.

- We denote by $Id_J(\mathcal{C})$ the set of $J$-ideals on $\mathcal{C}$ endowed with the subset-inclusion order relation. It is worth mentioning that an arbitrary union of $J$-ideals is not necessarily a $J$-ideal. For example, the union of two $J$-ideals may contain $dom(f)$ for every $f$ in a covering $R$ of $c$, but not contain $c$ itself. We need instead to consider a modification of taking unions to obtain an operation on $Id_J(\mathcal{C})$; taking the smallest $J$-ideal containing the union. On the other hand, arbitrary intersections of $J$-ideals is again a $J$-ideal. Given these two operations, as join and meet respectively, $Id_J(\mathcal{C})$ is a frame. This is implicit in Theorem 1.1.2.

- Let $\mathcal{C} = (\mathcal{K}, \leq)$ be a proset. A *J-prime filter* on $\mathcal{C}$ is a subset $\mathfrak{F} \subseteq ob(\mathcal{C})$ such that:
    1. $\mathfrak{F}$ is non-empty.
    2. $a \in \mathfrak{F}$ implies $b \in \mathfrak{F}$ whenever $a \leq b$.
    3. For any $a, b \in \mathfrak{F}$ there exists $c \in \mathfrak{F}$ such that $c \leq a$ and $c \leq b$.
    4. For any $J$-covering sieve $\{a_i \longrightarrow a\}_{i \in I}$ if $a \in \mathfrak{F}$ then there exists $i \in I$ such that $a_i \in \mathfrak{F}$.

- We are given in [xi] 3.1 that a functor $(C, J_C) \longrightarrow (D, J_D)$ is a morphism of sites if it is exactly a *flat* functor which is cover preserving. The latter means that given a covering sieve $R \in J_C$, then $F(R) = \{F(f) : f \in R\}$ generates a covering sieve $R' \in J_D$. By Proposition 3.2 of [xi], the former, when dealing with prosites, is a proset functor $F : (\mathcal{K}, \leq) \longrightarrow (\mathcal{H}, \leq')$ such that:
    1. For each object $h$ of $\mathcal{H}$ there exists and object $k$ of $\mathcal{K}$ such that $h \leq' F(k)$.
    2. For any object $h$ of $\mathcal{H}$ having $h \leq' F(c)$ and $h \leq' F(c')$ there exists an object $c''$ of $\mathcal{K}$ such that $h \leq' F(c'')$, $c'' \leq c$ and $c'' \leq c'$.

    In general, a *flat* functor is a filtered colimit of representable functors. We will call a prosite functor that is also a morphism of sites a *prosite map*. Notice that whenever involved prosites are finitely complete, the latter is coincident with a left exact functor that is cover preserving.

- A morphism $f : \mathscr{X} \longrightarrow \mathscr{Y}$ of topoi is called a *geometric morphism*. It is a pair of adjoint functors $f^* \dashv f_* : \mathscr{X} \longrightarrow \mathscr{Y}$ where $f^*$ commutes with finite limits (left exact). $f^*$ is referred to as the inverse image functor and $f_*$ as the direct image functor. An example of a geometric morphism is the localization $\mathbf{Sh}(\mathcal{C}, J) \hookrightarrow [\mathcal{C}^{op}, \mathrm{Set}]$ where the inverse image functor is the sheafification functor.

---

[1]That is, one whose category has a set for its objects and a set for its hom-objects.

- Given any topos $\mathscr{X}$, a *point* of $\mathscr{X}$ is a geometric morphism Set $\longrightarrow \mathscr{X}$. The eponymic example is the geometric morphism $\varinjlim_{U \ni x}(\bullet) \dashv Sky_x$ : Set $\longrightarrow \mathbf{Sh}(X)$ between Set and the category of sheaves on a topological space $X$. The direct image functor is the functor mapping a set $S$ to the sheaf defined as mapping an open subset of $X$ to $S$ if it contains a particular point $x$ of $X$ and to the empty set otherwise; this is usually called the *skyscraper sheaf*. Its inverse image functor takes the stalk of a sheaf at the point $x$. The collection of all points on $\mathscr{X}$, when indexed by some set, can be given a topology that gives rise to a *space of points*. Moreover, this construction is functorial. In particular, when the category of points of a topos $\mathscr{X}$ is small[1], the frame of subterminals (objects whose unique morphisms to the terminal object are monic) determines a topology on the space of points of $\mathscr{X}$. We will refer to this space simply as "the space of points of $\mathscr{X}$". Notably, for the topos of sheaves of sets on a sober topological space, the former space recovers the original space. The keen reader may work through [xi] Section 2 for a detailed exposition.

- We state the following theorems due to O. Caramello. Their proofs can be found in [xi] Proposition 2.7 and Theorem 3.1 (for 1.1.1 and 1.1.2 respectively). Theorem 1.1.2 is the basis from which emanate the key results of this section.

**Theorem 1.1.1.** *Let $\mathcal{C}$ be a proset endowed with a Grothendieck topology $J$. Then the space of points of $\mathbf{Sh}(\mathcal{C}, J)$ is homeomorphic to the space which has: as its set of points the collection $\mathscr{F}_\mathcal{C}^J$ of the $J$-prime filters on $\mathcal{C}$, and as open subsets the sets of the form*

$$U_I = \{\mathfrak{F} \in \mathscr{F}_\mathcal{C}^J : \mathfrak{F} \cap I \neq \varnothing\}$$

*where $I$ ranges among the $J$-ideals on $\mathcal{C}$. In particular, a sub-basis for this topology is given by the subsets*

$$B_c = \{\mathfrak{F} \in \mathscr{F}_\mathcal{C}^J : c \in \mathfrak{F}\}$$

*where $c$ varies among the elements of $\mathcal{C}$.*

**Theorem 1.1.2.** *Let $\mathcal{C}$ be a proset, $J$ a Grothendieck topology on $\mathcal{C}$ and, $\mathbf{Sh}(Id_J(\mathcal{C}))$ the topos of sheaves of sets on the locale with $Id_J(\mathcal{C})$ as its frame of opens. Then $\mathbf{Sh}(\mathcal{C}, J)$ and $\mathbf{Sh}(Id_J(\mathcal{C}))$ are equivalent as topoi.*

- Recall that the category of locales and locale maps, Loc, is dual to the category of frames and frame maps, Frm. Therefore, a locale morphism has as an inverse image, a morphism of lattices which commutes with arbitrary joins. Additionally, Loc is a reflective subcategory of the category of topoi and geometric morphisms (cf. [xviii] Locale 4.14). Thus, an adjunction $Sub \dashv \mathbf{Sh}$ : Loc $\hookrightarrow$ Topos follows. Here, $Sub$ takes each topos to the locale of its subterminals.

- A point of a locale $L$ is a locale morphism $p : 1 \longrightarrow L$, where 1 denotes the locale corresponding to the one-point space. Equivalently, it is a class of isomorphic points of the topos $\mathbf{Sh}(L)$. In the former case, an open subspace $U$ of $L$ contains a point $p$ of $L$ if and only if $p^*(U) = 1$; here $p^* : \mathcal{O}(L) \longrightarrow \{0, 1\}$ is the frame morphism corresponding to the locale morphism $p$.

- A locale $L$ is said to be *spatial* or to have *enough points* if it is locale isomorphic to the locale corresponding to the frame of opens of its space of points. More explicitly, a locale $L$ is said to have enough points when for any two opens $U$ and $V$ in $L$, $U \subseteq V$ if and only if every point of $L$ that belongs to $U$ also belongs to $V$ (cf. [xxi] II 1.5). It is worthwhile to mention that the category of spatial locales and locale maps ($\mathsf{Loc}^{\mathrm{sob}}$) is equivalent to the category of sober topological spaces and continuous maps (cf. [xxi] II Corollary 1.7).

---
[1] We would otherwise have to care about the indexing of the points by some set and use slightly different language; it is only for this purpose that we make this hypothesis.

- An element $a \in L$ of a locale $L$ is said to be *finite* if whenever $a = \bigvee_{\alpha \in I} b_\alpha$, then there exists a finite subset $J \subseteq I$ such that $a = \bigvee_{\alpha \in J} b_\alpha$. A locale $L$ is said to be *coherent* if:

  1. The top element is finite.
  2. Every element of $L$ can be expressed as a join of finite elements.
  3. The meet of any two finite elements is again finite.

  A *coherent frame* is a frame of opens of a coherent locale.

- A topological space $X$ is called *spectral* if it is sober, admits a base of quasicompact open subsets which is closed under finite intersections, and is quasicompact. The quintessential example is the topological space underlying any concentrated scheme.

- A map of topological spaces $f : X \longrightarrow Y$ is called a *spectral map* if for any open subset $U \subseteq Y$ which is quasicompact, $f^{-1}(U)$ is also quasicompact. Naturally, a subspace of a spectral space $X$ is a topological subspace $U \hookrightarrow X$ where inclusion is a spectral map. We define a category $\mathcal{S}p$ of spectral spaces and spectral maps.

- The category $\mathcal{S}p$ is equivalent to image of the functor $\mathrm{Spec} : \mathsf{CRing}^{\mathrm{op}} \longrightarrow \mathsf{Top}$. This lifts from M. Hochster's seminal paper on spectral spaces. In fact, per [xxii] Theorem 6, the full statement of the result is as follows:

  **Theorem 1.1.3.** *$\mathcal{S}p$ is equivalent to the image of* $\mathrm{Spec}$. *Furthermore,* $\mathrm{Spec}$ *is invertible at the following subcategories of $\mathcal{S}p$:*

  1. *The subcategory of all spectral spaces and surjective spectral maps.*
  2. *For a spectral space $X$, the subcategory of its spectral subspaces and inclusions of these.*
  3. *The full subcategory of $\mathcal{S}p$ whose objects are the $T_1$ spectral spaces.*

- Every topological space determines a locale through the assignment $X \longmapsto \mathcal{O}(X)$. On the other hand, each locale determines a topological space through the space-of-points construction, $pt : \mathsf{Loc} \longrightarrow \mathsf{Top}$. This description is an idempotent adjunction $\mathcal{O} \dashv pt : \mathsf{Loc} \longrightarrow \mathsf{Top}$ (cf. [xviii] Locale 4.4). Additionally, the adjunction restricts to an equivalence of categories $\mathsf{Loc}^{\mathrm{coh}} \longrightarrow \mathcal{S}p$, where the domain is the subcategory (of $\mathsf{Loc}$) of coherent locales and locale maps that preserve finite elements under preimages ($\mathsf{Frm}^{\mathrm{coh}} := (\mathsf{Loc}^{\mathrm{coh}})^{\mathrm{op}}$). In fact, this is but a shadow of Stone duality which states that the category of distributive lattices, $\mathsf{DLat}$, is dual to $\mathcal{S}p$ (cf. [xxi] II, 3.3). Thus, we have a functor $\mathrm{Spec} : \mathsf{DLat}^{\mathrm{op}} \longrightarrow \mathsf{Top}$.

- We lift from K. Aoki the following facts useful to us (cf. [xxiii] 3.3 and 3.4 respectively). Their proofs are found in the same place.

  **Proposition 1.1.4.** *The inclusion* $\mathsf{Loc}^{\mathrm{coh}} \hookrightarrow \mathsf{Loc}$ *preserves (small) limits.*

  **Proposition 1.1.5.** *The functor* $\mathrm{Spec}: \mathsf{DLat}^{\mathrm{op}} \longrightarrow \mathsf{Top}$ *preserves (small) limits.*

- Classical Stone duality establishes a representable functor $\mathsf{CH}^{\mathrm{op}} \longrightarrow \mathsf{BoolAlg}$ between the category of compact Hausdorff spaces and continuous maps and the category of Boolean algebras and Boolean homomorphisms.

- A *Stone space* is a compact, Hausdorff, and totally disconnected topological space. One has a full inclusion Stone $\hookrightarrow$ CH of the subcategory of Stone spaces. The representable functor of the previous bullet point restricts to an equivalence of categories Stone$^{\mathrm{op}} \longrightarrow$ BoolAlg (cf. [xviii] Boolean Algebra 3.1). Furthermore, via Stone-Čech compatification, one obtains a functor Top$^{\mathrm{op}} \longrightarrow$ BoolAlg which again restricts to the aforementioned equivalence. It is noted also that the category of Boolean rings and ring homomorphims, BoolRing, is equivalent to BoolAlg (cf. [xviii] Boolean Ring 3.1).

$$\infty\text{-Prosets}$$

We introduce $\infty$-prosets and prove several results—first, some to ground us in this new world, and later, others that carry importance relative to this project. Naturally, we have a forgetful functor $U :$ Proset $\longrightarrow$ Set from category of prosets and proset maps (which is faithful). And for ends which will become clear later in this article, it will be useful for us to fix three Grothendieck universes: $\mathbb{U}_0 \in \mathbb{U}_1 \in \mathbb{U}_2$. We will call the elements of $\mathbb{U}_0$ *small*, those of $\mathbb{U}_1$ *large* and, those of $\mathbb{U}_2$ *very large*. Cat$_\infty$ will indicate the $\infty$-category of small $\infty$-categories, while $\widehat{\mathsf{Cat}_\infty}$ will indicate the $\infty$-category of large $\infty$-categories. Henceforth, we shall also treat the nerves of all the ordinary 1-categories mentioned in the previous subsection as being objects of Cat$_\infty$. And while the results we obtain here (to apply to later parts of this paper) involve Cat$_\infty$, they hold as true if we instead focused on $\widehat{\mathsf{Cat}_\infty}$. We will assume all prosets to be small. We will write $X^\leq$ for a simplicial proset where the superscript signals *particular* pointwise proset structures. For example, $X^{\leq_1}$ and $X^{\leq_2}$ are two simplicial prosets that become the same simplicial set once acted on by the forgetful functor.

**Definition 1.2.1.** Let $X^\leq$ be a simplicial proset and $X$ its composition with $U$. Assume that $X$ is an $\infty$-category. We will say $X$ is *n-commensurate* with face maps $\partial_i : X_n^\leq \longrightarrow X_{n-1}^\leq$ if the homotopies in $X_{n-1}$ are created jointly by the former maps. By the latter, we mean that any two $n$-simplices, which map to homotopic $(n-1)$-simplices of $X$ under the action of $\partial_i$ for all $0 \leq i \leq n$, are equivalent as objects of $X_n^\leq$. This implies homotopic $(n-1)$-simplices are equivalent as objects of $X_{n-1}^\leq$. Taking the case $n = 1$ as an example, if two objects are equivalent as objects of $X$, then the maps expressing this equivalence are themselves equivalent as objects of $X_1^\leq$. We also say that a map $\mathsf{N}(\Delta)^{\mathrm{op}} \longrightarrow \mathcal{C}$ is a Segal space (of $\mathcal{C}$) if it satisfies the Segal condition and is complete (in the context of $\mathcal{C}$). Here, $\mathcal{C}$ is an $\infty$-category that admits pullbacks.

**Proposition 1.2.2.** *Let $X^\leq$ be a simplicial proset and $X$ its composition with $U$. Then the following are equivalent.*

1. *$X^\leq$ determines a Segal space $\mathsf{N}(\Delta)^{op} \longrightarrow \mathsf{Cat}_\infty$.*

2. *$X^\leq$ determines a map $\mathsf{N}(\Delta)^{op} \longrightarrow \mathsf{Cat}_\infty$ that satisfies the Segal condition, $X$ is an $\infty$-category and is 1-commensurate with the face maps $\partial_i : X_1^\leq \longrightarrow X_0^\leq$.*

*Proof.* To see that (1) $\Rightarrow$ (2), consider what completeness tells us about the map $\mathsf{N}(s_0) : \mathsf{N}(X_0^\leq) \longrightarrow \mathsf{N}((X_1^\leq)^{\mathrm{eq}})$. That it is a categorical equivalence implies that it is essentially surjective. Assuming that $X$ is an $\infty$-category, suppose that $f : x \longrightarrow y$ is an equivalence. Then $f \in \mathsf{N}((X_1^\leq)^{\mathrm{eq}})$, and this means that it is equivalent to some $\mathsf{N}(s_0)$. Recalling the characterization of 1-categorical equivalences, we conclude that $f$ is equivalent to $g : y \longrightarrow x$ where $f \circ g \simeq \mathrm{id}_x$ and $g \circ f \simeq \mathrm{id}_y$ in $X$. It remains to show that, indeed, $X$ is an $\infty$-category. Applying the groupoidification functor pointwise, we obtain a complete Segal space. Hence $X$ is determined as an $\infty$-category.

Given (2), it remains to show completeness. Notice that the nerve of the map $s_0 : X_0^{\leq} \longrightarrow X_1^{\leq}$ is automatically fully faithful. Ergo, considering $(X_1^{\leq})^{\text{eq}}$, the sub-proset whose objects correspond to equivalences in $X$, we need only show essential surjectivity to show categorical equivalence of the nerve of the map $s_0 : X_0^{\leq} \longrightarrow (X_1^{\leq})^{\text{eq}}$. This is given by the condition imposition on homotopies of $X_0$. Therefore, taking the nerve pointwise, we obtain the desired Segal space $\mathsf{N}(\Delta)^{\text{op}} \longrightarrow \mathsf{Cat}_\infty$.

□

We define $\infty$-prosets from the above conclusion. That is, as simplicial prosets that are at the same time Segal space objects of $\mathsf{Cat}_\infty$. Additionally, 1.2.2 allows us to define an $\infty$-category of $\infty$-prosets, $\mathsf{Pros}_\infty$, as the full subcategory of $\mathsf{Fun}(\mathsf{N}(\Delta)^{\text{op}}, \mathsf{Cat}_\infty)$ spanned by Segal space objects that, when composed with any $\Delta^0 \longrightarrow \mathsf{N}(\Delta)^{\text{op}}$, have nerves of prosets as their images. Now, given an $\infty$-proset $X^{\leq}$, we will refer to $X = U \circ X^{\leq}$ as the $\infty$-category underlying the $\infty$-proset structure. This suggests a projection $\mathsf{Pros}_\infty \longrightarrow \mathsf{Cat}_\infty$ that we exhibit below.

**Proposition 1.2.3.** *There exists a map* $\mathsf{u} : \mathsf{Pros}_\infty \longrightarrow \mathsf{Cat}_\infty$ *of* $\infty$-*categories.*

*Proof.* Applying Lemma 1.4.2 of [xxxi], starting with a functor $g : \mathsf{Cat}_\infty \longrightarrow \mathsf{Kan}$ right adjoint to the inclusion $\mathsf{Kan} \hookrightarrow \mathsf{Cat}_\infty$. Hence, we induce a map $\mathsf{Pros}_\infty \longrightarrow \mathsf{CSS}$. The latter is the $\infty$-category of complete Segal spaces. From here, the categorical equivalence $\mathsf{CSS} \simeq \mathsf{Cat}_\infty$, gives the desired map. □

**Proposition 1.2.4.** $\mathsf{Pros}_\infty$ *admits small limits.*

*Proof.* Because we have described $\mathsf{Pros}_\infty$ as a full subcategory of an $\infty$-category of $\infty$-functors into a complete $\infty$-category, we can compute small limits pointwise. By Example 2.2.4 of [xxx], we find that, in fact, we can compute limits in $\mathsf{Cat}$ since, pointwise, we have small 1-categories. Moreover, Proset admits small limits. It remains to show that Segal and completeness conditions are upheld by taking limits. But in [xxi], this is the implication of Corollary 1.3.4 given Theorem 1.4.1. □

Consider two sets $S_i$ and a pair of parallel surjective arrows $S_1 \rightrightarrows S_0$. Then if $S_i$ is a proset, there exists a proset structure on $S_j$, where $i \neq j$, such that the pair of parallel arrows are both proset maps. In the case $i = 0$, declare for $s$ and $s'$ in $S_1$, $s \leq s'$ if and only if each of the arrows is monotone on the relationship. On the flip side, if $i = 1$, the construction proceeds as follows. First, given $a$ and $b$ in $S_0$, write $a \sim b$ if and only if there exist $u$ and $v$ in $S_1$ such that $u \leq v$ and under the action of the parallel pair of arrows, $(u, v)$ maps to $(a, b)$ for at least one of the arrows. Second, declare $a \leq_0 b$ if and only if there exists a chain $a \sim r_0 \sim \cdots \sim r_N \sim b$ where $N$ is finite. In either case, we will call the resulting pair of prosets maps a *parallel proset structure*. In our situation, we will say that, given a small $\infty$-category $\mathcal{C}$, the pair

$$\mathcal{C}_1 \overset{\partial_1}{\underset{\partial_0}{\rightrightarrows}} \mathcal{C}_0$$

is commensurate with $\mathcal{C}$ if one of them admits a proset structure that results in a parallel proset structure that $\mathcal{C}$ is 1-commensurate with.

**Proposition 1.2.5.** *Let $\mathcal{C}$ be a small $\infty$-category and $0 \leq m \leq 1$. Suppose that $\mathcal{C}_m$ admits a proset structure $\mathcal{C}_1^{\leq} = (\mathcal{C}_1, \leq)$ that results in a parallel proset structure $\partial : \mathcal{C}_1 \longrightarrow \mathcal{C}_0$. Then there exists an $\infty$-proset $\mathcal{C}^{\leq}$ having $\mathcal{C}$ as the underlying $\infty$-category and such that $[m] \longmapsto \mathsf{N}(\mathcal{C}_m^{\leq})$.*

*Proof.* We show that for $n \geq 0$ each $\mathcal{C}_n$ admits a proset structure and that all degeneracy maps $s_n : \mathcal{C}_n \longrightarrow \mathcal{C}_{n+1}$ and face maps $\partial_n : \mathcal{C}_{n+1} \longrightarrow \mathcal{C}_n$ respect the proset structure. It suffices to show this for $m = 1$. We proceed inductively by showing that this is the case with $s_i : \mathcal{C}_1 \longrightarrow \mathcal{C}_2$ and $\partial_k : \mathcal{C}_2 \longrightarrow \mathcal{C}_1$ for $0 \leq k \leq 2$. Given two 2-simplices, $\sigma$ and $\sigma'$, we say $\sigma \leq_* \sigma'$ if and only if $\partial_k \sigma \leq \partial_k \sigma'$ for all $k$. Looking through the axioms that describe prosets, one concludes that $\mathcal{C}_2^{\leq} = (\mathcal{C}_2, \leq_1)$ is a proset. Furthermore, all face maps are a priori proset maps while for degeneracy maps, the simplicial identities guarantee the desired outcome. Mutatis mutandis, a replica of this argumentation can now be made in the case of $\mathcal{C}_n$ and $\mathcal{C}_{n+1}$ after the inductive step is taken. The key observation about the transition maps remains the same.

It remains to show that the resulting simplicial proset satisfies the Segal condition. Note that we can calculate limits of nerves of small categories in $\mathsf{Cat}$ (cf. [xxx] 2.2.4). More so in our case, given the localization $\mathsf{Proset} \longrightarrow \mathsf{Cat}$ creates limits. By unraveling the construction above, one notices that for $n \geq 2$, $n$-simplices are determined up to equivalence in $\mathcal{C}_n^{\leq}$ by their boundaries in the case of $n > 2$ and by the inner horn when $n = 2$. But this is exactly the statement of the Segal condition in our situation. $\square$

This conclusion informs us that $\infty$-prosets are determined entirely from the proset structures on their 0-simplices and 1-simplices. This should not come as a surprise given 1.2.2; which to begin with establishes that the 1-simplices determine the proset structure on all $n$-simplices where $n \geq 2$. What is left then is the proset structure on the 0-simplices. And this is more or less determined by their interactions with the proset of the 1-simplices through the face maps. Moving forward, we will say that an $\infty$-proset $\mathcal{C}^{\leq}$ is determined from $\mathcal{C}_0^{\leq}$ if it arises in the process of Proposition 1.2.5.

**Remark 1.2.6** The construction of the proset $\mathcal{C}_0^{\leq}$ from the proset $\mathcal{C}_1^{\leq}$ has the added advantage that if one starts with $\mathcal{C}_0^{\leq}$ and using the procedure highlighted in the passage before Proposition 1.2.5 one obtains a proset structure on $\mathcal{C}_1$, one ends up with $\mathcal{C}_1^{\leq}$ whenever $\mathcal{C}^{\leq}$ is determined from $\mathcal{C}$ as in 1.2.8.

**Remark 1.2.7.** If one takes any $\infty$-proset $\mathcal{C}^{\leq}$, isolates the proset structure of its $m$-simplices (where $0 \leq m \leq 1$) and then proceeds in the way of 1.2.5, one does not always end up with the original $\infty$-proset. However, as we shall see in 1.2.9, since the prosets of the objects of the now two $\infty$-prosets are identical, a map of $\infty$-prosets is induced between them; specifically, a map from the original to the new proset that is, moreover, an inclusion of $\infty$-prosets. This intimates a "local final" $\infty$-proset.

The following is the $\infty$-categorical reflection of the fact that every (small) ordinary category $C$ naturally determines a proset $(C, \leq)$ that encodes the relationships between its objects. In Theorem 1.2.14 we describe exactly what is meant by this "naturalness".

**Proposition 1.2.8.** *Let $\mathcal{C}$ be a small $\infty$-category. Then there exists an $\infty$-proset $\mathcal{C}^{\leq}$ such that $\mathcal{C}$ is the underlying $\infty$-category.*

*Proof.* By the previous proposition, it suffices to show that $\mathcal{C}_1 \rightrightarrows \mathcal{C}_0$ admits a parallel proset structure that $\mathcal{C}$ is commensurate with. We define this structure as follows. Given two 0-simplices, $a$ and $b$, we write $a \leq b$ if and only if there exists a 1-simplex $\sigma : \Delta^1 \longrightarrow \mathcal{C}$ such that $\partial_0 \sigma = b$ and $\partial_1 \sigma = b$: we say $\sigma$ witnesses $a \leq b$. It follows immediately that since $\partial_k s_0 = \mathrm{id}$ for $k \in \{0, 1\}$, we must have $a \leq a$ for any $a : \Delta^0 \longrightarrow \mathcal{C}$. Now suppose $a \leq b$ is witnessed by a 1-simplex $\sigma$ and $b \leq c$ is witnessed by a 1-simplex $\alpha$. Then let $\omega : \Lambda_1^2 \longrightarrow \mathcal{C}$ be the inner horn with $\partial_2 \omega = \sigma$ and $\partial_0 \omega = \alpha$. By the inner horn filling property of $\infty$-categories, we have a 2-simplex $\omega'$ that extends $\omega$ through the inclusion $\Lambda_1^2 \subseteq \Delta^2$. Therefore, we obtain a 1-simplex $\partial_1 \omega'$ that witnesses $a \leq c$. $\square$

**Proposition 1.2.9.** *Let $f : \mathcal{C} \longrightarrow \mathcal{D}$ be a map of small $\infty$-categories. Suppose that $\mathcal{C}_0$ and $\mathcal{D}_0$ admit proset structures and that $f_0 : \mathcal{C}_0^{\leq} \longrightarrow \mathcal{D}_0^{\leq}$ is a map of prosets. Then $f$ upgrades to a map $f : \mathcal{C}^{\leq} \longrightarrow \mathcal{D}^{\leq}$ of the $\infty$-prosets determined from $\mathcal{C}_0^{\leq}$ and $\mathcal{D}_0^{\leq}$.*

*Proof.* It suffices to show that $f_n : \mathcal{C}_n^{\leq} \longrightarrow \mathcal{D}_n^{\leq}$ is a map of prosets. We will do this inductively. For the base case, let $n = 1$. Suppose $\lambda \leq \lambda'$ in $\mathcal{C}_1^{\leq}$. Then we have $f_0(\partial_k \lambda) \leq f_0(\partial_k \lambda')$ for all $0 \leq k \leq 1$. But by definition, $(f_0(\partial_k \lambda), f_0(\partial_k \lambda')) = (\partial_k(f_1 \lambda), \partial_k(f_1 \lambda'))$. Hence $\partial_k(f_1 \lambda) \leq \partial_k(f_1 \lambda')$. Considering how $\mathcal{D}_1^{\leq}$ comes about, this implies $f_1 \lambda \leq f_1 \lambda'$. Indeed, for $(\lambda, \lambda') \in \mathcal{C}_n \times \mathcal{C}_n$, since $(f_n(\partial_k \lambda), f_n(\partial_k \lambda')) = (\partial_k(f_{n+1} \lambda), \partial_k(f_{n+1} \lambda'))$ for all $0 \leq k \leq n$ and $n \geq 0$, if $f_n$ is a map of prosets then $f_{n+1}$ is also a map of prosets given how the proset $\mathcal{D}_{n+1}^{\leq}$ is determined from $\mathcal{D}_n^{\leq}$. $\square$

**Lemma 1.2.10.** *Let $f : \mathcal{C} \longrightarrow \mathcal{D}$ be a categorical equivalence of small $\infty$-categories. Suppose $\mathcal{C}$ is the underlying $\infty$-category of some $\infty$-proset $\mathcal{C}^{\leq}$. Then there exists an $\infty$-proset $\mathcal{D}^{\leq}$ having $\mathcal{D}$ as its underlying category and such that $\mathsf{N}(f_k) : \mathsf{N}(\mathcal{C}_k^{\leq}) \longrightarrow \mathsf{N}(\mathcal{D}_k^{\leq})$ is a categorical equivalence for $k \in \{0, 1\}$.*

*Proof.* We begin with the equivalence of homotopy categories $\mathrm{h}f : \mathrm{h}\mathcal{C} \longrightarrow \mathrm{h}\mathcal{D}$. Notice that there is a proset structure on the objects of $\mathrm{h}\mathcal{C}$ because there is one on $\mathcal{C}_0$. If we declare that for any pair $(a, b)$ of objects of $\mathrm{h}\mathcal{D}$, $a \leq' b$ if and only there exists $a'$ and $b'$ in $\mathrm{h}\mathcal{C}$ such that $a' \leq b'$ while $\mathrm{h}f(a') \simeq a$ and $\mathrm{h}f(b') \simeq b$, then we realize a proset structure on the objects of $\mathrm{h}\mathcal{D}$. It follows that since $\mathrm{h}$ is essentially surjective, the resulting proset is equivalent to the initial one. Now, mutatis mutandis, a similar argument shows the existence of a proset structure on the set of morphisms of $\mathrm{h}\mathcal{D}$ and its equivalence to that on the set of morphisms of $\mathrm{h}\mathcal{C}$. It remains to observe that through [i] Proposition 1.2.3.9, the equivalences above can be lifted to the desired equivalences; one needs only treat objects (resp. 1-simplices) as related to each other if and only if they admit equivalent (resp. homotopic) counterparts which are related. To complete the proof, recall Proposition 1.2.5 and apply it to $\mathcal{D}$ for the case $m = 1$. The hypothesis therein is met since the categorical equivalences map parallel proset structures to each other: as witnessed by the equivalences described above induced by $\mathrm{h}f$. $\square$

**Proposition 1.2.11.** *Let $f : \mathcal{C} \longrightarrow \mathcal{D}$ be a categorical equivalence of small $\infty$-categories. Suppose $\mathcal{C}$ is the underlying $\infty$-category of some $\infty$-proset $\mathcal{C}^{\leq}$. Then there exists an $\infty$-proset $\mathcal{D}^{\leq}$ having $\mathcal{D}$ as its underlying category and such that $f : \mathcal{C}^{\leq} \longrightarrow \mathcal{D}^{\leq}$ is an equivalence of $\infty$-prosets.*

*Proof.* We deduce from Proposition 1.2.2 that the $n^{\text{th}}$-Segal map

$$\mu_n^{\text{Seg}}(\mathcal{C}) : \mathsf{N}(\mathcal{C}_n^{\leq}) \longrightarrow \underbrace{\mathsf{N}(\mathcal{C}_1^{\leq}) \times_{\mathsf{N}(\mathcal{C}_0^{\leq})} \cdots \times_{\mathsf{N}(\mathcal{C}_0^{\leq})} \mathsf{N}(\mathcal{C}_1^{\leq})}_{n \text{ factors}}$$

is a categorical equivalence. The same holds for $\mathcal{D}^{\leq}$ constructed in 1.2.10. We know already from Lemma 1.2.10 that $\mathsf{N}(f_k) : \mathsf{N}(\mathcal{C}_k^{\leq}) \longrightarrow \mathsf{N}(\mathcal{D}_k^{\leq})$ is a categorical equivalence for $k \in \{0, 1\}$. Hence, taking into consideration the Segal maps, an equivalence is induced $P_n : \mathsf{N}(\mathcal{C}_n^{\leq}) \longrightarrow \mathsf{N}(\mathcal{D}_n^{\leq})$ for all $n \geq 0$. If we then perform induction on $n$ with the base case being $n = 2$, we are led to conclude that $P_n \circ \mu_n^{\text{Seg}}(\mathcal{C}) \simeq \mu_n^{\text{Seg}}(\mathcal{D}) \circ f_n$ as maps of categories. Subsequently, by the 2-out-of-3 property of categorical equivalences, $f_n$ must be an equivalence. Proposition 1.2.9 now gives the result. $\square$

**Corollary 1.2.12.** *The following are equivalent.*

1. $f : \mathcal{C}^{\leq} \longrightarrow \mathcal{D}^{\leq}$ *is an equivalence of $\infty$-prosets.*

2. $f : \mathcal{C}^{\leq} \longrightarrow \mathcal{D}^{\leq}$ *is a map of $\infty$-prosets such that* $\mathsf{N}(f_k) : \mathsf{N}(\mathcal{C}_k^{\leq}) \longrightarrow \mathsf{N}(\mathcal{D}_k^{\leq})$ *is a categorical equivalence for $k \in \{0, 1\}$.*

*Proof.* That $1 \Rightarrow 2$ is immediate. Notice that the hypothesis in 2 achieves the same effect as Lemma 1.2.10 in the proof of Proposition 1.2.11. This is not generally true. $\square$

**Proposition 1.2.13.** *There exists a map* $\mathsf{v} : \mathsf{Cat}_\infty \longrightarrow \mathsf{Pros}_\infty$ *of $\infty$-categories.*

*Proof.* Starting with the inclusion $\mathsf{Kan} \longrightarrow \mathsf{Cat}_\infty$, and recalling $\mathsf{CSS} \simeq \mathsf{Cat}_\infty$, we obtain a map $\mathsf{Cat}_\infty \longrightarrow \mathsf{Fun}(\mathsf{N}(\Delta)^{\text{op}}, \mathsf{Cat}_\infty)$. For each small $\infty$-category $\mathcal{C}$, the image under this map acts on finite ordinals as $[n] \longmapsto \mathsf{Fun}(\Delta^n, \mathcal{C})^{\simeq}$ where the image denotes the largest simplicial subset of $\mathsf{Fun}(\Delta^n, \mathcal{C})$ that is a Kan complex (cf. [xvii] Complete Segal Spaces, Proposition 4.7). Invoking Proposition 1.2.8 after limiting ourselves to objects of the kind $\mathsf{Fun}(\Delta^n, \mathcal{C})^{\simeq}$ for all $n$ and $\mathcal{C} \in \mathsf{Cat}_\infty$, we obtain the desired map when we compose with $\mathsf{h} : \mathsf{Cat}_\infty \longrightarrow \mathsf{Cat}$. $\square$

We remark that in above argument, we need Proposition 1.2.8 to be certain that composition with h results in a Segal space object of $\mathsf{Cat}_\infty$.

**Theorem 1.2.14.** *The map* $\mathsf{v} : \mathsf{Cat}_\infty \longrightarrow \mathsf{Pros}_\infty$ *is left adjoint to* $\mathsf{u} : \mathsf{Pros}_\infty \longrightarrow \mathsf{Cat}_\infty$.

*Proof.* We begin by unwinding the definition of $\mathsf{v} : \mathsf{Cat}_\infty \longrightarrow \mathsf{Pros}_\infty$. This is in fact the map obtained by post-composition with the sequence

$$\mathsf{Cat}_\infty \xrightarrow{\mathsf{h}} \mathsf{N}(\mathsf{Cat}) \to \mathsf{N}(\mathsf{Proset})$$

where the last functor is the 1-categorical version of Proposition 1.2.8. We observe that each of these functors are left adjoint. Therefore, they preserve colimits. Moreover, the sequence above is that of presentable $\infty$-categories. Now, since the computation of colimits is done pointwise (cf. [xxvi] 7.1.7.2), we conclude that the induced map $\mathsf{Fun}(\mathsf{N}(\Delta)^{\text{op}}, \mathsf{Cat}_\infty) \longrightarrow \mathsf{Fun}(\mathsf{N}(\Delta)^{\text{op}}, \mathsf{N}(\mathsf{Proset}))$ preserves colimits. Applying Proposition 5.5.3.6 of [i], the adjoint functor theorem concludes that this map is left adjoint. Furthermore, its right adjoint necessarily preserves Segal space objects since it preserves limits. But the restriction of the map to the Segal space objects arising from the map $\mathsf{Cat}_\infty \longrightarrow \mathsf{Fun}(\mathsf{N}(\Delta)^{\mathrm{op}}, \mathsf{Cat}_\infty)$ is $\mathsf{v} : \mathsf{Cat}_\infty \longrightarrow \mathsf{Pros}_\infty$. By right adjointness of the map $g : \mathsf{Cat}_\infty \longrightarrow \mathsf{Kan}$, we see that $\mathsf{v} \dashv \mathsf{u}$. □

This deduction is an analogue of the classical adjunction between $\mathsf{Cat} \longrightarrow \mathsf{Proset}$, the functor that stimulates 1.2.8, and the inclusion $\mathsf{Proset} \longrightarrow \mathsf{Cat}$. More importantly, it gives validation that the construction of $\mathsf{Pros}_\infty$ bears an $\infty$-categorical nature. This in mind, we continue to specialized phenomena that we see to usefulness, one way or the other, in later parts of this article.

**Definition 1.2.15.** Let $\mathcal{C}^\leq$ be an $\infty$-proset. We say $\mathcal{C}^\leq$ is *connected* whenever the map of sets $(\partial_1 \times \partial_0) : \mathcal{C}_1 \longrightarrow \mathcal{C}_0 \times \mathcal{C}_0$ induces a categorical equivalence $\mathsf{N}(\partial_1 \times \partial_0) : \mathsf{N}(\mathcal{C}_1^\leq) \longrightarrow \mathsf{N}(\mathcal{C}_0^\leq) \times \mathsf{N}(\mathcal{C}_0^\leq)$.

Connected $\infty$-prosets are of interest because, as we shall see in the next subsection, their behavior under Stone-type dualities is more accessible than otherwise. Moreover, Proposition 1.2.15 informs us that when dealing with pointed $\infty$-categories, we should expect that the $\infty$-prosets they determine (via 1.2.13) are connected: recall that stable $\infty$-categories are pointed.

**Lemma 1.2.16.** *Let $\mathcal{C}$ be a small connected $\infty$-category that admits a prosite $\mathcal{C}_0^\leq$. Then the $\infty$-prosite $\mathcal{C}^\leq$ determined from $\mathcal{C}_0^\leq$ is connected.*

*Proof.* It is sufficient to presume that $\mathcal{C}$ is a Kan complex. To see this, evaluate the Kan fibrant replacement $\mathcal{C} \mapsto [\mathcal{C}]$ and take Lemma 1.2.10 into account. With this consideration, recall that $\mathcal{C}$ is connected precisely if the equivalence relation on $\mathcal{C}_0$ induced by

$$\mathcal{C}_1 \underset{\partial_0}{\overset{\partial_1}{\rightrightarrows}} \mathcal{C}_0 \times \mathcal{C}_0$$

is the singleton set. This means that given two objects of $\mathcal{C}$, say $a$ and $b$, there exists a 1-simplex $\alpha : \Delta^1 \longrightarrow \mathcal{C}$ with $(\partial_1 \alpha, \partial_0 \alpha) = (a, b)$. This means that the map $(\partial_1 \times \partial_0) : \mathcal{C}_1 \longrightarrow \mathcal{C}_0 \times \mathcal{C}_0$ is surjective. Consequently, $\mathsf{N}(\partial_1 \times \partial_0) : \mathsf{N}(\mathcal{C}_1^\leq) \longrightarrow \mathsf{N}(\mathcal{C}_0^\leq) \times \mathsf{N}(\mathcal{C}_0^\leq)$ admits a section which is essentially surjective. This follows from the axiom of choice and how the proset structure on $\mathcal{C}_1^\leq$ arises from that on $\mathcal{C}_0^\leq$. □

**Proposition 1.2.17.** *Let $\mathcal{C}$ be a small connected $\infty$-category. Suppose that $\mathcal{C}^\leq$ is determined from a finitely complete prosite $\mathcal{C}_0^\leq$. Then each $\mathcal{C}_n^\leq$ is finitely complete, and the transition maps are left exact for all $n \geq 0$. Specifically, face maps jointly create finite limits.*

*Proof.* Let us begin with the assumption that $\mathcal{C}_1^\leq$ is finitely complete and that $s_0 : \mathcal{C}_0^\leq \longrightarrow \mathcal{C}_1^\leq$ is left exact. We then notice that by the Segal condition, each $\mathcal{C}_n^\leq$ is finitely complete since each is a finite limit of finitely complete small categories and left exact functors. It also follows that the face maps $\partial : \mathcal{C}_{n+1}^\leq \longrightarrow \mathcal{C}_n^\leq$ are equivalent to projections $\mathcal{C}_n^\leq \times_{\mathcal{C}_0^\leq} \mathcal{C}_1^\leq \longrightarrow \mathcal{C}_n^\leq$ and are therefore left exact since the limits involved are computed from pointwise limits. Moreover, the latter informs us that face maps create finite limits jointly. On the other hand, each degeneracy map $\mathcal{C}_n^\leq \longrightarrow \mathcal{C}_{n+1}^\leq$ is equivalent to a map $(\mathrm{id} \times s_0) : \mathcal{C}_n^\leq \times \{*\} \longrightarrow \mathcal{C}_{n+1}^\leq$ and is, therefore, left exact.

It is left to show that the assumption we opened with is necessarily true given the premise about $\mathcal{C}$. Lemma 1.2.16 reduces this task so that we need only check that $\mathcal{C}_0^\leq$ is finitely complete. □

**Proposition 1.2.18.** *Let $\mathcal{C}^{\leq}$ be an $\infty$-proset determined from its underlying $\infty$-category $\mathcal{C}$ through the map $\mathsf{v} : \mathsf{Cat}_\infty \longrightarrow \mathsf{Pros}_\infty$. The following are equivalent.*

1. *$\mathcal{C}^{\leq}$ is connected.*

2. *$\mathcal{C}$ is connected.*

*Proof.* It is sufficient to assume that $\mathcal{C}$ is a Kan complex. By default, the categorical equivalence $\mathsf{N}(\partial_1 \times \partial_0) : \mathsf{N}(\mathcal{C}_1^{\leq}) \longrightarrow \mathsf{N}(\mathcal{C}_0^{\leq}) \times \mathsf{N}(\mathcal{C}_0^{\leq})$ tells us that given any object $(a,b) \in \mathsf{N}(\mathcal{C}_0^{\leq}) \times \mathsf{N}(\mathcal{C}_0^{\leq})$, there exists an $f \in \mathsf{N}(\mathcal{C}_1^{\leq})$ such that $(\partial_1(f), \partial_0(f)) \leq (a,b)$. Unwinding the map $\mathsf{v} : \mathsf{Cat}_\infty \longrightarrow \mathsf{Pros}_\infty$, we see that the former implies the existence of two 1-simplices of $\mathcal{C}$, $\lambda'$ and $\lambda$, such that $(\partial_1 \lambda, \partial_0 \lambda) = (\partial_1(f), a)$ and $(\partial_1 \lambda', \partial_0 \lambda') = (\partial_0(f), b)$. If we recall that $\mathcal{C}$ is a Kan complex, we obtain a 1-simplex $\alpha : \Delta^1 \longrightarrow \mathcal{C}$ with $(\partial_1 \alpha, \partial_0 \alpha) = (a, b)$. Now for $2 \Rightarrow 1$, notice that $\mathcal{C}^{\leq}$ is determined from the prosite $\mathcal{C}_0^{\leq}$ as described by the map $\mathsf{v} : \mathsf{Cat}_\infty \longrightarrow \mathsf{Pros}_\infty$. Thus, we may apply Lemma 1.2.16. $\square$

**Proposition 1.2.19.** *Let $f : \mathcal{C} \longrightarrow \mathcal{D}$ be a map of small connected $\infty$-categories that admits an equivalence of prosets $f_0 : \mathcal{C}_0^{\leq} \longrightarrow \mathcal{D}_0^{\leq}$. Then the map $f : \mathcal{C}^{\leq} \longrightarrow \mathcal{D}^{\leq}$ between the $\infty$-prosets determined from $\mathcal{C}_0^{\leq}$ and $\mathcal{D}_0^{\leq}$ is an equivalence of $\infty$-prosets.*

*Proof.* Combine Lemma 1.2.16 and Corollary 1.2.12. $\square$

### Localic Quantization of Space

In this subsection, we present the $\infty$-category of $\infty$-prosites. In doing so, we hope to welcome into the world of $\infty$-prosets the practices of pointless topology—in the sense of classical Stone-type dualities and O. Caramello's work in [xi]. In culmination, once the latter technologies are suited in place, we should find ourselves a framework from which should emerge a couple functors into animated $S$-stacks.

**Definition 1.3.1.** Let the dyad $(\mathcal{C}^{\leq}, J)$ be the data constituted as follows. For each $n \geq 0$ and all $0 \leq k \leq n$, we have that:

1. $\mathcal{C}_n^{\leq}$ is endowed with a Grothendieck topology $J_n$. Hence a prosite $(\mathcal{C}_n^{\leq}, J_n)$ which we will call the $n^{th} - prosite$.

2. The transition maps $\partial_k : \mathcal{C}_n^{\leq} \longrightarrow \mathcal{C}_{n-1}^{\leq}$ and $s_k : \mathcal{C}_n^{\leq} \longrightarrow \mathcal{C}_{n+1}^{\leq}$ are prosite maps.

We will call such a dyad an $\infty$-*prosite* and say a map of $\infty$-prosites $f : \mathcal{C}^{\leq} \longrightarrow \mathcal{D}^{\leq}$ is a *map of $\infty$-prosites* if it is pointwise a prosite map. Moreover, we will abuse language and attribute to an $\infty$-prosite, the properties of the $\infty$-category associated with it.

**Proposition 1.3.2.** *Let $\mathcal{C}$ be a small connected $\infty$-category. Suppose that we have a finitely complete prosite $(\mathcal{C}_0^{\leq}, J_0)$. Then there exists an $\infty$-prosite $(\mathcal{C}^{\leq}, J)$ where the zeroth-prosite is $(\mathcal{C}_0^{\leq}, J_0)$.*

*Proof.* We begin by upgrading $\mathcal{C}$ to an $\infty$-prosite $\mathcal{C}^{\leq}$ through Proposition 1.2.5. We argue from induction. We start by declaring a sieve $\{\omega_i \longrightarrow \sigma\}_{i \in I}$ in $\mathcal{C}_1^{\leq}$ to be a covering sieve if and only if $\{\partial_k \omega_i \longrightarrow \partial_k \sigma\}_{i \in I}$ is a covering sieve for $k \in \{0, 1\}$. Invoking Proposition 1.2.17, we see that the collection of these sieves is a Grothendieck topology $J_1$ on $\mathcal{C}_1^{\leq}$. Furthermore, and again following from 1.2.17, all transition maps are left exact. This means that since these maps are a priori cover preserving (observing the identity $\partial_k s_0 = \text{id}$), they meet the criteria for prosite maps. Now taking the inductive step, we notice that in a manner similar to that argued, a prosite structure is obtained on $(\mathcal{C}_{n+1}^{\leq}, J_{n+1})$ from $(\mathcal{C}_n^{\leq}, J_n)$. It remains to show that transition maps are cover preserving. This is a priori true for face maps and is underwritten by simplicial identities for the degeneracy maps. □

**Proposition 1.3.3.** *Let $f : \mathcal{C} \longrightarrow \mathcal{D}$ be a map of small connected $\infty$-categories. Suppose that we have finitely complete prosites $(\mathcal{C}_0^{\leq}, J_0)$ and $(\mathcal{D}_0^{\leq}, J_0')$ and that $f_0 : (\mathcal{C}_0^{\leq}, J_0) \longrightarrow (\mathcal{D}_0^{\leq}, J_0')$ is a map of prosites. Then we have a map $f : (\mathcal{C}^{\leq}, J) \longrightarrow (\mathcal{D}^{\leq}, J')$ of the $\infty$-prosites determined from the zeroth-prosites.*

*Proof.* It suffices to show that for any $n \geq 0$, $f_n : (\mathcal{C}_n^{\leq}, J_n) \longrightarrow (\mathcal{D}_n^{\leq}, J_n')$ is a map of prosites. We will argue for the case where $n = 1$ and it will become apparent that the argument can be transplanted to the case where $n > 1$. Taking into account $f_0 \partial = \partial f_1$ and that $\partial$ are left exact and create finite limits jointly (cf. Proposition 1.2.17), we deduce that $f_1$ commutes with finite limits. Moreover, given how we have defined the topology, it again follows from $f_0 \partial = \partial f_1$ that $f_1$ preserves the covering sieves whenever $f_0$ preserves them. □

**Corollary 1.3.4.** *Let $f : \mathcal{C} \longrightarrow \mathcal{D}$ be a map of small connected $\infty$-categories. Suppose that we have finitely complete prosites $(\mathcal{C}_0^{\leq}, J_0)$ and $(\mathcal{D}_0^{\leq}, J_0')$ and that $f_0 : (\mathcal{C}_0^{\leq}, J_0) \longrightarrow (\mathcal{D}_0^{\leq}, J_0')$ is an equivalence of prosites. Then we have an equivalence $f : (\mathcal{C}^{\leq}, J) \longrightarrow (\mathcal{D}^{\leq}, J')$ of the $\infty$-prosites determined from the zeroth-prosites.*

*Proof.* By definition, an equivalence of prosites is an equivalence of prosites that is also a bijection of the covering sieves. Combining Proposition 1.3.3 and Proposition 1.2.19 and seeing how zeroth-prosites induce $\infty$-prosites, we obtain the result. □

It is now suitable for us to introduce an $\infty$-category of $\infty$-prosites. We will limit ourselves to when we have, as objects, $\infty$-prosites arising from connected $\infty$-prosets having pointwise, finitely complete prosites. Taking into account Corollary 1.3.4, this supposed $\infty$-category appears to be fully determined by a certain category having as objects prosites that have finitely complete underlying categories and as morphisms, morphisms of prosites. We would also like these categories to be the homotopy categories of small connected $\infty$-categories.

**Definition 1.3.5.** Consider the category whose objects are finitely complete prosites with underlying categories of the kind $(Obj(h\mathcal{C}), \leq)$, and whose morphisms are morphisms of prosites; where $\mathcal{C}$ is a small connected $\infty$-category. We will call the nerve of this category the $\infty$-category of finitely complete connected $\infty$-prosites. We will label it as such: $\beth_\infty$.

**Proposition 1.3.6.** *There exists a map of $\infty$-categories $q : \beth_\infty \longrightarrow \mathsf{Fun}(\mathsf{N}(\Delta), \mathsf{N}(\mathsf{Loc}))^{op}$.*

*Proof.* Theorem 1.1.2 informs us that given any prosite $(C, J)$, we can find a locale $Id_J(C)$ such that $\mathbf{Sh}(Id_J(C)) \simeq_{\mathsf{Topos}} \mathbf{Sh}(C, J)$. Now, given the reflective localization $\mathsf{Loc} \longrightarrow \mathsf{Topos}$, we obtain the assignment $(C, J) \longmapsto Id_J(C, J)$ which is part of a contravariant functor $\mathsf{Prosite} \longrightarrow \mathsf{Loc}$. Therefore, for every $\infty$-prosite we derive a cosimplicial locale $\Delta \longrightarrow \mathsf{Loc}$ and for every map of $\infty$-prosites, a map of cosimplicial locales in the opposite direction. Observing Proposition 1.3.3 in light of the definition given in 1.3.5 leads to the desired conclusion. Implicit in this is the use of Proposition 3.13 of the entry "Nerves" in [xvii] that establishes the categorical equivalence $\mathsf{N}(\mathsf{Hom}_{\mathsf{Cat}}(\Delta, \mathsf{Loc})^{\mathrm{op}}) \simeq \mathsf{Fun}(\mathsf{N}(\Delta), \mathsf{N}(\mathsf{Loc}))^{\mathrm{op}}$. □

**Proposition 1.3.7.** *Let $\beth_\infty^{sub}$ be the full subcategory of $\beth_\infty$ spanned by the $\infty$-prosites whose Grothendieck topologies are subcanonical. Then $q|\beth_\infty^{sub}$ is a monomorphism of $\infty$-categories.*

*Proof.* It is observed in Theorem 3.5 of [xi] that the map $\mathsf{Prosite} \longrightarrow \mathsf{Loc}$ is faithful up to naturally isomorphic functors. Therefore, $q|\beth_\infty^{sub}$ is faithful pointwise up to homotopy. □

Fix $\Sigma(\mathcal{C}) \subseteq \mathsf{Cat}$ as the category of subobjects of a finitely complete proset $\mathcal{C}$ and inclusions of these into each other. Consider a functor $\Gamma : \mathcal{C}^{\mathrm{op}} \longrightarrow \Sigma(\mathcal{C})$ that is characterized as follows.

1. $c \in \Gamma(c)$.
2. Given $d \in \Gamma(c)$, $\Gamma(d) \subseteq \Gamma(c)$.
3. $\Gamma(c)$ is finitely complete and the inclusion $\Gamma(c) \subseteq \mathcal{C}$ is left exact.

We will call such a functor a *packeting* of $\mathcal{C}$. For $a \in \Gamma(c)$, let $S_{a,c}^\Gamma : \mathcal{C}^{\mathrm{op}} \longrightarrow \mathsf{Set}$ be the sieve mapping each $d \in \mathcal{C}$ to $\mathsf{Hom}_\mathcal{C}(d,c)$ if $d \leq a$ and to the empty set otherwise. By abuse of notation, we will label the family of all arrows $d \longrightarrow c$ such that $S_{a,c}^\Gamma(d) \neq \varnothing$ as $S_{a,c}^\Gamma$ and the collection of all $\bigcup_{a \in \Gamma(c)} S_{a,c}^\Gamma$ for all objects of $\mathcal{C}$ as $S^\Gamma$.

**Proposition 1.3.8.** *Suppose that $\Gamma : \mathcal{C}^{\mathrm{op}} \longrightarrow \Sigma(\mathcal{C})$ is a packeting of a finitely complete proset $\mathcal{C}$. Then $S^\Gamma$ is a coverage.*

*Proof.* If $a \in \Gamma(c)$, then $g \in S_{a,c}^\Gamma$ if and only if it factors through $a \times c \longrightarrow c$. Given any $f : d \longrightarrow c$ in $\mathcal{C}$, we consider the collection of pullbacks of maps in $S_{a,c}^\Gamma$ along $f$, $f^*S_{a,c}^\Gamma$. We claim $f^*S_{a,c}^\Gamma = S_{\delta(a),d}^\Gamma$ where $\delta(a) = (a \times c) \times_c d$. By definition, $\Gamma(c) \subseteq \Gamma(d)$ and given the finite completeness of $\Gamma(d)$, $\delta(a) \in \Gamma(d)$. Now, since each map in $f^*S_{a,c}^\Gamma$ factors through $\delta(a) \longrightarrow d$, by default, $f^*S_{a,c}^\Gamma \subseteq S_{\delta(a),d}^\Gamma$. In the other direction, given $g \in S_{\delta(a),d}^\Gamma$, there exists $h : dom(g) \longrightarrow \delta(a)$ whose composition with the projection $\delta(a) \longrightarrow a \times c \longrightarrow c$ gives a map in $S_{a,c}^\Gamma$ whose pullback along $f$ recovers $g$. Thus $S_{\delta(a),d}^\Gamma \subseteq f^*S_{a,c}^\Gamma$. □

**Definition 1.3.9.** Let $(\mathcal{C}, J)$ be a finitely complete prosite such that the Grothendieck topology is determined from a packeting $\Gamma : \mathcal{C}^{\mathrm{op}} \longrightarrow \Sigma(\mathcal{C})$. We will say that such a prosite is *packeted* and similarly call an $\infty$-prosite determined from a *packeted* zeroth-prosite. Maps of prosites (resp. $\infty$-prosites) between *packeted* prosites (resp. *packeted* $\infty$-prosites) will be called maps of *packeted* prosites (resp. maps of *packeted* $\infty$-prosites). $\beth_\infty^\square$ will be the full subcategory of $\beth_\infty$ spanned by *packeted* $\infty$-prosites.

**Definition 1.3.9.** Let $(\mathcal{C}, J)$ be a finitely complete prosite such that the following are true.

1. The Grothendieck topology is determined from a packeting $\Gamma : \mathcal{C}^{\mathrm{op}} \longrightarrow \Sigma(\mathcal{C})$.
2. Whenever $V \subseteq U$ is an inclusion of $J$-ideals where $U$ is finite and $\langle V^c \rangle$ is the smallest $J$-ideal containing $U \backslash V$, if $\langle V^c \rangle \cup W = U$ for some $V \supseteq W \in Id_J(\mathcal{C})$, then $W = V$.

We will say that such a prosite is *retro-packeted* and similarly call an $\infty$-prosite determined from a *retro-packeted* zeroth-prosite: choice of the latter language for $\infty$-prosites will make much sense after Corollary 1.3.19. Maps of prosites (resp. $\infty$-prosites) between *retro-packeted* prosites (resp. *retro-packeted* $\infty$-prosites) will be called maps of *retro-packeted* prosites (resp. maps of *retro-packeted* $\infty$-prosites). $\beth_\infty^\sharp$ will be the full subcategory of $\beth_\infty$ spanned by *retro-packeted* $\infty$-prosites.

**Remark 1.3.10.** It seems natural to define a map of packeted prosites $f : (\mathcal{C}, J^{\Gamma_0}) \longrightarrow (\mathcal{D}, J^{\Gamma_1})$ as a map of prosites such that $f^{-1}\Gamma_1 = \Gamma_0$. However, in practice this is unnecessarily restrictive since prosite maps need only *generate* covering sieves; this amounts to the requirement that $f(\Gamma_0(c)) \subseteq \Gamma_1(f(c))$ for all $c \in \mathcal{C}$. Furthermore, in accordance with our needs, we find that the definition adapted for the said maps is adequate; this is exemplified in Lemma 1.3.12 and Lemma 1.3.13.

**Proposition 1.3.11.** *Suppose that $(\mathcal{C}, J^\Gamma)$ is a finitely complete prosite such that the Grothendieck topology is determined from a packeting $\Gamma : \mathcal{C}^{op} \longrightarrow \Sigma(\mathcal{C})$ where for all $c \in \mathcal{C}$, $d \in \Gamma(c)$ for any $d \leq c$. Then $(\mathcal{C}, J^\Gamma)$ is a retro-packeted prosite.*

*Proof.* It suffices to show that given any inclusion of $J^\Gamma$-ideals $V \subseteq U$, $V^c = \langle V^c \rangle$. For a subset $M$ of a proset $X$ if we denote by $M \downarrow$ the set $\{x \in X : x \leq m, m \in M\}$, then we note that for any prosite $\langle V^c \rangle = V^c \downarrow$. But in our case, if $r \in V$ and $r \leq s$ where $s \in V^c$, then $r \in \Gamma(s)$. This implies $s \in V$. □

**Lemma 1.3.12.** *Suppose that $(\mathcal{C}, J)$ is a packeted prosite. Then the functor $\mathsf{Prosite}^{op} \longrightarrow \mathsf{Loc}$ restricts on packeted prosites to a map with codomain $\mathsf{Loc}^{sob}$.*

*Proof.* We show that $(\mathcal{C}, J)$ maps to a spatial locale. Through II Definition 1.5 of [xxi], we are made aware that a locale $L$ is spatial if and only if for any opens $U$ and $V$ of $L$, $U \leq V$ if and only every point of $L$ that belongs to $U$ belongs to $V$. In the framework of prosites and the resultant locale of Theorem 1.1.2, this amounts to showing that whenever it is the case that any $J$-prime filter intersecting $U$ must intersect $V$, $U \subseteq V$ (cf. Proof of Proposition 2.7 of [xi]). Suppose that for any object $a \in \mathcal{C}$, $\Gamma(a)$ is a prime filter and that $a \in U$. By the hypothesis above, $\Gamma(a) \cap V \neq \varnothing$. Assume $b \in \Gamma(a) \cap V$, then all $c \in \mathcal{C}$ such that $c \leq a \times b$ are contained in $V$. This describes a covering sieve $\{f_i : c_i \longrightarrow a\}$.

$\Gamma(a)$ admits finite limits. This means it is non-empty and that it is closed under taking meets. Furthermore, by the contravariance of $\Gamma$, it is upward closed. It remains to observe that given any covering sieve $\{f_i : c_i \longrightarrow b\}_{i \in I}$ for any $b \in \Gamma(a)$, we can find $c_i \in \Gamma(a)$ for some $i \in I$. But this is guaranteed by the fact that if $b \in \Gamma(a)$, $\Gamma(b) \subseteq \Gamma(a)$.

Now observe that maps of spatial locales are simply maps of locales. That is, we have a fully faithful inclusion $\mathsf{Loc}^{sob} \hookrightarrow \mathsf{Loc}$. □

**Lemma 1.3.13.** *Suppose that $(\mathcal{C}, J)$ is a retro-packeted prosite. Then the functor $\mathsf{Prosite}^{op} \longrightarrow \mathsf{Loc}$ maps $(\mathcal{C}, J)$ to a coherent locale.*

*Proof.* Theorem 1.1.2 establishes that the locale $Id_J(C)$ is has the frame whose objects are $J$-ideals and whose morphisms are subset inclusions. Immediately, we have that any open $V \in Id_J(C)$ is determined as the minimal $J$-ideal containing $\bigcup_{c \in V}(c) \downarrow_J$. Furthermore, given that $\mathcal{C}$ admits a final object $*$, $(*) \downarrow_J = obj(\mathcal{C})$. Hence, $(*) \downarrow_J = 1$, the top element of $Id_J(C)$. Subsequently, it suffices to show that principal ideals are finite, a union of $J$-ideals is a $J$-ideal, and that the intersection of two finite $J$-ideals is again a finite $J$-ideal.

For some set $\Omega$, let $I = \bigcup_{\alpha \in \Omega} I_\alpha$ where $I_\alpha \in Id_J(\mathcal{C})$. If $c \in I$, then $c \in I_\alpha$ for some $\alpha \in \Omega$. Therefore, whenever $b \leq c$, $b \in I_\alpha$ and hence $b \in I$. If $R$ is some covering sieve of $c$, suppose that for all $f \in R$ we have $dom(f) \in I$. Keeping in mind that the Grothendieck topology of $\mathcal{C}$ is induced by a packeting, there exists $f : b \longrightarrow c$ belonging to some $S^\Gamma_{a,c}$ such that $f^*R = S^\Gamma_{a',b}$. This means that for every $h \in f^*R$, $dom(h) \in I$ since $dom(h) \leq dom(h')$ for some $h' \in R$ such that $f^*h' = h$. Specifically, when we take $h : a' \times b \longrightarrow b$, that is $dom(h) = a' \times b$, this indicates that $b \in I_\alpha$ for some $\alpha \in \Omega$. In fact, it is enough to recall the instance $b = a \times c$ so that the domains of all $g \in S^\Gamma_{a,c}$ are in $I_\alpha$. But this implies $c \in I_\alpha \Rightarrow c \in I$. Therefore, $I \in Id_J(\mathcal{C})$.

To show that $(c) \downarrow_J$ is finite, it suffices to show that if $(c) \downarrow = I \vee J$, then $c \in I$ or $c \in J$. Notice that the conclusion of the previous paragraph tells us that $I \vee J = I \cup J$. On the other hand, assuming that $I$ and $I'$ are finite elements of $Id_J(\mathcal{C})$, suppose that $I \cap I' = \bigvee_{\alpha \in \Omega} W_\alpha$ where each $W_\alpha \in Id_J(\mathcal{C})$. Let $\langle I^c \rangle$ be the smallest $J$-ideal containing the complement of $I \cap I'$ relative to $I$ and set $F_\alpha = \langle I^c \rangle \cup I_\alpha$. We know thanks to the previous paragraph that $F_\alpha$ is a $J$-ideal. Thus, because $I = \bigvee_{\alpha \in \Omega} F_\alpha$ and $I$ is finite, there exists a finite subset $S \subseteq \Omega$ such that $I = \bigvee_{\alpha \in S} F_\alpha$. But this means $I = \langle I^c \rangle \bigvee_{\alpha \in S} I_\alpha$. Recalling the definition of retro-packeted prosites, we deduce that $I \cap I' = \bigvee_{\alpha \in S} I_\alpha$. $\square$

**Lemma 1.3.14.** *Suppose that $f : (\mathcal{C}, J) \longrightarrow (\mathcal{D}, J')$ is a map of retro-packeted prosites. Then the functor $\mathsf{Prosite}^{op} \longrightarrow \mathsf{Loc}$ maps $f$ to a map of coherent locales.*

*Proof.* From the remarks after the proof of Proposition 3.2 in [xi], we understand that the map of frames $Id_J(\mathcal{C}) \longrightarrow Id_{J'}(\mathcal{D})$ corresponding to the locale map that is the image of $f$ under the functor $\mathsf{Prosite}^{op} \longrightarrow \mathsf{Loc}$ acts as follows: it maps $I \in Id_J(\mathcal{C})$ to the smallest $J'$-ideal containing $f(I)$, $\langle f(I) \rangle$. Therefore, we need to show that when $I$ is a finite element, $\langle f(I) \rangle$ is also a finite element. Suppose $\langle f(I) \rangle = \bigvee_{\alpha \in \Omega} J_\alpha$. Recall that when working with retro-packeted prosites, joins are the same as unions. Furthermore, generally, the preimage of an ideal of a prosite is an ideal of the codomain prosite; that images of covering sieves generate covering sieves, which is a criterion of prosite maps, guarantees this. Resuming our case, consider $I = \bigvee_{\alpha \in \Omega} \widetilde{J}_\alpha$ where $\widetilde{J}_\alpha = f^{-1}(J_\alpha) \cap I$. By the finite condition, we find a finite $S \subseteq \Omega$ such that $I = \bigvee_{\alpha \in S} \widetilde{J}_\alpha$. Hence $f(I) = \bigcup_{\alpha \in \Omega} f(\widetilde{J}_\alpha)$ and this implies $\langle f(I) \rangle = \bigvee_{\alpha \in S} J_\alpha$. $\square$

**Proposition 1.3.15.** *The category of packeted prosites is finitely complete.*

*Proof.* This is the first paragraph of the proof of Proposition 1.3.18 after observing that the final object of Prosite is packeted; that is, it can only admit the packeting that is constant on itself since a packeting cannot take the empty set as a value by the finite completeness requirement. $\square$

**Corollary 1.3.16.** *If $(\mathcal{C}^\leq, J)$ is an $\infty$-prosite determined from a zeroth-prosite that is packeted, then each $(\mathcal{C}^\leq_n, J_n)$ is a packeted prosite for all $n \geq 0$.*

*Proof.* Proposition 1.3.3 through Proposition 1.2.17 makes us aware that each $(\mathcal{C}_{\overline{n}}^{\leq}, J_n)$ arises as a finite limit of a diagram of packeted prosites. Taking into account Proposition 1.3.15 gives the result. $\square$

**Theorem 1.3.17.** *There exists a map of $\infty$-categories $q^{\square} : \beth_{\infty}^{\square} \longrightarrow \mathsf{Fun}(\mathsf{N}(\Delta), \mathsf{N}(\mathsf{Loc}^{\mathrm{sob}}))^{\mathrm{op}}$.*

*Proof.* Combine Lemma 1.3.12 and Corollary 1.3.16. $\square$

**Proposition 1.3.18.** *The category of retro-packeted prosites is finitely complete relative to diagrams with surjective maps.*

*Proof.* It is immediate that the proset with one object and the trivial Grothendieck topology is retro-packeted. It remains to show that pullbacks exist. Considering the categories involved, the pullbacks of prosets are computed at the level of sets and therefore exist. Moreover, because the maps involved are left exact, the pullbacks obtained are also finitely complete; limits are computed pointwise. Therefore, if we are given maps of prosites $f_i : (\mathcal{D}_i, J^{\Gamma_i}) \longrightarrow (\mathcal{C}, J^{\Gamma})$ for $i \in \{0, 1\}$, the functor $\Gamma_{12}^{\mathrm{op}} : \mathcal{D}_0 \times_{\mathcal{C}} \mathcal{D}_1 \longrightarrow \Sigma^{\mathrm{op}}(\mathcal{D}_0 \times_{\mathcal{C}} \mathcal{D}_1)$ that acts as $(d_0, d_1) \longmapsto \Gamma_0(d_0) \times_{\Gamma(c)} \Gamma_1(d_1)$ where $c = f_i(d_i)$ meets the requirements for a packeting. Ergo, this exhibits half the criteria for the pullback being a retro-packeted prosite.

Suppose that $U$ is a finite element of the locale determined by the ideals of the prosite induced by the above packeting. Suppose also that we are given an inclusion $V \subseteq U$ and an equivalence $\langle V^c \rangle \cup W = U$ where $W \subseteq V$ is an ideal of the aforementioned prosite. Moreover, let $p_i : \mathcal{D}_0 \times_{\mathcal{C}} \mathcal{D}_1 \longrightarrow \mathcal{D}_i$ be projections from the pullback. For the reason that taking unions is the same as taking joins in our setup, in order to show $W = V$, it suffices to show that $\langle p_i(\langle V^c \rangle) \rangle = \langle \langle p_i(V) \rangle^c \rangle$ where the right-hand side is the complement relative to $\langle p_i(U) \rangle$. This is because from Lemma 1.3.14 we are made aware that projections, being left exact, preserve finite elements, and thus $\langle p_i(U) \rangle$ is finite; taking into account $\langle p_i(\langle V^c \rangle) \rangle = \langle \langle p_i(V) \rangle^c \rangle$, the former would in turn imply $\langle p_i(W) \rangle = \langle p_i(V) \rangle$ given that each $\mathcal{D}_i$ is a retro-packeted prosite. And since in the category Set pullbacks of epimorphisms are always epimorphisms, the projections $p_i$ are surjective, and hence, taking into account $W \subseteq V$, we would then arrive at the desired conclusion. We observe $p_i(\langle V^c \rangle) \downarrow = (p_i(U) \backslash p_i(V)) \downarrow$. Further, we note that adding to $p_i(\langle V^c \rangle) \downarrow$ (resp. $(p_i(U) \backslash p_i(V)) \downarrow$) any $d \in \mathcal{D}_i$ which has some covering sieve $R$ such that $dom(f) \in p_i(\langle V^c \rangle) \downarrow$ (resp. $dom(f) \in (p_i(U) \backslash p_i(V)) \downarrow$) for all $f \in R$, gives $\langle p_i(\langle V^c \rangle) \rangle$ (resp. $\langle \langle p_i(V) \rangle^c \rangle$). $\square$

**Corollary 1.3.19.** *If $(\mathcal{C}^{\leq}, J)$ is an $\infty$-prosite determined from a zeroth-prosite that is retro-packeted, then each $(\mathcal{C}_{\overline{n}}^{\leq}, J_n)$ is a retro-packeted prosite for all $n \geq 0$.*

*Proof.* Proposition 1.3.3 through Proposition 1.2.17 makes us aware that each $(\mathcal{C}_{\overline{n}}^{\leq}, J_n)$ arises as a finite limit of a diagram of retro-packeted prosites where the maps involved are surjective. Taking into account Proposition 1.3.18 gives the result. $\square$

**Theorem 1.3.20.** *There exists a map of $\infty$-categories $q^{\sharp} : \beth_{\infty}^{\sharp} \longrightarrow \mathsf{Fun}(\mathsf{N}(\Delta), \mathsf{N}(\mathsf{Frm}^{\mathrm{coh}})^{\mathrm{op}})^{\mathrm{op}}$.*

*Proof.* Combine Lemma 1.3.13 and Lemma 1.3.14 with Corollary 1.3.19. □

We hope to establish a couple of functors $\mathsf{Frm}^{\mathrm{op}} \longrightarrow [\mathrm{Poly}_S, \mathrm{Set}]$. Here, $\mathrm{Poly}_S$ is the full subcategory of commutative rings spanned by rings of the kind $S[T_0, \ldots, T_n]$ where $n$ is finite and $S$ is a commutative ring. We will then exploit these functors in combination with Proposition 1.3.6 and Theorem 1.3.20 to obtain functors of $\infty$-categories into animated $S$-stacks.

The functor $\mathrm{Spec} : \mathrm{Poly}_S^{\mathrm{op}} \longrightarrow (\mathsf{Frm}^{\mathrm{coh}})^{\mathrm{op}}$ provides the basis for the first of these functors. Notice that Theorem 1.1.3 combined with the equivalence $\mathsf{Loc}^{\mathrm{coh}} \longrightarrow \mathcal{Sp}$ gives us this functor. Furthermore, we are informed that it is full and essentially surjective. Keeping in mind the functor now obtained, the inclusion $i^{\mathrm{op}} : (\mathsf{Frm}^{\mathrm{coh}})^{\mathrm{op}} \subseteq \mathsf{Frm}^{\mathrm{op}}$ and the coYoneda embedding $\sharp^{\mathrm{op}} : \mathrm{CRing}^{\mathrm{op}} \hookrightarrow [\mathrm{CRing}, \mathrm{Set}]$, induce a functor $\mu_{\mathrm{Hoc}}^0(S) : \mathsf{Frm}^{\mathrm{op}} \longrightarrow [\mathrm{Poly}_S, \mathrm{Set}]$. It is the sequence

$$\mathsf{Frm}^{\mathrm{op}} \longrightarrow [\mathsf{Frm}, \mathrm{Set}] \xrightarrow{(i^{\mathrm{op}})^*} [\mathsf{Frm}^{\mathrm{coh}}, \mathrm{Set}] \xrightarrow{(\mathrm{Spec}^{\mathrm{op}})^*} [\mathrm{Poly}_S, \mathrm{Set}]$$

On the other hand, through the forgetful functor $\mathsf{Frm} \longrightarrow \mathsf{DLat}$, Stone duality, and the functor $\mathsf{Top}^{\mathrm{op}} \longrightarrow \mathsf{BoolRing}$, we obtain a functor $\mathsf{Frm} \longrightarrow \mathsf{BoolRing}$. Again, through the maps induced by the coYoneda embedding, we arrive at a functor $\mathrm{Sm}_0^{\flat}(S) : \mathsf{Frm}^{\mathrm{op}} \longrightarrow [\mathrm{Poly}_S, \mathrm{Set}]$. Given the inclusion $j^{\mathrm{op}} : \mathrm{Poly}_S^{\mathrm{op}} \subseteq \mathrm{CRing}^{\mathrm{op}}$, the derivation is as follows

$$\mathsf{Frm}^{\mathrm{op}} \longrightarrow \mathsf{BoolRing}^{\mathrm{op}} \xrightarrow{\sharp} [\mathrm{CRing}, \mathrm{Set}] \xrightarrow{(j^{\mathrm{op}})^*} [\mathrm{Poly}_S, \mathrm{Set}]$$

**Corollary 1.3.21.** $\mu_{\mathrm{Hoc}}^0(S)$ *is faithful when restricted to the categories of locales of open subsets corresponding to the following subcategories of* $\mathrm{Spec}(\mathrm{Poly}_S^{\mathrm{op}})$:

1. *The subcategory of all spectral spaces and surjective spectral maps.*
2. *For a spectral space $X$, the subcategory of its spectral subspaces and inclusions of these.*
3. *The full subcategory whose objects are the $T_1$ spectral spaces.*

*Proof.* This is immediate from Theorem 1.1.3. □

**Proposition 1.3.22.** $\mu_{\mathrm{Hoc}}^0(S)$ *and* $\mathrm{Sm}_0^{\flat}(S)$ *are coincident when restricted to the category of locales of open subsets corresponding to the full subcategory of compact and totally disconnected Hausdorff spaces in* $\mathrm{Spec}(\mathrm{Poly}_S^{\mathrm{op}})$.

*Proof.* Stone duality (for Stone spaces) establishes an equivalence between the category of frames of open subsets corresponding to Stone spaces and $\mathsf{BoolRing}$. Consequently, and since spectral spaces are sober, $\mathrm{Sm}_0^{\flat}$ is fully faithful when restricted as in the claim. On the other hand, Proposition 1.3.21 (3) establishes that $\mu_{\mathrm{Hoc}}^0(S)$ is faithful on a certain full subcategory of locales of open subsets corresponding to $\mathrm{Spec}(\mathrm{Poly}_S^{\mathrm{op}})$. But those objects of $\mathrm{Spec}(\mathrm{Poly}_S^{\mathrm{op}})$ that are $T_1$ spaces are exactly Stone spaces. □

The functors from $\mathsf{Frm}^{\mathrm{op}}$ to $[\mathrm{Poly}_S, \mathrm{Set}]$ previously determined act pointwise to substantiate a dyad of functors $\mathsf{Fun}(\mathsf{N}(\Delta), \mathsf{N}(\mathsf{Frm})^{\mathrm{op}})^{\mathrm{op}} \longrightarrow \mathsf{Fun}(\mathsf{N}(\Delta), \mathsf{Stk}_S)^{\mathrm{op}}$. Furthermore, if we label the geometric corealization functor as $\mathrm{colim}_\Delta : \mathsf{Fun}(\mathsf{N}(\Delta), \mathsf{Stk}_S) \longrightarrow \mathsf{Stk}_S$ and recall the map of Proposition 1.3.6, we should obtain maps $\beth_\infty \longrightarrow \mathsf{Stk}_S^{\mathrm{op}}$.

**Definition 1.3.23.** Let $\widetilde{\mu}_{\text{Hoc}}(S) : \text{Fun}(\mathsf{N}(\Delta), \mathsf{N}(\text{Frm})^{\text{op}})^{\text{op}} \longrightarrow \text{Fun}(\mathsf{N}(\Delta), \text{Stk}_S)^{\text{op}}$ be the map obtained from $\mu_{\text{Hoc}}^0$ acting pointwise on simplicial frames. Then we have a map $\mu_{\text{Hoc}}(S) : \beth_\infty \longrightarrow \text{Stk}_S^{\text{op}}$ defined as $\text{colim}_\Delta^{\text{op}} \circ \widetilde{\mu}_{\text{Hoc}}(S) \circ q$. We will call this map the *S-Hochster spectrum*. We will drop the '$S$' when $S = \mathbb{Z}$.

**Definition 1.3.24.** Let $\widetilde{\text{Sm}}^\flat(S) : \text{Fun}(\mathsf{N}(\Delta), \mathsf{N}(\text{Frm})^{\text{op}})^{\text{op}} \longrightarrow \text{Fun}(\mathsf{N}(\Delta), \text{Stk}_S)^{\text{op}}$ be the map obtained from $\text{Sm}_0^\flat$ acting pointwise on simplicial frames. Then we have a map $\text{Sm}^\flat(S) : \beth_\infty \longrightarrow \text{Stk}_S^{\text{op}}$ defined as $\text{colim}_\Delta^{\text{op}} \circ \widetilde{\text{Sm}}^\flat(S) \circ q$. We will call this map the *S-Boolean smashing spectrum* and drop the '$S$' when $S = \mathbb{Z}$. The naming here is informed by the smashing spectrum of condensed mathematics which is, informally speaking, determined by objects of symmetric monoidal $\infty$-categories such that $x^2 = x$: this is exactly how elements of a Boolean ring behave.

**Corollary 1.3.25.** *Let $\beth_\infty^0$ be the full subcategory of $\beth_\infty$ that maps under $q$ to cosimplicial Stone locales. Then $\text{Sm}^\flat(S)$ and $\mu_{\text{Hoc}}(S)$ are coincident when restricted to $\beth_\infty^0$.*

*Proof.* This follows from Proposition 1.3.22. □

**Proposition 1.3.26.** *Suppose that $\text{Aff}_S$ is the $\infty$-category dual to the $\infty$-category of simplicial commutative S-algebras and $X \in \beth_\infty$. Then setting $\mu(S) = \text{colim}_\Delta^{\text{op}} \circ \widetilde{\mu}(S) \circ q$ as a stand-in for both $Sm^\flat(S)$ and $\mu_{Hoc}(S)$, the following are equivalent.*

1. $\mu(X) \in \text{Aff}_S$.
2. $(\widetilde{\mu}(S) \circ q)_0(X) : \mathsf{N}(\Delta) \longrightarrow \text{Pro}(\mathsf{N}(\text{Poly}_S)^{\text{op}})$.

*Proof.* This is a direct application of Lemma 5.5.8.14 in [i]. □

# Geometric Factorization of Six Functor Formalisms

In this section, we prove Theorem A and Theorem B. The outset of this endeavor is the act of demystifying the notion of an *inverse Tannakian formalism*. So far, we have seen that in the definition offered, it is transparent why this notion encapsulates a "geometric space": being that, roughly speaking, factorization through classical stacks is involved. We would like to begin at the place where the latter is a priori unknown. Here, the immediate task is finding out the material conditions under which a map into an appropriate category of algebraic objects factors through stacks. Theorem A is, of course, an answer to this question. And, furthermore, Theorem B, taken in conjunction with Theorem A, makes concrete, albeit in a limited sense, what is alluded to in occurrences such as Example 0.2.2 and the discussion following Theorem 1.6 in [iii]. That is, it expresses some 6-functor formalisms as $\infty$-algebras in an appropriate $\infty$-topos; as discussed earlier under "What is a Space?".

   What is implicit in the above theorems, and eventually what arises as an object central to our concern, is a *suprematic space*. We will begin by proving results with an eye towards formulating these objects; which after formulating, we will supply "real-world" examples. It is worth mentioning that the study of *suprematic spaces* is founded upon the work of J. Lurie in [iv] and [i]. And when time demands, we will only minimally reintroduce the necessary tenets and outsource an in-depth exposition to a review of the primary sources.

## Suprematic Spaces

[1] We will assume that all $\infty$-topoi live in the universe of very large $\infty$-categories and we will use *small* to mean $\mathbb{U}_1$-small. The main aim of this subsection is to give a proof of Theorem A. When visualizing diagrams of the kind $\Delta^n \times \Delta^n \longrightarrow \mathcal{C}$, we adapt the convention in which the first factor is the horizontal direction and the second factor is vertical direction.

**Definition 2.1.1.** Consider a map $q : \mathcal{D} \longrightarrow \mathcal{C}$ of $\infty$-categories where $\mathcal{C}$ has pullbacks. Let $P_q$ be the collection of all the squares $\Delta^1 \times \Delta^1 \longrightarrow \mathcal{D}$ that map to pullback squares when composed with $q$ and that meet the following criteria:

1. Suppose that $\sigma : \Delta^2 \times \Delta^1 \longrightarrow \mathcal{D}$ is depicted as the diagram

$$
\begin{array}{ccccc}
\bullet & \longrightarrow & \bullet & \xrightarrow{f} & \bullet \\
\downarrow & & \downarrow & & \downarrow \\
\bullet & \longrightarrow & \bullet & \longrightarrow & \bullet
\end{array}
$$

    where $f$ is degenerate. If both the outer square and the right square are in $P_q$, the left square is in $P_q$. If the left square is in $P_q$ and the edge opposite to $f$ spans a square in $P_q$, then the right square is in $P_q$.

2. Given $\sigma : \Delta^2 \times \Delta^1 \longrightarrow \mathcal{D}$, if $\sigma|\Delta^{\{i,i+1\}} \times \Delta^1$ are in $P$ for $0 \leq i \leq 1$, then $\sigma|\Delta^{\{0,2\}} \times \Delta^1$ is in $P_q$. Informally speaking, this means that $P_q$ is closed under "composition/pasting".

3. The restriction of $q$ to the subcategory spanned by edges of squares in $P_q$ is (-1)-truncated.

We drop the subscript from $P_q$ where the context allows us to do so and to simplify our notation.

**Definition 2.1.2.** Given $P$ as above and under the same hypothesis, we define $\xi_P(\mathcal{C}; q)$ as the collection of morphisms in $\mathcal{C}$ with the following properties.

1. $\mathrm{id}_C \in \xi_P(\mathcal{C}; q)$ for all $C \in q(\mathcal{D})$.

2. Given $f \in \xi_P(\mathcal{C}; q)$ and any $\alpha : \Lambda_2^2 \longrightarrow \mathcal{D}$ with $q(\partial_0 \alpha) = f$, there exists $\widetilde{\alpha} : (\Lambda_2^2)^\triangleleft \longrightarrow \mathcal{D}$ such that $\widetilde{\alpha} \in P$.

3. Given $f \in \xi_P(\mathcal{C}; q)$ and the rectangle $\Delta^2 \times \Delta^1 \longrightarrow \mathcal{D}$

$$
\begin{array}{ccccc}
\bullet & \longrightarrow & \bullet & \longrightarrow & \bullet \\
\downarrow & & \downarrow & & \downarrow {\scriptstyle f'} \\
\bullet & \longrightarrow & \bullet & \longrightarrow & \bullet
\end{array}
$$

    where the right square is in $P$ and $q(f') = f$, the left square is in $P$ if the outer square is in $P$.

---

[1]The naming of these spaces is inspired by the twentieth century art movement, *suprematism*. Here, feeling is taken as the supreme quality of painting; in practice, the stripping of painterly forms to their bare minimum is insisted. And what do we know—geometric shapes are found to be the building blocks of said forms! The analogy is carried over in the embrace of pointless topology as the study of spaces, geometric or otherwise, par excellence. It is in this manner, bundled with the 'motivic shape' of Theorem C, that we adapt the name.

When we do not suffer the chance of imprecision, we will simply write $\xi_P(\mathcal{C}; q)$ as $\xi_P$. Note that in general, the existence of $\xi_P(\mathcal{C}; q)$ is non-trivial.

**Proposition 2.1.3.** $\xi_P(\mathcal{C}; q)$ *is stable under base change along any* $f \in q(\mathcal{D})$.

*Proof.* Suppose that we have the following pullback square in $\mathcal{C}$ where $h \in \xi_P$.

$$\begin{array}{ccc} \cdot & \xrightarrow{f'} & \cdot \\ {\scriptstyle h'}\downarrow & & \downarrow{\scriptstyle h} \\ \cdot & \xrightarrow{f} & \cdot \end{array}$$

Then it follows that there exists an $\alpha : \Lambda_2^2 \longrightarrow \mathcal{D}$ with $q(\partial_0 \alpha) = h$ and $q(\partial_1 \alpha) = f$. This means that we have a square in $P$ whose image under $q$ is the square above, given property (1) of $\xi_P$; call this square $Q$ and observe that it is unique up to contractible choice by $q$ being (-1)-truncated when restricted as specified by condition (3) of $P$. Given any $\alpha : \Lambda_2^2 \longrightarrow \mathcal{D}$ with $q(\partial_0 \alpha) = h'$, consider instead $\widetilde{\alpha} : \Lambda_2^2 \longrightarrow \mathcal{D}$ with $q(\partial_0 \widetilde{\alpha}) = h$ and $q(\partial_1 \widetilde{\alpha}) = f \circ q(\partial_1 \alpha)$. By property (2) of $\xi_P$ and the fact that squares in $P$ map to pullback squares, we have a square in $\mathcal{C}$ where the vertical arrow on the right is $h$ and the horizontal map on the bottom is $f \circ q(\partial_1 \alpha)$. But by the pasting law for pullbacks (cf. [i] Lemma 4.4.2.1), we have a pasting diagram $\Delta^2 \times \Delta^1 \longrightarrow \mathcal{D}$ where the right square is $Q$ and the left square is some $\lambda : (\Lambda_2^2)^\triangleleft \longrightarrow \mathcal{D}$ such that $\lambda | \Lambda_2^2 = \alpha$. Property (3) of $\xi_P$ then tells us that the right square is in $P$.

On the other hand, given a rectangle $\sigma : \Delta^2 \times \Delta^1 \longrightarrow \mathcal{D}$ where the right square is in $P$ and $q(u) = h'$,

$$\begin{array}{ccccc} \cdot & \longrightarrow & \cdot & \longrightarrow & \cdot \\ \downarrow & & \downarrow & & \downarrow{\scriptstyle u} \\ \cdot & \longrightarrow & \cdot & \longrightarrow & \cdot \end{array}$$

consider instead the pasting diagram whose right square is $Q$ and whose left square is the right square of $\sigma$; observe that the implicit assumption fixing $u$ as the left vertical map of $Q$ does not impede the generality sought. It follows from property (3) of $\xi_P$ that the outer square of this pasting diagram is in $P$; call it $Q'$. We now obtain a new rectangle whose right square is $Q'$ and whose left square is the left square of $\sigma$. Property (3) of $\xi_P$ then obtains the conclusion sought. □

**Proposition 2.1.4.** *Suppose that* $\xi_P(\mathcal{C}; q)$ *is stable under retracts in* $\mathsf{Fun}(\Delta^1, q(\mathcal{D}))$. *Then* $\xi_P(\mathcal{C}; q)$ *contains all equivalences in* $q(\mathcal{D})$.

*Proof.* An equivalence between $\partial_1 f = x$ and $\partial_1 g = y$ in $q(\mathcal{D})$ is made up of two pieces. That is, the piece $\sigma : \Delta^2 \longrightarrow q(\mathcal{D})$ such that $\partial_2 \sigma = f$, $\partial_0 \sigma = g$ and $\partial_1 \sigma = \mathrm{id}_x$ and the piece dual to the latter. When taken together, these data give the following rectangle.

$$\begin{array}{ccccc} x & \xrightarrow{\mathrm{id}_x} & x & \xrightarrow{\mathrm{id}_x} & x \\ {\scriptstyle f}\downarrow & & \downarrow{\scriptstyle \mathrm{id}_x} & & \downarrow{\scriptstyle f} \\ y & \xrightarrow{g} & x & \xrightarrow{f} & y \end{array}$$

Now recall property (1) in 2.1.2. □

**Proposition 2.1.5.** *Suppose we are given the following 2-simplex in $q(\mathcal{D})$ where $g \in \xi_P(\mathcal{C}; q)$.*

$$\begin{array}{ccc} & B & \\ {}^{f}\nearrow & & \searrow^{g} \\ A & \xrightarrow{g \circ f} & C \end{array}$$

*Then $g \circ f \in \xi_P(\mathcal{C}; q)$ if and only if $f \in \xi_P(\mathcal{C}; q)$.*

*Proof.* In one direction, assume that $g \circ f \in \xi_P$. Suppose $\alpha : \Lambda_2^2 \longrightarrow \mathcal{D}$ is some outer horn with $q(\partial_0 \alpha) = f$ and $\widetilde{\alpha} : \Lambda_2^2 \longrightarrow \mathcal{D}$ is the outer horn with $q(\partial_0 \widetilde{\alpha}) = g$ and $q(\partial_1 \widetilde{\alpha}) = g \circ q(\partial_1 \alpha)$. By $g \in \in \xi_P$ we obtain a square in $\mathcal{D}$ whose image under $q$ is the pullback diagram of $g$ along $q(\partial_1 \widetilde{\alpha})$. Furthermore, by $g \circ f \in \xi_P$, obtain another square in $P$ whose image is the pullback of $g \circ f$ along $q(\partial_1 \widetilde{\alpha})$. In short, by second half of property (1) of squares in $P$, we obtain a pasting diagram in $\mathcal{D}$

$$\begin{array}{ccccc} \bullet & \longrightarrow & \bullet & \xrightarrow{\ell} & \bullet \\ \downarrow & & \downarrow^{\partial_1 \alpha} & & \downarrow \\ \bullet & \xrightarrow{\partial_0 \alpha} & \bullet & \longrightarrow & \bullet \end{array}$$

where the outer square is in $P$, the right square is in $P$ and $\ell$ is degenerate. Therefore, from Definition 2.1.1, we conclude that the left square is in $P$.

To see that $f$ obeys property (3) of $\xi_P$, consider any pasting diagram $\lambda : \Delta^2 \times \Delta^1 \longrightarrow \mathcal{D}$ where the image of the rightmost vertical arrow is $f$ and the right square is in $P$, and the outer square is in $P$. Then we observe the emergence of a map $\sigma : \Delta^2 \times \Delta^2 \longrightarrow \mathcal{D}$ such that $\sigma | \Delta^2 \times \Delta^{\{0,1\}} = \lambda$. In fact, $q(\sigma)$ is depicted as the diagram

$$\begin{array}{ccccc} \bullet & \longrightarrow & \bullet & \longrightarrow & \bullet \\ \downarrow & & \downarrow & & \downarrow^{f} \\ \bullet & \xrightarrow{l} & \bullet & \xrightarrow{k} & \bullet \\ \downarrow^{\mathrm{id}_{\partial_1(l)}} & & \downarrow^{\mathrm{id}_{\partial_1(k)}} & & \downarrow^{g} \\ \bullet & \xrightarrow{l} & \bullet & \xrightarrow{g \circ k} & \bullet \end{array}$$

The preimage of the bottom right square lies in $P$ since $g \in \xi_P$ by second half of property (1) of squares in $P$. Thus, Definition 2.1.1 property (2) informs us that the outer square of $q(\sigma)|\Delta^{\{1,2\}} \times \Delta^2$ (rectangle composed of top and bottom right squares) is in $P$. By means of a similar argument, observing that $\mathrm{id}_{\partial_1(k)} \in \xi_P$, we deduce that the preimage of the outer square of $q(\sigma)|\Delta^2 \times \Delta^{\{1,2\}}$ is in $P$. Hence, since the preimage of $q(\sigma)|\Delta^2 \times \Delta^{\{0,1\}}$ is $\lambda$ and its outer square is in $P$, recalling Definition 2.1.1 property (2), we find that the preimage of the outer square of $q(\sigma)$ is also in $P$. Now because $g \circ f \in \xi_P$, we conclude that the preimage of $\sigma|\Delta^2 \times \Delta^{\{0,1\}}$ is in $P$. Consequently, since $q$ is a monomorphism in this context, it remains to note Definition 2.1.1 property (1).

In the other direction, assume that $f \in \xi_P$. Suppose $\alpha : \Lambda_2^2 \longrightarrow \mathcal{D}$ is some outer horn with $q(\partial_0 \alpha) = g \circ f$. Then utilizing the pasting law for pullbacks, we obtain a pasting diagram $\lambda : \Delta^2 \times \Delta^1 \longrightarrow \mathcal{D}$ where the image of the right square is the pullback square of $g$ along $q(\partial_1 \alpha)$ and the left square is the cartesian square of the pullback of the pullback of $q(\partial_1 \alpha)$ along $g$, along $f$. It follows that both the left and right squares of $\sigma$ are in $P$. Ergo, by characterization of $P$, the outer square of $\sigma$ is also in $P$.

Suppose we are given a rectangle $\lambda : \Delta^2 \times \Delta^1 \longrightarrow \mathcal{D}$ with the rightmost vertical arrow mapping to $g \circ f$ and where both the outer and right square are in $P$. Taking into account that under $q$ this maps to a pasting diagram of pullback squares, and given the pasting law for pullbacks, we obtain a map $\widetilde{\sigma} : \Delta^2 \times \Delta^2 \longrightarrow \mathcal{C}$ that factors through $q$. We can display this as the following diagram:

$$\begin{array}{ccccc}
\bullet & \xrightarrow{l''} & \bullet & \xrightarrow{k''} & \bullet \\
\downarrow{\scriptstyle f''} & & \downarrow{\scriptstyle f'} & & \downarrow{\scriptstyle f} \\
\bullet & \xrightarrow{l'} & \bullet & \xrightarrow{k'} & \bullet \\
\downarrow{\scriptstyle g''} & & \downarrow{\scriptstyle g'} & & \downarrow{\scriptstyle g} \\
\bullet & \xrightarrow{l} & \bullet & \xrightarrow{k} & \bullet
\end{array}$$

Notice that each square in this diagram is a pullback square. Proposition 2.1.4 bundled with property (2) of $\xi_P$ then guarantee that each of these squares has a preimage in $P$. Focusing on $\widetilde{\sigma}|\Delta^{\{0,1\}} \times \Delta^2$, we see that its preimage is a rectangle with both inner squares in $P$. Property (2) of $P$ implies the outer square of this rectangle is also in $P$. $\square$

**Lemma 2.1.6.** *Suppose that $\xi_P(\mathcal{C}; q)$ is stable under retracts in $\mathsf{Fun}(\Delta^1, q(\mathcal{D}))$ and that we have $h : \Delta^1 \times \Delta^1 \longrightarrow q(\mathcal{D})$, a homotopy between $\{f_i : \Delta^1 \longrightarrow q(\mathcal{D})\}_{i \in \{0,1\}}$. Then $f_0 \in \xi_P(\mathcal{C}; q)$ if and only if $f_1 \in \xi_P(\mathcal{C}; q)$.*

*Proof.* Homotopies between 1-simplices are identifiable with equivalences in $\mathsf{Fun}(\Delta^1, q(\mathcal{D}))$; hence, they are retracts. $\square$

**Proposition 2.1.7.** *Suppose that $\xi_P(\mathcal{C}; q)$ is stable under retracts in $\mathsf{Fun}(\Delta^1, q(\mathcal{D}))$ and that we are given the following 2-simplex in $q(\mathcal{D})$ where $g \in \xi_P(\mathcal{C}; q)$.*

$$\begin{array}{ccc}
 & B & \\
{\scriptstyle f} \nearrow & & \searrow {\scriptstyle g} \\
A & \xrightarrow{h} & C
\end{array}$$

*Then $h \in \xi_P(\mathcal{C}; q)$ if and only if $f \in \xi_P(\mathcal{C}; q)$.*

*Proof.* Let $\alpha : \Lambda^2_1 \longrightarrow q(\mathcal{D})$ be the inner horn with $\partial_0 \alpha = g$ and $\partial_2 \alpha = f$. By characterization of $\infty$-categories, all 2-simplices $\sigma : \Delta^2 \longrightarrow q(\mathcal{D})$ filling $\alpha$ are homotopic (cf. [i] 2.3.2.2). Thus, we can find an homotopy of 2-simplices that restricts to a homotopy between $\partial_1 \sigma$ and $g \circ f$. Combine Lemma 2.1.6 and Proposition 2.1.5. $\square$

Henceforth, we shall call a collection of morphisms in $\mathcal{C}$ satisfying the properties inhered in Proposition 2.1.3 and Proposition 2.1.5, *q-preadmissible*. We see instantly that if $E$ is such a collection of morphisms, then $(q(\mathcal{D}), E)$ is a geometric setup whenever $E$ is closed under retracts of morphisms. Furthermore, whenever the subcategory of $\mathcal{C}$ spanned by morphisms in $E$ has a Grothendieck topology, and $E$ is closed under retracts of morphisms, we obtain a collection of *admissible* structures on $q(\mathcal{D})$ (cf. [vi] 1.2.2). We will call a collection of morphisms that is *q-preadmissible* and closed under retracts *q-admissible*.

**Proposition 2.1.8.** *Let $E$ be the collection of morphisms in $q(\mathcal{D})$ that can be written, up to equivalence of arrows in $\mathsf{Fun}(\Delta^1, q(\mathcal{D}))$, as $f_1 f_0$ where for $i \in \{0,1\}$, $f_i \in E_i$ and $E_i$ is q-preadmissible. Then $E$ is a q-preadmissible whenever it is closed under precomposition with morphisms in $E_1$.*

*Proof.* It follows from the pasting law for pullbacks that $E$ is stable under base change in $q(\mathcal{D})$. It remains to show the analogue of Proposition 2.1.6 for $E$. In one direction, assume that $E$ is closed under composition and that we are given a 2-simplex $\sigma : \Delta^2 \longrightarrow q(\mathcal{D})$ with $\partial_0 \sigma = f$, $\partial_2 \sigma = h$ and $\partial_1 \sigma = f \circ h$ where $f \in E$ and $f \circ h \in E$. Let $g := f \circ h$ and consider the following limit diagram of a cartesian square

$$\begin{array}{c}\bullet \xrightarrow{h} \\ \downarrow^{h'} \searrow \\ \text{id} \quad \bullet \xrightarrow{g'} \bullet \\ \searrow \downarrow^{f'} \quad \downarrow^{f} \\ \bullet \xrightarrow{g} \bullet\end{array}$$

and note that $g' \in E$ and $f' \in E$. Observe also that it suffices to show that $h' \in E$ in order to show $h \in E$. Because $f'$ is in $E$, it can be written as $f' \simeq f'_1 f'_0$ where $f'_i \in E_i$. Thus, $\text{id} \simeq f'_1 f'_0 h'$. This in turn implies that $f'_0 h' \in E_1$. Let $p \simeq f'_0 h'$ and again, consider the following limit diagram of a cartesian square

$$\begin{array}{c}\bullet \xrightarrow{h'} \\ \downarrow^{j} \searrow \\ \text{id} \quad \bullet \xrightarrow{p'} \bullet \\ \searrow \downarrow^{f''_0} \quad \downarrow^{f'_0} \\ \bullet \xrightarrow{p} \bullet\end{array}$$

and notice that $f''_0 \in E_0$ and $p' \in E_1$. Moreover, since $\text{id} \simeq f''_0 j$, we deduce that $j \in E_0$. But $h' \simeq p'j$.

In the other direction, assume that whenever we have a 2-simplex $\sigma : \Delta^2 \longrightarrow q(\mathcal{D})$ with $\partial_1 \sigma = \partial_0 \sigma \circ \partial_2 \sigma$ where $\partial_0 \sigma \in E$ and $\partial_1 \sigma \in E$, then $\partial_2 \sigma \in E$. Consider a 2-simplex $\alpha : \Delta^2 \longrightarrow q(\mathcal{D})$ with $\partial_0 \alpha = f$, $\partial_2 \alpha = h$ and $\partial_1 \alpha = f \circ h =: g$ where both $f$ and $h$ are in $E$. From the latter, we may rewrite $h \simeq h_1 h_0$ so that $g \simeq f \circ h_1 h_0$. But the hypothesis on $E$ implies $f \circ h_1 \in E$. Ergo, we obtain a 2-simplex $\beta : \Delta^2 \longrightarrow q(\mathcal{D})$ with $\partial_1 \sigma = f \circ h_1$, $\partial_0 \beta = f_1$, and $\partial_2 \beta = f_0 \circ h_1$. By assumption, $f_0 \circ h_1$ is in $E \Rightarrow f_0 \circ h_1 \simeq k_1 k_0$ with $k_i \in E_i$. Hence, $g \simeq f_1 k_1 \circ k_0 h_0$. □

**Corollary 2.1.9.** *Suppose that $E$ is the collection of morphisms in $q(\mathcal{D})$ that can be written, up to equivalence of arrows in $\mathsf{Fun}(\Delta^1, q(\mathcal{D}))$, as $f_1 f_0$ where for $i \in \{0, 1\}$, $f_i \in E_i$ and $E_i$ is q-preadmissible. Let $E$ have the property that given any $f_0 f_1 \simeq f$, there exists $h_i \in E_i$ such that $f \simeq h_1 h_0$.*

*Proof.* The only obstruction to giving a proof that $E$ is q-preadmissible as is (described using $E_i$'s only and imposing no further condition(s)), is that it is not possible to do what the hypothesis allows. That is, composition generally looks like $k_1 k_0 r_1 r_0$ where $k_i$ and $r_i$ are in $E_i$. Being able to rewrite $k_0 r_1$ as $u_1 u_0$ gives $(k_1 u_1) \circ (u_0 r_0)$ as desired. □

**Definition 2.1.10.** Suppose that $\widehat{\mathcal{D}}$ is an essentially small presentable self-dual subcategory of $\mathsf{Pr}^L$ and whose inclusion is left exact and preserves small colimits. Moreover, suppose that it admits an $\infty$-symmetric monoidal structure $\widehat{\mathcal{D}}^{\boxtimes} \longrightarrow \mathsf{N}(\mathcal{F}\mathrm{in}_*)$ whereby the tensor product $\boxtimes : \widehat{\mathcal{D}} \times \widehat{\mathcal{D}} \longrightarrow \widehat{\mathcal{D}}$ preserves small colimits separately in each variab, andnd that its algebra objects are pointed and their tensor products are closed symmetric monoidal. Recall, as discussed earlier, $\mathsf{Pr}^L$ is the subcategory of $\widehat{\mathsf{Cat}_\infty}$ spanned by presentable $\infty$-categories and cocontinuous maps. We define $\Sigma_I(\widehat{\mathcal{D}})$ as a small collection of morphisms in $\widehat{\mathcal{D}}$ with the following characteristics:

1. $\Sigma_I(\widehat{\mathcal{D}})$ contains all identity maps and each morphism it contains is a categorical fibration. That is, a fibration in the Joyal model structure on $\mathsf{Set}_\Delta$. Moreover, if $q_0 \simeq q_2 q_1$ and $q_2 \in \Sigma_I(\widehat{\mathcal{D}})$, then $q_0 \in \Sigma_I(\widehat{\mathcal{D}}) \iff q_1 \in \Sigma_I(\widehat{\mathcal{D}})$.

2. Given any $\sigma : \Lambda_2^2 \longrightarrow \widehat{\mathcal{D}}$ with $\partial_0 \sigma \in \Sigma_I(\widehat{\mathcal{D}})$ and $\partial_1 \sigma \in \mathsf{CAlg}(\widehat{\mathcal{D}}^{\boxtimes})^{\mathrm{op}}$. There exists an horizontally right-adjointable square $\alpha : (\Lambda_2^2)^{\triangleleft} \longrightarrow \widehat{\mathcal{D}}$ with $\alpha|\Delta^{\{0\}} \times \Delta^1$ in $\Sigma_I(\widehat{\mathcal{D}})$, $\alpha|\Delta^1 \times \Delta^{\{0\}}$ in $\mathsf{CAlg}(\widehat{\mathcal{D}}^{\boxtimes})^{\mathrm{op}}$ and such that $\alpha|\Lambda_2^2 = \sigma$.

3. The subcategory of $\widehat{\mathcal{D}}$ spanned by the right adjoints of the morphisms in $\Sigma_I(\widehat{\mathcal{D}})$ is a subcategory of $\mathsf{CAlg}(\widehat{\mathcal{D}}^{\boxtimes})$.

4. Given $q : A \longrightarrow B$ in $\Sigma_I(\widehat{\mathcal{D}})$ and its right adjoint $q^*$, we have the following square in $\widehat{\mathcal{D}}$.

$$\begin{array}{ccc} A \times B & \xrightarrow{q \times 1} & B \times B \\ \downarrow{\scriptstyle 1 \otimes q^*} & & \downarrow{\scriptstyle \otimes} \\ A & \xrightarrow{q} & B \end{array}$$

5. Consider the following pasting diagram of pullback squares where $q \in \Sigma_I(\widehat{\mathcal{D}})$ and $f \in \mathsf{CAlg}(\widehat{\mathcal{D}}^{\boxtimes})$.

$$\begin{array}{ccc} \bullet & \longrightarrow & \bullet \\ \downarrow{\scriptstyle g'} & & \downarrow{\scriptstyle g} \\ \bullet & \longrightarrow & \bullet \\ \downarrow & & \downarrow{\scriptstyle q} \\ \bullet & \xrightarrow{f} & \bullet \end{array}$$

Then whenever $g$ belongs to $\Sigma_I(\widehat{\mathcal{D}})$, $g'$ also belongs to $\Sigma_I(\widehat{\mathcal{D}})$.

6. Let $\sigma : \Delta^1 \times \Delta^2 \longrightarrow \widehat{\mathcal{D}}$ be the pasting diagram of pullback squares visualized below

$$\begin{array}{ccc} \bullet & \longrightarrow & \bullet \\ \downarrow & & \downarrow g \\ \bullet & \longrightarrow & \bullet \\ \downarrow & & \downarrow q \\ \bullet & \xrightarrow{f} & \bullet \end{array}$$

where $q \in \Sigma_I(\widehat{\mathcal{D}})$ is fully faithful and $f \in \mathsf{CAlg}(\widehat{\mathcal{D}}^{\boxtimes})$. Then whenever $g$ belongs to $\Sigma_I(\widehat{\mathcal{D}})$, the square given as $\sigma|\Delta^1 \times \Delta^{\{0,1\}}$ is vertically right-adjointable.

7. It is closed under retracts.

We will write $\Sigma_I$ for $\Sigma_I(\widehat{\mathcal{D}})$ when there is no danger of ambiguity. When we strengthen condition (2) above such that all the squares obtained are pullbacks, we will write $\overline{\Sigma}_I$. We will write $\Sigma_I^\triangleleft$ (resp. $\overline{\Sigma}_I^\triangleleft$) for the collection of all squares in $\widehat{\mathcal{D}}$ spanned by vertical edges in $\Sigma_I$ (resp. $\overline{\Sigma}_I$) and horizontal edges in $\mathsf{CAlg}(\widehat{\mathcal{D}}^{\boxtimes})^{\mathrm{op}}$ and which are obtained from the process of condition(2) above.

**Definition 2.1.11.** We will define $\Sigma_I^{\mathfrak{p}}(\widehat{\mathcal{D}})$ as the subset of $\Sigma_I(\widehat{\mathcal{D}})$ that meets the criterion that follows. Let $q : A \longrightarrow B$ be in $\Sigma_I^{\mathfrak{p}}(\widehat{\mathcal{D}})$. Then the square of its pullback along $\otimes : B \times B \longrightarrow B$ is the following square where $q^*$ is right adjoint to $q$.

$$\begin{array}{ccc} A \times B & \xrightarrow{q \times 1} & B \times B \\ \downarrow 1 \otimes q^* & & \downarrow \otimes \\ A & \xrightarrow{q} & B \end{array}$$

We will write $\Sigma_I^{\mathfrak{p}}$ for $\Sigma_I^{\mathfrak{p}}(\widehat{\mathcal{D}})$ when there is no danger of ambiguity and $\overline{\Sigma}_I^{\mathfrak{p}}$ for $\Sigma_I^{\mathfrak{p}} \cap \overline{\Sigma}_I$.

**Remark 2.1.12.** At its heart, $\overline{\Sigma}_I$ is devised to behave simultaneously like the image of immersions of schemes under a reasonably well-behaved sheaf theory (such as Betti cohomology) and the image of an *admissibility structure* $\mathcal{T}^{\mathrm{ad}}$ on a *pregeometry* $\mathcal{T}$ under some $\mathcal{T}$-structure $\mathcal{O} : \mathcal{T} \longrightarrow \mathcal{Y}$. The latter is expressed most in criterion (2) and the former in criteria (3) and (4). In time, we should make these allusions abundantly clear. If, however, we loosen condition (2) so that we have $\Sigma_I$, then this is an imitation of immersions of schemes under nice enough sheaf theories. On the other hand, $\Sigma_I^{\mathfrak{p}}$ is meant to emulate the formulation of primes given in Definition 0.1.1, but with some variation. It is possible to use the exact definition here, and by demanding that Proposition 2.1.15 holds, we obtain the same kind of results as we obtain for the rest of the paper, with the only difference being the spaces reconstructed in Corollary 2.1.25.

**Lemma 2.1.13.** *Consider* $f : \Delta^1 \longrightarrow \mathsf{CAlg}(\widehat{\mathcal{D}}^{\boxtimes})$ *and* $q \in \Sigma_I$. *If* $q$ *is fully faithful and* $\partial_0 f = \partial_0 q$, *then* $f'$, *the pullback of* $f$ *along* $q$, *commutes with tensor products.*

*Proof.* From property (5), we know that $\Sigma_I$ is closed under pullbacks along maps in $\mathsf{CAlg}(\widehat{\mathcal{D}}^{\boxtimes})$. Additionally, $q'$ is fully faithful (cf. [xxvi] 4.6.2.7). Hence, we have homotopies $q'^* \circ q' \longrightarrow 1$ and $q' \circ (1 \otimes q'^*) \longrightarrow q' \otimes 1$. Combining these two, we obtain a homotopy $(\otimes \circ (q' \times q') \longrightarrow q' \circ (1 \otimes 1)$. The same is true for $q$. But $f' \circ q' \simeq q \circ f''$ and $q$ is conservative, so $f'$ commutes with tensor products. □

**Lemma 2.1.14.** *Suppose that we are given the following square in $\widehat{\mathcal{D}}$ where $S$ is a finite indexing set*

$$\begin{array}{ccc} \prod_{\alpha \in S} X'_\alpha & \xrightarrow{\prod_{\alpha \in S} f'_\alpha} & \prod_{\alpha \in S} X_\alpha \\ \downarrow{\prod_{\alpha \in S} q'_\alpha} & & \downarrow{\prod_{\alpha \in S} q_\alpha} \\ \prod_{\alpha \in S} Y'_\alpha & \xrightarrow{\prod_{\alpha \in S} f_\alpha} & \prod_{\alpha \in S} Y_\alpha \end{array} \quad (1)$$

*where for each $\alpha \in S$, $q_\alpha$ is a categorical fibration, and we have the following pullback square.*

$$\begin{array}{ccc} X'_\alpha & \xrightarrow{f'_\alpha} & X_\alpha \\ \downarrow{q'_\alpha} & & \downarrow{q_\alpha} \\ Y'_\alpha & \xrightarrow{f_\alpha} & Y_\alpha \end{array}$$

*Then (1) is a pullback square.*

*Proof.* We first observe that because $S$ is finite, $\prod_{\alpha \in S} q_\alpha$ is a categorical fibration; the objects involved are all fibrant in the Joyal model structure on simplicial sets. Moreover, the sequence of inclusions $\widehat{\mathcal{D}} \subseteq \mathsf{Pr}^L \subseteq \widehat{\mathsf{Cat}}_\infty$ is left exact. Altogether, this means, by Proposition A.2.4.4 of [i], we can compute the pullbacks in $\mathsf{Set}_\Delta$; precisely, in the category of small simplicial sets. In turn, we obtain the result by computing the relevant pullbacks pointwise. $\square$

**Proposition 2.1.15.** *Suppose we are given a pasting diagram of pullback squares below*

$$\begin{array}{ccc} A' & \xrightarrow{f''} & A \\ \downarrow{g'} & & \downarrow{g} \\ C & \xrightarrow{f'} & B \\ \downarrow{q'} & & \downarrow{q} \\ X & \xrightarrow{f} & Y \end{array}$$

*where $q \in \Sigma_I(\widehat{\mathcal{D}})$ is fully faithful and $f \in \mathsf{CAlg}(\widehat{\mathcal{D}}^{\boxtimes})$. Then $g \in \Sigma_I^{\mathfrak{p}}(\widehat{\mathcal{D}})$ implies $g' \in \Sigma_I^{\mathfrak{p}}(\widehat{\mathcal{D}})$.*

*Proof.* It is sufficient to show existence of a cube $\sigma : \Delta^1 \times \Delta^1 \times \Delta^1 \longrightarrow \widehat{\mathcal{D}}$ visualized below

$$\begin{array}{c} \begin{array}{ccc} A \times B & \longrightarrow & B \times B \\ \nwarrow_{f'' \times f'} & & \nwarrow_{f' \times f'} \\ A' \times C & \longrightarrow & C \times C \\ \downarrow & & \downarrow \\ A & \longrightarrow & B \\ \nwarrow_{f''} & & \nwarrow_{f'} \\ A' & \xrightarrow{g'} & C \end{array} \end{array} \quad (2)$$

where the front and back faces of the cube are due to property (4) of Definition 2.1.10; by property (5), $q'$ is in $\Sigma_I(\widehat{\mathcal{D}})$. If this were the case, then one can extract that the diagram $\lambda_1 : \Delta^2 \times \Delta^1 \longrightarrow \widehat{\mathcal{D}}$

$$\begin{array}{ccccc}
A' \times C & \xrightarrow{f'' \times f'} & A \times B & \xrightarrow{1 \otimes g^*} & A \\
\downarrow{\scriptstyle g' \times 1} & & \downarrow{\scriptstyle g \times 1} & & \downarrow{\scriptstyle g} \\
C \times C & \xrightarrow{f' \times f'} & B \times B & \xrightarrow{\otimes} & B
\end{array}$$

and the diagram $\lambda_2 : \Delta^2 \times \Delta^1 \longrightarrow \widehat{\mathcal{D}}$

$$\begin{array}{ccccc}
A' \times C & \xrightarrow{1 \otimes g'^*} & A' & \xrightarrow{f''} & A \\
\downarrow{\scriptstyle g' \times 1} & & \downarrow{\scriptstyle g'} & & \downarrow{\scriptstyle g} \\
C \times C & \xrightarrow{\otimes} & C & \xrightarrow{f'} & B
\end{array}$$

are homotopic; this homotopy, when one focuses only on the outer squares, is precisely $\sigma$. Now, observe that by the pasting law for pullbacks, it is enough to show that $\lambda_1$ is a pasting diagram of pullback squares to get the desired result. Notice that $g \times 1$ is a categorical fibration since it can be realized as a pullback of one. Subsequently, we may apply Lemma 2.1.14 to deduce that $\lambda_1|\Delta^{\{0,1\}}$ is a pullback square.

Contemplating the possibility of the existence of $\sigma$, we observe that it only remains to show that the left facing square of (2) exists. Combining Lemma 2.1.13 with property (6) of $\Sigma_I$, we obtain an equivalence between $f'' \circ (1 \otimes g'^*)$ and $(1 \otimes g^*) \circ (f'' \times f')$ as objects of $\mathsf{Fun}(\Delta^1, \widehat{\mathcal{D}})$. $\square$

**Proposition 2.1.16.** $\overline{\Sigma}_I^{\triangleleft}(\widehat{\mathcal{D}})$ *satisfies condition (1) and (2) of Definition 2.1.1 if all the vertical edges are (-1)-truncated.*

*Proof.* We imagine that the diagrams involved in Definition 2.1.1 are stood on their sides so that the horizontal arrows therein become vertical. Condition (2) is immediate given the pasting law for pullbacks. The first half of condition (1) is also a consequence of the pasting law for pullbacks and the 2-out-of-3 property spelled out in Definition 2.1.10 property (1).

It remains to show the second half of condition (1). But this reduces to the assertion that the pullback of the edge opposite to the degenerate edge along itself gives a (-2)-truncated map. By Lemma 5.5.6.15 of [i] the diagonal of the pullback of a map along itself is $(k-1)$-truncated if and only if the map is $k$-truncated for $k \geq -1$. $\square$

**Remark 2.1.17.** The proofs we have given above do not utilize the adjointness of the morphisms involved in the definition of $\Sigma_I$. Therefore, we can 'dualize' $\Sigma_I$ by switching every instance of "right adjoint" with "left adjoint" in the definition, and the preceding results will remain true. We shall write $\Sigma_P(\widehat{\mathcal{D}})$ for this dual ($\Sigma_P$ in short). In fact, as tools to build the framework that follows, the preference of $\Sigma_I$ over $\Sigma_P$ comes purely from considerations of convention. In a very real sense, it is the equivalent of choosing open subsets over closed subsets to describe topological spaces.

We will now use part of the data provided by $\Sigma_I$ and $\Sigma_I^{\mathfrak{p}}$ to introduce $\infty$-prosites and subsequently to construct a map $\mathsf{CAlg}(\widehat{\mathcal{D}}^{\boxtimes}) \longrightarrow \mathbb{J}_{\infty}^{\square}$.

**Definition 2.1.18.** Suppose that $\widehat{\mathcal{D}}$ is as described in 2.1.10 (henceforth, we shall take this to be the case). For $D \in \widehat{\mathcal{D}}$, isolate any two objects $d$ and $d'$. Say $d \leq d'$ if and only whenever $d$ lies in the essential image of a map in $\Sigma_I$, $d'$ lies in that essential image as well. This describes a proset $D_0^{\leq}$, and Proposition 1.2.5 constructs an $\infty$-proset $D^{\leq}$. Notice that this $\infty$-proset keeps track of homotopies; homotopic $n$-simplices are automatically equivalent as objects of $D_n^{\leq}$. This guarantees that $D$ is commensurate with the parallel proset structure $D_1 \rightrightarrows D_0$.

**Proposition 2.1.19.** *Let $f : C \longrightarrow D$ be a map in $\mathsf{CAlg}(\widehat{\mathcal{D}}^{\boxtimes})$. Then $f$ updates to a map $f : C^{\leq} \longrightarrow D^{\leq}$ of $\infty$-prosets.*

*Proof.* It suffices to show $f_0$ upgrades to a map $f_0 : C_0^{\leq} \longrightarrow D_0^{\leq}$ of prosites (cf. Proposition 1.2.9). Suppose $d \leq d'$, but $f_0(d) \not\leq f_0(d')$. This means that there exists a map $q \in \Sigma_1$ with target $D$ and whose essential image contains $f_0(d)$ but not $f_0(d')$. But the pullback of $q$ along $f$ is a map in $\Sigma_I$ whose essential image must contain both $d$ and $d'$. □

Let $q_{\alpha}(C_{\alpha})_0$ be all objects in the essential image of a morphism $q_{\alpha}$ in $\Sigma_I^{\mathfrak{p}}$ that lands in $D$. Given an object $d$ of $D \in \widehat{\mathcal{D}}$, let $[d] = \bigcap q_{\alpha}(C_{\alpha})_0$ where the intersection is taken over all essential images containing $d$. The association $d \longmapsto [d]$ describes a contravariant functor $\Gamma : D_0^{\leq} \longrightarrow \widehat{\mathsf{Cat}}$; here, $\widehat{\mathsf{Cat}}$ is the category of large categories.

**Lemma 2.1.20.** *$D_0^{\leq}$ is finitely complete and $\Gamma : (D_0^{\leq})^{op} \longrightarrow \widehat{\mathsf{Cat}}$ determines a packeted prosite $(D_0^{\leq}, J^{\Gamma})$.*

*Proof.* It is enough to show that $D_0^{\leq}$ contains pullbacks and an initial object. Because each map $q_{\alpha} : C_{\alpha} \longrightarrow D$ in $\Sigma_I$ is left adjoint, they preserve colimits. Therefore, the essential image of such a map contains the zero object of $D$ and admits coproducts. Unwinding what is required of pullbacks, we note that coproducts are satisfactory.

For $d \in D$, we notice that $\Gamma(d)$ contains the zero object and coproducts exist as seen in the previous paragraph; both of which coincide with those of $D_0^{\leq}$. If $c \in \Gamma(d)$, then $\Gamma(c) \subseteq \Gamma(d)$ by virtue of the definition involving taking intersections. □

**Lemma 2.1.21.** *Let $f : C \longrightarrow D$ be a map in $\mathsf{CAlg}(\widehat{\mathcal{D}}^{\boxtimes})$. Then $f$ updates to a map of $\infty$-prosites $f : (C^{\leq}, J) \longrightarrow (D^{\leq}, J')$.*

*Proof.* Proposition 2.1.15 informs us that $\Sigma_I^{\mathfrak{p}}$ is stable under base change along edges in $\mathsf{CAlg}(\widehat{\mathcal{D}}^{\boxtimes})$. Therefore, $f_0^{-1}(\Gamma'(f_0(d)) \supseteq \Gamma(d)$. This means $f_0$ is cover preserving. Recalling Proposition 1.3.3, it remains to show that $f_0$ is left exact. But this is guaranteed because $f$ is cocontinuous. □

**Remark 2.1.22.** The Grothendieck topology defined above is inspired by the reconstruction procedure in the appendix (cf. Theorem A.1.11). To recap, we obtain the space underlying a nice scheme $X$, for instance one where every quasicompact open immersion is quasi-perfect, by gluing the topological spaces obtained from coverings $\{K_{\alpha} \hookrightarrow D_{\text{perf}}(X)\}_{\alpha \in \Omega}$ where each $K_{\alpha}$ looks like the image of $j_! : D_{\text{perf}}(U) \longrightarrow D_{\text{perf}}(X)$ for some quasicompact open immersion $j : U \hookrightarrow X$ (of said nice schemes). Specifically, we want each $K_{\alpha}$ to behave like a prime, as in Definition 0.1.1. In the new context above, the topology on $D_0^{\leq}$ is set up to be the coarsest topology such that for all with $q \in \Sigma_I^{\mathfrak{p}}$ $codom(q) = D$, $q^* : D_0^{\leq} \longrightarrow dom(q)_0^{\leq}$ is a morphism of prosites when $dom(q)_0^{\leq}$ is given the trivial Grothendieck topology; that is, its packeting is the functor associating each object with $dom(q)_0^{\leq}$.

**Remark 2.1.23.** Suppose that each $q : C \longrightarrow D$ in $\Sigma_I$ is fully faithful. It is immediate that $Q = qq^*$ is a colocalization. Moreover, from the equivalence $q(1 \otimes q^*) \simeq q \otimes 1$, we conclude that $Q \simeq Q(1 \otimes \mathbf{1}_\otimes) \simeq 1 \otimes Q(\mathbf{1}_\otimes)$. Therefore, $Q$ is a *smashing colocalization*. A similar argument used in the context of $\Sigma_P$ produces *smashing localizations*. On the other hand, $\Sigma_I^{\mathfrak{p}}$ under the fully faithful assumption, gives the association $d \longmapsto [d]$ a characterization akin to the Zariski spectrum of Remark 3.10 in [xxiii]. In fact, when dealing with higher enhancements of the usual derived categories of complexes of $O_X$-modules with quasicoherent cohomology, these constructions "almost" collide when looking at a nice enough scheme $X$ and make slight variations on the definition of primeness. That is, we obtain an embedding of the Zariski spectrum into the space obtained by the method we use. Furthermore, although we obtain a larger space, we are still able to fully access the initial scheme since it turns out the naturally occuring sheaves of rings, constructed after Section 7 of [vii], promote the topological embedding to a ringed space embedding (cf. Theorem A.2.8 and Remark A.2.9).

**Proposition 2.1.24.** *There exists a map* $\mathsf{v}_I : \mathsf{CAlg}(\widehat{\mathcal{D}}^{\boxtimes}) \longrightarrow \mathsf{J}_\infty^{\square}$ *of $\infty$-categories.*

*Proof.* Begin with the map $\mathsf{CAlg}(\widehat{\mathcal{D}}^{\boxtimes}) \longrightarrow \mathsf{N}(\mathsf{hCAlg}(\widehat{\mathcal{D}}^{\boxtimes}))$ induced by the adjunction $\mathsf{h} \dashv \mathsf{N}$. Taking the homotopy categories of each of the objects of $\mathsf{N}(\mathsf{hCAlg}(\widehat{\mathcal{D}}^{\boxtimes}))$, and after unwinding Definition 1.3.5 and recalling Theorem 1.3.17, we realize a map $\mathsf{N}(\mathsf{hCAlg}(\widehat{\mathcal{D}}^{\boxtimes})) \longrightarrow \mathsf{J}_\infty^{\square}$. The composition of all of the above maps is the desired map. $\square$

**Corollary 2.1.25.** *There exists a map* $\mu_I(S) : \mathsf{CAlg}(\widehat{\mathcal{D}}^{\boxtimes}) \longrightarrow \mathsf{Stk}_S^{op}$ *of $\infty$-categories.*

*Proof.* This follows from Definition 1.3.23 and Definition 1.3.24. $\square$

We will omit $S$ and write $\mu_I$ to simplify our notation. In fact, as we proceed, we will still simply write $\mu_I$ even though this map will become dependent on more than one factor; doing otherwise would overburden the notation. Lucky for us, the results we hope to obtain will not depend on these factors.

We say that a collection of edges in $\widehat{\mathcal{D}}$ admits a $\Sigma_I$-structure (resp. $\overline{\Sigma}_I$ structure) if that collection of morphisms meets all the criteria laid out in Definition 2.1.10 (resp. Definition 2.1.10 with (2) strengthened accordingly). Crucially, we see that such a collection of morphisms also induces its own map into animated $S$-stacks; the $I$ in the subscript of $\mu_I$ is meant to keep track of this observation.

**Theorem 2.1.26.** *Suppose $\mathcal{C}$ is an essentially small $\infty$-category and that $\pi^{op} : \mathcal{C} \longrightarrow \widehat{\mathcal{D}}$ is a map with an essential image whose edges admit a $\Sigma_I$-structure. Furthermore, assume that the restriction of $\mu_I$ to the opposite subcategory of the subcategory spanned by the edges of the squares in*

$\Sigma_I^\triangleleft$ is (-1)-truncated. Then we have the following diagram of $\infty$-categories

$$\begin{array}{ccc} \mathcal{C}^{op} & \xrightarrow{\pi} & \mathsf{CAlg}(\widehat{\mathcal{D}}^\boxtimes) \\ {\scriptstyle (\_)^\pi} \downarrow & \nearrow {\scriptstyle \pi_0} & \\ \mathsf{Stk}^{op}_{S\,|_\pi} & & \end{array}$$

where $\mathsf{Stk}^{op}_{S\,|_\pi}$ is the largest subcategory of the essential image of $\mu_I$ with pushouts and such that the inclusion $\pi(\mathcal{C}^{op}) \subseteq \mathsf{CAlg}(\widehat{\mathcal{D}})$ admits a left Kan extension along $\mu_I$ and the inclusion of the essential image of $\mu_I \circ \pi$ is right exact. Furthermore, $\pi_0$ admits a section $\mu_\pi : \pi(\mathcal{C}^{op}) \longrightarrow \mathsf{Stk}^{op}_{S\,|_\pi}$.

*Proof.* Beginning with Corollary 2.1.25 and taking into consideration the fact that $\pi(\mathcal{C}^{\mathrm{op}})$ is a subcategory of $\mathsf{CAlg}(\widehat{\mathcal{D}}^\boxtimes)$ by virtue of (3) of Definition 2.1.10, one obtains $(\_)^\pi$ as $\mu_I \circ \pi$. In fact, since $\pi(\mathcal{C}^{\mathrm{op}})$ is a subcategory of the subcategory of $\widehat{\mathcal{D}}$ spanned by the edges of the squares in $(\Sigma_I^\triangleleft)^{\mathrm{op}}$, the restriction $\mu_I | \pi(\mathcal{C}^{\mathrm{op}})$ is (-1)-truncated. Thus, it is fully faithful as a map into its essential image. We will label it $\mu_\pi$. By Corollary 3.2.3.5 of [xvi], it follows that since $\widehat{\mathcal{D}}$ is a presentable symmetric monoidal $\infty$-category and the tensor product commutes with small colimits in each variable, small colimits in $\mathsf{CAlg}(\widehat{\mathcal{D}}^\boxtimes)$ exist. Consequently, the left Kan extension of the inclusion functor $\pi(\mathcal{C}^{\mathrm{op}}) \subseteq \mathsf{CAlg}(\widehat{\mathcal{D}}^\boxtimes)$ along $\mu_\pi$ exists (cf. [i] 4.3.2.2 and 4.3.2.6). Hence, we are guaranteed both the existence and non-triviality of $\mathsf{Stk}^{\mathrm{op}}_{S\,|_\pi}$. Label $\pi_0$ the left Kan extension that arises from its description. $\square$

**Definition 2.1.27.** Suppose $\pi^{\mathrm{op}} : \mathcal{C} \longrightarrow \widehat{\mathcal{D}}$ satisfies the same conditions as those of Theorem 2.1.26 initially, but with a slight variation. We have a subcategory $\mathcal{C}' \subseteq \mathcal{C}$ spanned by edges of covering sieves generating a topology on $\mathcal{C}$, and such that the essential image of $\pi^{\mathrm{op}} | \mathcal{C}'$ admits a $\overline{\Sigma_I}$-structure. Then we will say that the map $\pi^{\mathrm{op}}$ is a *quasi-suprematic space*.

**Definition 2.1.28.** We lift the definitions of structured spaces from [iv]. Let $\mathcal{G}$ be an essentially small $\infty$-category; that is, its minimal model is a small $\infty$-category. An *admissibility structure* on $\mathcal{G}$ is a subcategory $\mathcal{G}^{\mathrm{ad}} \subseteq \mathcal{G}$ such that

1. It contains every object of $\mathcal{G}$. The morphisms that belong to it will be called *admissible* morphisms.

2. $\mathcal{G}^{\mathrm{ad}}$ has a Grothendieck topology.

Additionally, admissible morphisms satisfy the following conditions.

1. Stability under base change.

2. Given any 2-simplex $\sigma : \Delta^2 \longrightarrow \mathcal{G}$ such that $\partial_0 \sigma$ is admissible, then if $\partial_1 \sigma$ is admissible, $\partial_2 \sigma$ is admissible.

3. Closure under retracts of morphisms.

We say $\mathcal{G}$ is a *geometry* if it admits finite limits and idempotent complete. An essentially small $\infty$-category $\mathcal{T}$ with an admissibility structure $\mathcal{T}^{\mathrm{ad}}$ is called a *pregeometry* if it admits finite products.

Given an ∞-topos $\mathcal{X}$ and a geometry $\mathcal{G}$, a map $\mathcal{O} : \mathcal{G} \longrightarrow \mathcal{X}$ is called a $\mathcal{G}$-structure if it is left exact and given any admissible cover $\{U_i \longrightarrow X\}_{i \in I}$, the induced map $\coprod_{i \in I} \mathcal{O}(U_i) \longrightarrow \mathcal{O}(X)$ is an effective epimorphism. On the other hand, given an ∞-topos $\mathcal{X}$ and a pregeometry $\mathcal{T}$, a map $\mathcal{O} : \mathcal{T} \longrightarrow \mathcal{X}$ is called a $\mathcal{T}$-structure if it preserves finite products, preserves pullbacks when restricted to $\mathcal{T}^{\text{ad}}$, and given any admissible cover $\{U_i \longrightarrow X\}_{i \in I}$, the induced map $\coprod_{i \in I} \mathcal{O}(U_i) \longrightarrow \mathcal{O}(X)$ is an effective epimorphism.

$\mathsf{Str}_\mathcal{G}(\mathcal{X})$ is the full subcategory of $\mathsf{Fun}(\mathcal{G}, \mathcal{X})$ spanned by $\mathcal{G}$-structures. We define $\mathsf{Str}_\mathcal{T}(\mathcal{X})$ similarly using $\mathcal{T}$-structures. $\mathsf{Str}_\mathcal{T}^{\text{loc}}(\mathcal{X}) \subseteq \mathsf{Str}_\mathcal{T}(\mathcal{X})$ is the subcategory spanned by edges that result in pullback squares when evaluated at $\mathcal{T}^{\text{ad}}$. That is, given such an edge $\mathcal{T}^{\text{ad}} \times \Delta^1 \longrightarrow \mathcal{X}$, when one restricts to a particular 1-simplex of $\mathcal{T}^{\text{ad}}$, one obtains a pullback square $\Delta^1 \times \Delta^1 \longrightarrow \mathcal{X}$.

**Example 2.1.29.** We lift from [xxiv] the following instance of a pregeometry. We have a pregeometry $\mathcal{T}_{\text{an}}$ as follows:

1. The underlying category of $\mathcal{T}_{\text{an}}$ is the category of smooth $k$-analytic spaces;
2. A morphism in $\mathcal{T}_{\text{an}}$ is *admissible* if and only if it is étale.
3. The topology on $\mathcal{T}_{\text{an}}$ is the étale topology.

**Definition 2.1.30.** Let $\mathcal{X}^L$ be an ∞-topos meeting the same conditions as $\widehat{\mathcal{D}}$ and whose symmetric monoidal structure is the one carried by the finite product. Henceforth, we shall use this designation for only such an ∞-topos. Given a pregeometry $\mathcal{T}$ where all admissible covering maps are (-1)-truncated and a $\mathcal{T}$-structure $\mathcal{O}^{\text{op}} : \mathcal{T} \longrightarrow \mathcal{X}^L$, $\mathcal{O}^{\text{op}}$ is a *suprematic* space over $\mathcal{X}^L$ if and only if $\mathcal{O}^{\text{op}}|\mathcal{T}^{\text{ad}}$ is a quasi-suprematic space. We will often infer to a suprematic space with the assumption that it is over an ∞-topos already made. We can now state Theorem A in its fullness as follows.

**Theorem 2.1.31.** *Let $\mathcal{O}^{op} : \mathcal{T} \longrightarrow \mathcal{X}^L$ be a suprematic space. Then there exist maps of ∞-categories $(\_)^\mathcal{O} : (\mathcal{T}^{ad})^{op} \longrightarrow \mathsf{Stk}^{op}_{S\,|_\mathcal{O}}$ and $\mathcal{O}_0 : \mathsf{Stk}^{op}_{S\,|_\mathcal{O}} \longrightarrow \mathcal{X}^L$, extending $\mathcal{O}|(\mathcal{T}^{ad})^{op} : (\mathcal{T}^{ad})^{op} \longrightarrow \mathcal{X}^L$ as $\mathcal{O}_0 \circ (\_)^\mathcal{O}$. Furthermore, $\mathcal{O}_0$ admits a section $\mu_\mathcal{O} : \mathcal{O}((\mathcal{T}^{ad})^{op}) \longrightarrow \mathsf{Stk}^{op}_{S\,|_\mathcal{O}}$.*

*Proof.* Theorem 2.1.26 and the forgetful functor $\mathsf{CAlg}(\mathcal{X}^L) \longrightarrow \mathcal{X}^L$. □

Theorem A is a trivial consequence of Theorem 2.1.26 and does not, of itself, command the attention we have paid it thus far. However, speaking for the ambition of this paper, concerns of thematic unity come together to grant it crucial importance; indeed, suprematic spaces are the bonds that cleave together the seemingly far-flung enterprises of this project. Its statement here is simply an exercise in documenting a cohesive whole.

### A Parametrization of Six Functor Formalisms

In this subsection, our main aim is to prove Theorem B. The outline of the proof is as follows. First, we will show that for any quasi-suprematic space $\pi^{\text{op}} : \mathcal{C} \longrightarrow \widehat{\mathcal{D}}$, $\mathsf{Stk}_{S|_\pi}$ is imbued with a $\mu_\pi$-admissibility structure that is at the same time a geometric setup. Furthermore, this geometric setup conforms to the types desired in [ii] Proposition A.5.10. Therefore, the map $\pi_0 : \mathsf{Stk}^{\text{op}}_{S\,|_\pi} \longrightarrow \widehat{\mathcal{D}} \subseteq \mathsf{Pr}^L$ canonically updates to a 6-functor formalism as soon as $\mathsf{Stk}^{\text{op}}_{S\,|_\pi}$ admits finite coproducts. In effect, this suggests that collections of 6-functor formalisms on appropriately chosen subcategories of $\mathsf{Stk}_S$ are tracked by quasi-suprematic spaces. Taking the effort to sharpen and concretize this intuition is how we arrive at Theorem B. In giving its proof, the result of Lemma 2.4.2 in [xxv] is pivotal.

**Proposition 2.2.1.** *Suppose $\mathcal{C}$ is an essentially small $\infty$-category and that $\pi^{op} : \mathcal{C} \longrightarrow \widehat{\mathcal{D}}$ is a map with an essential image whose edges admit a $\Sigma_I$-structure. Furthermore, assume that the restriction of $\mu_I$ to the opposite subcategory of the subcategory spanned by the edges of the squares in $\Sigma_I^\triangleleft$ is (-1)-truncated. Then there exists a geometric setup $(\mathsf{Stk}_{S|_\pi}, E_\pi^I)$ with the property that given any 2-simplex $\sigma : \Delta^2 \longrightarrow \mathsf{Stk}_{S|_\pi}$ with $\partial_0 \sigma \in E_\pi^I$, then $\partial_2 \sigma \in E_\pi^I$ if and only if $\partial_1 \sigma \in E_\pi^I$.*

*Proof.* Observe that $\mathsf{Stk}_{S|_\pi}$ has pullbacks by definition and products since inclusion of $\mu_\pi$ into $\mathsf{Stk}_{S|_\pi}^{op}$ is right exact; the essential image of $\pi$ has coproducts since $\pi^{op}$ is a $\mathcal{T}$-structure. Let $P$ be the largest subset of $\Sigma_I^\triangleleft$ that meets all the requirements of Definition 2.1.1. Looking through the conditions set out by this definition, one notices that the subset of squares spanned by degenerate edges is a subset of $P$. This tells us that $\xi_P(\widehat{\mathcal{D}}; \mu_I^{op})$ is non-empty since degenerate edges once in the image of the subcategory spanned by edges in $P$, automatically populate $\xi_P(\widehat{\mathcal{D}}; \mu_I^{op})$ since they meet conditions (2) and (3) of Definition 2.1.2. Condition (3) is met due to the 2-out-of-3 property of $\Sigma_I$ and condition (2) due to property (2) of $\Sigma_I$.

Set $E_\pi^I$ as the collection of all edges in $\xi_P$ up to equivalence in $\mathsf{Fun}(\Delta^1, \mathsf{Stk}_{S|_\pi})$ and observe that by virtue of Proposition 2.1.3 and Proposition 2.1.5, we obtain a geometric setup that satisfies the 2-out-of-3 property spelled out. $\square$

**Definition 2.2.2.** In Remark 2.1.17, we mentioned $\Sigma_P$ as the "dual" of $\Sigma_I$. It seems yet again that we can also arrive at a geometric setup $(\mathsf{Stk}_{S|_\pi}, E_\pi^P)$ from $\Sigma_P$. However, we need to be careful since suprematic spaces are defined with a particular $\Sigma_I$ structure in mind; the same applies to how $\mathcal{C}$ above is set up. We still want to have some $\Sigma_P$-structure to play the role of "proper maps" and therefore we proceed as follows to formulate it. We begin by defining $\Sigma_P$ as if "dualizing" $\Sigma_I$; that is, we switch the "right-adjointness" in Definition 2.1.10 with "left-adjointness". Then we demand that $\Sigma_P^\triangleleft$ is constituted by squares each with one pair of opposite edges in $\Sigma_I$ and the other pair in $\Sigma_P$. Furthermore, we add the extra condition that whenever $f \in \Sigma_I \cap \Sigma_P$ then $f$ is $n$-truncated for some $n \geq 2$. If we proceed as we did to arrive at $E_\pi^I$, we arrive at $E_\pi^P$. At the very least, $E_\pi^P$ is populated by isomorphisms.

**Proposition 2.2.3.** *Suppose $\mathcal{C}$ is an essentially small $\infty$-category and that $\pi^{op} : \mathcal{C} \longrightarrow \mathcal{X}^L$ is a map with an essential image whose edges admit a $\Sigma_I$-structure. Furthermore, assume that the restriction of $\mu_I$ to the opposite subcategory of the subcategory spanned by the edges of the squares in $\Sigma_I^\triangleleft$ is (-1)-truncated. Then there exists a 6-functor formalism $\mathcal{D}^\pi : \mathrm{Corr}(\mathsf{Stk}_{S|_\pi}, E_\pi^I) \longrightarrow \mathcal{X}^L$ such that $\mathcal{D}^\pi | \mathrm{Corr}(\mathsf{Stk}_{S|_\pi}, \mathrm{isom}) = \pi_0$.*

*Proof.* We observe that since we begin with a map $\pi_0 : \mathsf{Stk}_{S|_\pi}^{op} \longrightarrow \mathsf{CAlg}(\mathcal{X}^L)$, we need only to see that $E_\pi^I$ meets certain conditions. Specifically, those laid out in Proposition A.5.10 of [ii]. By formulation, the edges in $E_\pi^I$ correspond to those in $\Sigma_I$ under the action of $\pi_0^{op}$. Hence, under $\pi_0^{op}$, they admit right adjoints, the Beck-Chevalley transformation is an equivalence, and they obey the projection formula with respect to their right adjoints. To reflect the notation used therein, we set $I = E_\pi^I$ and for $P$, we collect the set of all isomorphisms. $\square$

Let $E_1$ be the smallest $\mu_I$-preadmissible subset of $E_\pi^P$ with the property that if $\overline{E}$ is a collection of edges determined as in Proposition 2.1.8 with $E_0 = E_\pi^I$ and $E_1$, precomposition of $f \in \overline{E}$ with any $g \in E_1$ is again in $\overline{E}$. Because the collection of isomorphisms is a subset of $E_\pi^P$, $E_1$ is non-empty.

**Theorem 2.2.4.** *Suppose $\mathcal{C}$ is an essentially small $\infty$-category and that $\pi^{op} : \mathcal{C} \longrightarrow \mathcal{X}^L$ is a map with an essential image whose edges admit a $\Sigma_I$-structure. Furthermore, assume that the restriction of $\mu_I$ to the opposite subcategory of the subcategory spanned by the edges of the squares in $\Sigma_I^\lhd$ is (-1)-truncated. Then there exists a 6-functor formalism $\mathcal{D}^\pi : \mathrm{Corr}(\mathsf{Stk}_{S|_\pi}, \overline{E}) \longrightarrow \mathcal{X}^L$ such that $\mathcal{D}^\pi | \mathrm{Corr}(\mathsf{Stk}_{S|_\pi}, \mathrm{isom}) = \pi_0$. Furthermore, the following are true.*

1. *Each $f \in \overline{E}$ can be decomposed as $f \simeq f_1 f_0$ where $f_i \in E_i$.*

2. *For all $f \in E_1$, $f_!$ is left adjoint to $f^*$, obeys the projection formula, and the Beck-Chevalley transformation is an equivalence.*

3. *For all $f \in E_1$, $f_*$ is right adjoint to $f^*$, obeys the projection formula, and the Beck-Chevalley transformation is an equivalence.*

4. *If $f \in E_0 \cap E_1$, then $f_* \simeq f_!$.*

*Proof.* We proceed as in the previous case. However, we need to show that $\overline{E}$ is sound relative to the requirements of Proposition A.5.10 of [ii]. That it is a geometric setup follows from Proposition 2.1.8. Now observe that by description, if $f \in E_1 \cap E_0$, it is $n$-truncated for some $n \geq -2$. Finally, given any pullback of a map in $f \in E_0$ along a map in $g \in E_1$, we know there exists a square in $E_P^\lhd$ corresponding to this pullback square. $\square$

In what is to follow, we will modify $\overline{E}$ above demanding that $E_1 \subseteq E_0$ and when necessary, we will indicate the dependence of $E$ on a suprematic space $\pi$ denoting it as $E_\pi$.

**Definition 2.2.5.** We will say that two suprematic spaces $\pi_i^{op} \in \mathsf{Str}_\mathcal{T}(\mathcal{X}^L)$ have the same *geometric content* if they induce the same geometric setup up to categorical equivalence. That is, we say $(\mathsf{Stk}_{S|_\pi}, E_\pi) \simeq (\mathsf{Stk}_{S|_{\pi'}}, E_{\pi'})$ if and only if $\mathsf{Stk}_{S|_\pi} \simeq \mathsf{Stk}_{S|_{\pi'}}$, $E_\pi = E_{\pi'}$, and $(\_)^\pi \simeq (\_)^{\pi'}$. In general, for a simplicial set $K$, we will say that a map $\Delta^1 \longrightarrow \mathsf{Fun}(K, \mathcal{X}^L)$ is *locally vertically right/left adjointable* if upon restriction to a specific 1-simplex of $K$, the square obtained $\Delta^1 \times \Delta^1 \longrightarrow \mathcal{X}^L$ has a dual square in $\mathsf{CAlg}(\mathcal{X}^L)$ that is vertically right/left adjointable. we will denote the subcategory of $\mathsf{Fun}(K, \mathcal{X}^L)$ spanned by edges that are vertically right/left adjointable when restricted to edges $E$ of $\mathcal{T}^{ad}$ by $\mathsf{Fun}_{\mathrm{Adj}}((K, E), \mathcal{X}^L)$. We define $\mathsf{Sup}_\mathcal{T}^\otimes(E, \mathcal{X}^L)$ as follows.

1. It is a full subcategory of $\mathsf{Fun}_{\mathrm{Adj}}((\mathcal{T}^{ad}, E), \mathcal{X}^L)$ spanned by suprematic spaces with the same geometric content.

2. Given any two $n$-simplices of $\mathsf{Fun}(\mathsf{Stk}_{S|}^{op}, \mathsf{CAlg}(\mathcal{X}^L))$ spanned by vertices that are left Kan extensions of objects $(\mathcal{T}^{ad})^{op} \longrightarrow \mathsf{CAlg}(\mathcal{X}^L)$ corresponding to objects of $\mathsf{Sup}_\mathcal{T}^\otimes(E, \mathcal{X}^L)$, the two $n$-simplices agree when restricted to $(\underline{(\mathcal{T}^{ad})^{op}})$ if and only if they are homotopic.

Where there is no need to specify $E$, we will simply write $\mathsf{Sup}_\mathcal{T}^\otimes(\mathcal{X}^L)$. And to indicate the subcategory of animated stacks associated with $\mathsf{Sup}_\mathcal{T}^\otimes(E, \mathcal{X}^L)$, we will simply write $\mathsf{Stk}_{S|}$, $E$ to stand in for $E_\pi$, and () to stand in for $()^\pi$. There is the possibility of numerous subcategories of $\mathsf{Fun}_{\mathrm{Adj}}((\mathcal{T}^{ad}, E), \mathcal{X}^L)$ that fit our description. The results obtained here do not depend on a choice of them.

**Proposition 2.2.6.** *There exists a map $f : \mathsf{Sup}^{\otimes}_{\mathcal{T}}(E, \mathcal{X}^L)^{\mathrm{op}} \longrightarrow \mathsf{Fun}^{\otimes,\mathrm{lax}}(\mathsf{Stk}_{S_|}^{\mathrm{op}}, \mathcal{X}^L)$ which is a monomorphism of simplicial sets.*

*Proof.* We understand from Theorem 2.1.31 that every suprematic space $\pi^{op} : \mathcal{T}^{\mathrm{ad}} \longrightarrow \mathcal{X}^L$ determines a map $\pi_0 : \mathsf{Stk}_{S_{|\pi}}^{\mathrm{op}} \longrightarrow \mathsf{CAlg}(\mathcal{X}^L)$. Furthermore, this determination occurs as the left Kan extension of the inclusion $\pi((\mathcal{T}^{\mathrm{ad}})^{\mathrm{op}}) \subseteq \mathsf{CAlg}(\mathcal{X}^L)$ along $\mu_\pi : \pi((\mathcal{T}^{\mathrm{ad}})^{\mathrm{op}}) \longrightarrow \mathsf{Stk}_{S_{|\pi}}^{\mathrm{op}}$ which is another faithful map. This is to say, $\pi : (\mathcal{T}^{\mathrm{ad}})^{\mathrm{op}} \longrightarrow \mathsf{CAlg}(\mathcal{X}^L)$ is extended along $(\_)^\pi : (\mathcal{T}^{\mathrm{ad}})^{\mathrm{op}} \longrightarrow \mathsf{Stk}_{S_{|\pi}}^{\mathrm{op}}$ via $\pi_0$. Therefore, given a subcategory $\mathcal{B} \subseteq \mathsf{Fun}((\mathcal{T}^{\mathrm{ad}})^{\mathrm{op}}, \mathsf{CAlg}(\mathcal{X}^L))$ spanned by suprematic spaces with the same geometric content, we obtain a map $\mathcal{B} \longrightarrow \mathsf{Fun}(\mathsf{Stk}_{S_|}^{\mathrm{op}}, \mathsf{CAlg}(\mathcal{X}^L))$. Ergo, we obtain a map $f : \mathsf{Sup}^{\otimes}_{\mathcal{T}}(E, \mathcal{X}^L)^{\mathrm{op}} \longrightarrow \mathsf{Fun}^{\otimes,\mathrm{lax}}(\mathsf{Stk}_{S_|}^{\mathrm{op}}, \mathcal{X}^L)$ taking into account Theorem 2.4.3.18 of [xvi].

Now, every two $n$-simplices $\sigma$ and $\sigma'$ of $\mathsf{Sup}^{\otimes}_{\mathcal{T}}(E, \mathcal{X}^L)$ that agree (up to homotopy) under $f$ must agree on their vertices. But by characterization of Kan extensions, these two vertices must agree pointwise when restricted to $(\mathcal{T}^{\mathrm{ad}})^{\mathrm{op}}$. Subsequently, the two $n$-simplices must agree, up to homotopy, on $(\mathcal{T}^{\mathrm{ad}})^{\mathrm{op}}$. $\square$

We will write, by abuse of notation, $\mathsf{Fun}^{\otimes,\mathrm{lax}}_{\mathrm{Adj}}((\mathsf{Stk}_{S_|}^{\mathrm{op}}, E_0 \circ E_1), \mathcal{X}^L)$ to indicate the subcategory of $\mathsf{Fun}(\mathsf{Stk}_{S_|}^{\mathrm{op}}, \mathsf{CAlg}(\mathcal{X}^L))$ with edges that become locally vertically right adjointable when restricted to morphisms in $E_1$ and locally vertically left adjointable when restricted to morphisms in $E_0$. Set $I := E_0$ and $P := E_1$. We write $\mathsf{Sup}^{\otimes}_{\mathcal{T}}(I \circ P, \mathcal{X}^L)$ for the subcategory of $\mathsf{Sup}^{\otimes}_{\mathcal{T}}(E, \mathcal{X}^L)$ that has the same characteristics with respect to the preimages of $I$ and $P$ under all $(\_)^{\mathrm{op}} : \mathcal{T}^{\mathrm{ad}} \longrightarrow \mathsf{Stk}_{S_|}$ induced by the suprematic spaces in question. And when the context is sufficiently clear, to simplify notation, we will write $\mathsf{Sup}^{\otimes}_{\mathcal{T}}(E, \mathcal{X}^L)$ for $\mathsf{Sup}^{\otimes}_{\mathcal{T}}(I \circ P, \mathcal{X}^L)$.

**Corollary 2.2.7.** *There exists a map $f : \mathsf{Sup}^{\otimes}_{\mathcal{T}}(I \circ P, \mathcal{X}^L)^{\mathrm{op}} \longrightarrow \mathsf{Fun}^{\otimes,\mathrm{lax}}_{\mathrm{Adj}}((\mathsf{Stk}_{S_|}^{\mathrm{op}}, I \circ P), \mathcal{X}^L)$ which is a fully faithful map of $\infty$-categories.*

*Proof.* Keeping in mind the map $f$ obtained in Proposition 2.2.6 restricted to $\mathsf{Sup}^{\otimes}_{\mathcal{T}}(I \circ P, \mathcal{X}^L)^{\mathrm{op}}$, we consider a homotopy $h : \partial \Delta^n \times \Delta^1 \longrightarrow \mathsf{Fun}^{\otimes,\mathrm{lax}}_{\mathrm{Adj}}((\mathsf{Stk}_{S_|}^{\mathrm{op}}, I \circ P), \mathcal{X}^L)$ such that $h|\partial \Delta^n \times \{0\} = \sigma$ and $h|\partial \Delta^n \times \{1\} = f(\sigma')$ where $\sigma'$ is an $n$-simplex of the domain of $f$. Therefore, this homotopy restricts to a homotopy $h : \partial \Delta^n \times \Delta^1 \longrightarrow \mathsf{Fun}^{\otimes,\mathrm{lax}}_{\mathrm{Adj}}(((\mathcal{T}^{\mathrm{ad}})^{\mathrm{op}}, I \circ P), \mathcal{X}^L)$. But this means that $\sigma$ restricts to an $n$-simplex in the image of $f$. It is left to show that if two $n$-simplices spanned by suprematic spaces agree once restricted to $(\mathcal{T}^{\mathrm{ad}})^{\mathrm{op}}$ they agree everywhere. But this is guaranteed by property (2) of Definition 2.2.5. $\square$

Given a multi-simplicial set $X$ of order $m$, let $\delta^* X$ be the simplicial set whose $n$-simplices are the maps $\Delta^n \times \cdots \times \Delta^n \longrightarrow X$ where the product on the right is taken $m$-times. In our situation, $m$ will at most be 3. In the case where we are given a simple set $S$, we can define an $m$-simplicial set whose $(n_1, \ldots, n_m)$-simplices are maps $\Delta^{n_1} \times \cdots \times \Delta^{n_m} \longrightarrow S$ that have all $k$-dimensional cubes being Cartesian for $2 \leq k \leq m$. If we wish some of the maps to be contravariant and others covariant, we will denote the resulting $m$-simplicial set as $\overline{\mathrm{Corr}}(S, [k_1, k_2])$ where $0 \leq k_1 \leq k_2 \leq m$ indicates the range of the covariant maps. For example, when $m = 2$, the bisimplicial set whose $(n_1, n_2)$-simplices are the maps $(\Delta^{n_1})^{\mathrm{op}} \times \Delta^{n_m} \longrightarrow S$ is written as $\overline{\mathrm{Corr}}(S, [2, 2])$. In the case where covariant maps always land in a particular collection of arrows, we will substitute $[k_1, k_2]$ with a list of these collections. For example, in the case where $m = 3$ and the second and third factors are covariant and always land in $E_0$ and $E_1$, respectively, we write $\overline{\mathrm{Corr}}(S, E_0, E_1)$. In the case that there are no covariant maps, we will write $\overline{\mathrm{Corr}}(S)$ and if there are any strict subsets of the collection of all edges in $S$ that any of the factors always land.

Theorem 4.8 of [iii] argues that for $m = 2$ where the covariant factor lands in $E$, given a geometric setup $(C, E)$, we have a categorical equivalence $\mathrm{Corr}(C, E) \simeq \delta^*\overline{\mathrm{Corr}}(S, E)$. In fact, Theorem 4.10 of [iii] establishes that if every $f \in E$ can be decomposed as $f_1 f_0$ where $f_i \in E_i$, then we have a categorical equivalence $\delta^*\overline{\mathrm{Corr}}(S, E) \simeq \delta^*\overline{\mathrm{Corr}}(S, E_0, E_1)$ for simplicial sets.

**Theorem 2.2.8.** *There exists a map $f : \mathsf{Sup}_{\mathsf{J}}^{\otimes}(E, \mathfrak{X}^L)^{\mathrm{op}} \longrightarrow \mathsf{Fun}^{\otimes, \mathrm{lax}}(\mathrm{Corr}(\mathsf{Stk}_{S_|}^{\mathrm{op}}, E), \mathfrak{X}^L)$ which is a fully faithful map of $\infty$-categories.*

*Proof.* It suffices, applying Corollary 2.2.7, to show the existence of a fully faithful map of $\infty$-categories $\mathsf{Fun}_{\mathrm{Adj}}^{\otimes, \mathrm{lax}}((\mathsf{Stk}_{S_|}^{\mathrm{op}}, I \circ P), \mathfrak{X}^L) \longrightarrow \mathsf{Fun}^{\otimes, \mathrm{lax}}(\mathrm{Corr}(\mathsf{Stk}_{S_|}^{\mathrm{op}}, E), \mathfrak{X}^L)$.

Consider the trisimplicial set $\delta^*\overline{\mathrm{Corr}}(\mathsf{Stk}_{S_|})$ where the second factor lands in $I$ and the factor always lands in $P$. That is, maps are of the kind $(\Delta^1)^{\mathrm{op}} \times (\Delta^1)^{\mathrm{op}} \times (\Delta^1)^{\mathrm{op}} \longrightarrow \mathsf{Stk}_{S_|}$. We have the map that takes the diagonal $d : \delta^*\overline{\mathrm{Corr}}(\mathsf{Stk}_{S_|}) \longrightarrow \mathsf{Stk}_{S_|}^{\mathrm{op}}$. We claim that $d$ admits a section $s : \mathsf{Stk}_{S_|}^{\mathrm{op}} \longrightarrow \delta^*\overline{\mathrm{Corr}}(\mathsf{Stk}_{S_|})$. If this were the case, then we would obtain a fully faithful map $\mathsf{Fun}(\mathsf{Stk}_{S_|}^{\mathrm{op}}, \mathsf{CAlg}(\mathfrak{X}^L)) \longrightarrow \mathsf{Fun}(\delta^*\overline{\mathrm{Corr}}(\mathsf{Stk}_{S_|}), \mathsf{CAlg}(\mathfrak{X}^L))$ of $\infty$-categories. And, restricting ourselves, to those edges that are locally vertically left adjointable on $I$ and locally vertically right adjointable on $P$, we obtain $\mathsf{Fun}_{\mathrm{Adj}}^{\otimes, \mathrm{lax}}((\mathsf{Stk}_{S_|}^{\mathrm{op}}, I \circ P), \mathfrak{X}^L) \longrightarrow \mathsf{Fun}_{\mathrm{Adj}}^{\otimes, \mathrm{lax}}((\delta^*\overline{\mathrm{Corr}}(\mathsf{Stk}_{S_|}), I \circ P), \mathfrak{X}^L)$ which is yet again fully faithul.

On the other hand, by the description of $\mathsf{Fun}_{\mathrm{Adj}}^{\otimes, \mathrm{lax}}((\delta^*\overline{\mathrm{Corr}}(\mathsf{Stk}_{S_|}), I \circ P), \mathfrak{X}^L)$, there is a passage to adjoints on the second and third factors. This is expressed as the following map map $\mathsf{Fun}_{\mathrm{Adj}}^{\otimes, \mathrm{lax}}((\delta^*\overline{\mathrm{Corr}}(\mathsf{Stk}_{S_|}), I \circ P), \mathfrak{X}^L) \longrightarrow \mathsf{Fun}^{\otimes, \mathrm{lax}}(\delta^*\overline{\mathrm{Corr}}(\mathsf{Stk}_{S_|}, I, P), \mathfrak{X}^L)$. Applying Lemma 2.4.2 of [xxv], this is a categorical equivalence. But, as we have discussed, there is a categorical equivalence $\mathrm{Corr}(\mathsf{Stk}_{S_|}, E) \simeq \delta^*\overline{\mathrm{Corr}}(\mathsf{Stk}_{S_|}, I, P)$ of simplicial sets.

Let $s : \mathsf{Stk}_{S_|}^{\mathrm{op}} \longrightarrow \delta^*\overline{\mathrm{Corr}}(\mathsf{Stk}_{S_|})$ be the map assigning each $n$-simplex $\sigma : \Delta^n \longrightarrow \mathsf{Stk}_{S_|}^{\mathrm{op}}$ to the $n$-simplex $(\Delta^n)^{\mathrm{op}} \times (\Delta^n)^{\mathrm{op}} \times (\Delta^n)^{\mathrm{op}} \longrightarrow \mathsf{Stk}_{S_|}$ determined by the cube spanned in one direction by $\sigma^{\mathrm{op}}$ and in the next two directions by $n$-simplices that are constant on the $n$-vertex of $\sigma$. Furthermore, the faces for this cube are sliced by the diagonal into 2-simplices that are images of degeneracy maps. That is, the faces of the resulting cubes are the "obvious" Cartesian squares. For example, a morphism $f \in \mathsf{Stk}_{S_|}^{\mathrm{op}}$ is assigned to the cube spanned in one direction with $f^{op}$ and in the other two directions with $\mathrm{id}_{\partial_1 f^{\mathrm{op}}}$; notice that the (0,0,0) vertex of this cube is $\partial_0 f^{op}$. That $d \circ s \simeq 1$, is seen because the compositions in all other directions except the one in which the first factor of $s(\sigma)$ lands, are the compositions of edges in constant diagrams of vertices common to $\sigma$; hence $d(s(\sigma))$ is necessarily homotopic to $\sigma$. □

**Remark 2.2.9.** The result above, beyond the description of suprematic spaces with the same geometric content, does not depend on the choice of the geometric setup; that is, of the collection $E$. Therefore, mutatis mutandis, they hold as true when we replace $E$ with any of its subsets that constitutes a geometric setup. We will find this observation of some use in the next section.

**Remark 2.2.10.** It is worth mentioning that how we defined $\mathsf{Stk}_{S_{|\pi}}$ given a suprematic space $\pi^{op}$ affects how we defined $\mathsf{Sup}_{\mathsf{J}}^{\otimes}(E, \mathfrak{X}^L)$. In turn, this affects the scope of Theorem 2.2.8. If we take $\mathsf{Stk}_{S|_\pi}$, recalling the setup, simply as the image of $\mu_I \circ \pi$, then we only need property (1) of Definition 2.2.5 for the above results to hold; this amounts to having an essentially surjective () and hence condition (2) is "unconsciously" fulfilled. Because of this, whenever we have any collection of suprematic spaces with the same geometric content, a version of Theorem 2.2.8 always holds. In the case where property (1) holds but not property (2), Theorem 2.2.8 can always be weakened so that $f$ is (-1)-truncated.

# A "Universal" Six Functor Formalism

In this section, our main task is to prove Theorem C. To do so, we will utilize the theory of structured spaces developed in [iv]. Indeed, in this section, our choice to formulate suprematic spaces using structured spaces eventually pays off. We will also highlight how, in using suprematic spaces, the problem of 'efficacy' of the factorization of particular maps through animated stacks is resolved, so to speak. To be precise, we will show that up to "lax" diagrams of $\infty$-categories, there is an initial factorization. We will use small to mean $\mathbb{U}_1$-small.

## A "Universal" Self-Dualizing $\infty$-Category

In this subsection, we establish for given any small $\infty$-topos $\mathcal{X}$, a "universal" self-dual category in the following sense. There exists a fully faithful map $\chi : \mathcal{X} \longrightarrow L_{\text{ét}}(\mathcal{X}) \subseteq \mathsf{Pr}^L$ of $\infty$-categories such that given any map $f : \mathcal{X} \longrightarrow \mathcal{D} \subseteq \widehat{\mathsf{Cat}}_\infty$ where $\mathcal{D}$ is a self-dual $\infty$-category, there exists a map $\tilde{f} : L_{\text{ét}}(\mathcal{X}) \longrightarrow \mathcal{D}$ which is unique up to contractible choice, and such that we have the following 2-simplices of $\widehat{\mathsf{Cat}}_\infty$.

$$\begin{array}{ccc} & L_{\text{ét}}(\mathcal{X}) & \\ \chi \nearrow & & \searrow \tilde{f} \\ \mathcal{X} & \xrightarrow{f} & \mathcal{D} \end{array} \qquad \begin{array}{ccc} & L_{\text{ét}}(\mathcal{X}) & \\ \chi^{\text{op}} \nearrow & & \searrow \tilde{f} \\ \mathcal{X}^{\text{op}} & \xrightarrow{f^{\text{op}}} & \mathcal{D} \end{array}$$

We speak of $\chi : \mathcal{X} \longrightarrow L_{\text{ét}}(\mathcal{X})$ as not having exactly a universal property. This is because, as we shall see, the pullback functor $\chi^* : \mathsf{Fun}(L_{\text{ét}}(\mathcal{X}), \mathcal{D}) \longrightarrow \mathsf{Fun}(\mathcal{X}, \mathcal{D})$ is generally not a fully faithful map for all self-dual $\mathcal{D}$. However, it is a monomorphism of $\infty$-categories under some mild restrictions.

In what is to follow, we will use the terminology fully faithful in reference to maps of simplicial sets to mean that the maps are fully faithful once subjected to the map $\mathfrak{C}[\_] : \mathsf{Set}_\Delta \longrightarrow \mathsf{Cat}_\Delta$.

**Proposition 3.1.1.** *Suppose that $\mathcal{C}$ is a small $\infty$-category which admits a categorical equivalence $p : \mathcal{C} \longrightarrow \widehat{\mathcal{A}} \subseteq \mathsf{Pr}^L$. Then there exists a simplicial set $\mathcal{K}$ and a fully faithful map $\bar{p} : \mathcal{C} \longrightarrow \mathcal{K}$ factoring through $p$ such that:*

1. *For any map $g : \mathcal{C} \longrightarrow \mathcal{B}$ where $\mathcal{B}$ is self-dual, there exists up to contractible choice, a unique map $\tilde{g} : \mathcal{K} \longrightarrow \mathcal{B}$ such that $g \simeq \tilde{g} \circ \bar{p}$.*

2. *For any map $g^{op} : \mathcal{C}^{op} \longrightarrow \mathcal{B}$ where $\mathcal{B}$ is self-dual, there exists a unique map $\widetilde{g} : \mathcal{K} \longrightarrow \mathcal{B}$ such that $g^{op} \simeq \tilde{g} \circ (\bar{p})^{op}$.*

3. *There exist monomorphisms of simplicial sets $i : \widehat{\mathcal{A}} \hookrightarrow \mathcal{K}$ and $i^{op} : \widehat{\mathcal{A}}^{op} \hookrightarrow \mathcal{K}$.*

*Proof.* Let $E$ be the discrete simplicial set on the objects of $\widehat{\mathcal{A}}$. Consider $\mathcal{K} := \widehat{\mathcal{A}} \coprod_E \widehat{\mathcal{A}}^{\mathrm{op}}$ which is the pushout of the respective inclusions of $E$ along each other. It is immediate that (3) is satisfied. If we take any map $g : \mathcal{C} \longrightarrow \mathcal{B}$ as described in (1), then we instantly obtain another map $g^{\mathrm{op}} : \mathcal{C}^{\mathrm{op}} \longrightarrow \mathcal{B}$. In turn, since $p$ is a categorical equivalence, we obtain a map $f : \widehat{\mathcal{A}} \longrightarrow \mathcal{B}$ such that $f \circ p \simeq g$. By duality, we have $f^{\mathrm{op}} \circ p^{\mathrm{op}} \simeq g^{\mathrm{op}}$. Hence, by the universal property of pushouts, we obtain a unique map $\tilde{g} : \mathcal{K} \longrightarrow \mathcal{B}$ with the properties desired for (2) and (1). Moreover, composing the covariant projection of (3) with $p$, we obtain $\bar{p}$. Subsequently, $\bar{p}$ is fully faithful. □

**Corollary 3.1.2.** *Suppose that $\mathcal{C}$ is a small $\infty$-category which admits a categorical equivalence $p : \mathcal{C} \longrightarrow \widehat{\mathcal{A}} \subseteq \mathsf{Pr}^L$. Then there exists a self-dual $\infty$-category $\mathcal{K}^{\mathrm{An}} \subseteq \mathcal{P}^L$ and a map $\bar{p}^{\mathrm{An}} : \mathcal{C} \longrightarrow \mathcal{K}^{\mathrm{An}}$ factoring through $p$ such that:*

1. *For any map $g : \mathcal{C} \longrightarrow \mathcal{B}$ where $\mathcal{B}$ is self-dual, there exists, up to contractible choice, a unique map $\tilde{g}^{\mathrm{An}} : \mathcal{K}^{\mathrm{An}} \longrightarrow \mathcal{B}$ such that $g \simeq \tilde{g}^{\mathrm{An}} \circ \bar{p}^{\mathrm{An}}$.*

2. *For any map $g^{op} : \mathcal{C}^{op} \longrightarrow \mathcal{B}$ where $\mathcal{B}$ is self-dual, there exists up to contractible choice, a unique map $\widetilde{g}^{\mathrm{An}} : \mathcal{K}^{\mathrm{An}} \longrightarrow \mathcal{B}$ such that $g^{\mathrm{op}} \simeq \tilde{g}^{\mathrm{An}} \circ (\bar{p}^{\mathrm{An}})^{\mathrm{op}}$.*

3. *There exist (-1)-truncated maps of $\infty$-categories $i^{\mathrm{An}} : \widehat{\mathcal{A}} \hookrightarrow \mathcal{K}^{\mathrm{An}}$ and $(i^{\mathrm{An}})^{\mathrm{op}} : \widehat{\mathcal{A}}^{\mathrm{op}} \hookrightarrow \mathcal{K}^{\mathrm{An}}$.*

*Proof.* Lemma 2.2.5.2 of [i] allows one to take the inner fibrant replacement $\mathcal{K} \longrightarrow [\mathcal{K}]$ to obtain a categorically equivalent $\infty$-category. We label the latter as $\mathcal{K}^{\mathrm{An}}$. The composition of all relevant maps in Proposition 3.1.1 with the inner fibrant replacement gives the desired results. Furthermore, $\mathcal{K}$ is evidently self-dual, and hence the same is said of $\mathcal{K}^{\mathrm{An}}$. □

The descriptions above establish the pullback maps $\bar{p}^{\mathrm{An}}_! : \mathsf{Fun}(\mathcal{C}, \mathcal{D}) \longrightarrow \mathsf{Fun}(\mathcal{K}^{\mathrm{An}}, \mathcal{D})$ and $(\bar{p}^{\mathrm{An}^{\mathrm{op}}})_! : \mathsf{Fun}(\mathcal{C}^{\mathrm{op}}, \mathcal{D}) \longrightarrow \mathsf{Fun}(\mathcal{K}^{\mathrm{An}}, \mathcal{D})$ be the pullback of $\bar{p}^{\mathrm{An}}_!$ along $(\bar{p}^{\mathrm{An}^{\mathrm{op}}})_!$ in $\widehat{\mathsf{Cat}_\infty}$.

**Proposition 3.1.3.** *The projections $\mathsf{Fun}^+(\mathcal{K}^{\mathrm{An}}, \mathcal{D}) \longrightarrow \mathsf{Fun}(\mathcal{C}, \mathcal{D})$ and $\mathsf{Fun}^+(\mathcal{K}^{\mathrm{An}}, \mathcal{D}) \longrightarrow \mathsf{Fun}(\mathcal{C}^{\mathrm{op}}, \mathcal{D})$ are (-1)-truncated.*

*Proof.* Putting to use Proposition 1.2.8.3 of [i], it suffices to replace $\mathcal{K}^{\mathrm{An}}$ with $\mathcal{K}$. Note that the maps $\bar{p}^{\mathrm{An}}_!$ and $(\bar{p}^{\mathrm{An}^{\mathrm{op}}})_!$ are each a composition of a categorical equivalence and a monomorphism of simplicial cells. Therefore, both are (-1)-truncated. Proposition 5.5.6.12 of [i] gives the desired conclusion. □

**Proposition 3.1.4.** *Suppose that $\mathcal{X}$ is a small $\infty$-topos. Then there exists a categorical equivalence $p : \mathcal{C} \longrightarrow \widehat{\mathcal{A}} \subseteq \mathsf{Pr}^L$.*

*Proof.* This is a direct application of Remark 6.3.5.10. In particular, given a small $\infty$-topos $\mathcal{X}$, one takes into consideration the map $p : \mathsf{Fun}(\Delta^1, \mathcal{X}) \longrightarrow \mathcal{X}$ that is described by evaluating at the end point $\{1\}$. This map is a Cartesian fibration and is therefore classified by a map $\mathcal{X}^{op} \longrightarrow \widehat{\mathsf{Cat}_\infty}$. The codomain of this map is the projection of the $\infty$-category $(\mathcal{L}\mathsf{Top}_{\text{ét}})_{\mathcal{X}/}$ into $\mathcal{L}\mathsf{Top}_{\text{ét}}$. The latter is the $\infty$-category of small $\infty$-topoi and étale geometric morphisms (cf. [i] 6.3.5.3). It follows that $(\mathcal{L}\mathsf{Top}_{\text{ét}})_{\mathcal{X}/}$ is identifiable with a subcategory of $\mathsf{Pr}^L$. □

In the context of Corollary 3.1.2, we let $L_{\text{ét}}(\mathcal{X})$ be the free finite product completion of $\mathcal{K}^{\text{An}}$ obtained from Proposition 3.1.4. That is, the smallest subcategory of $\mathsf{Fun}(\mathcal{K}^{\text{An}}, \widehat{\mathcal{S}})$ that contains finite products and $\mathcal{K}^{\text{An}}$; by virtue of being self-dual, the coYoneda embedding provides an embedding $\mathcal{K}^{\text{An}} \hookrightarrow \mathcal{P}((\mathcal{K}^{\text{An}})^{\text{op}})$.

**Theorem 3.1.5.** *Suppose that $\mathcal{X}$ is a small $\infty$-topos. Then there exists an $\infty$-category $L_{\text{ét}}(\mathcal{X})$ admitting finite products and a fully faithful map $\chi : \mathcal{X} \longrightarrow L_{\text{ét}}(\mathcal{X})$ which factors through $(\mathcal{L}\mathcal{T}\text{op}_{\text{ét}})_{\mathcal{X}/}$ such that:*

1. *For any map $g : \mathcal{C} \longrightarrow \mathcal{D}$ where $\mathcal{D}$ is self-dual and admits finite products, there exists, up to contractible choice, a unique map $\widetilde{g} : L_{\text{ét}}(\mathcal{X}) \longrightarrow \mathcal{D}$ such that $g \simeq \widetilde{g} \circ \chi$.*

2. *For any map $g^{op} : \mathcal{X}^{op} \longrightarrow \mathcal{D}$ where $\mathcal{D}$ is self-dual and admits finite products, there exists up to contractible choice, a unique map $\widetilde{g} : L_{\text{ét}}(\mathcal{X}) \longrightarrow \mathcal{D}$ such that $g^{op} \simeq \widetilde{g} \circ \chi^{op}$.*

*Proof.* This follows from the universal property of the free finite product completion and by applying Corollary 3.1.2 and Remark 6.3.5.10 of [i]. The universal property of the free finite completion is derived from applying Proposition 5.3.6.2 of [i] to $(\mathcal{K}^{\text{An}})^{\text{op}}$ for the case of finite coproducts. □

**Proposition 3.1.6.** *Let $\mathsf{Fun}^{+}(L_{\text{ét}}(\mathcal{X}, \mathcal{D}))$ indicate the $\infty$-category of product preserving maps. Then the pullback functor induced by coYonneda embedding restricts to a categorical equivalence $\mathsf{Fun}^{+}(L_{\text{ét}}(\mathcal{X}), \mathcal{D}) \longrightarrow \mathsf{Fun}^{+}(\mathcal{K}^{An}, \mathcal{D})$ for all $\mathcal{D}$ that admit finite products.*

*Proof.* This is seen from the combination of Proposition 5.3.6.2 of [i] in the case of finite coproducts for $(\mathcal{K}^{\text{An}})^{\text{op}}$. □

**Corollary 3.1.7.** $\mathsf{Fun}^{+}(L_{\text{ét}}(\mathcal{X}), \mathcal{D}) \longrightarrow \mathsf{Fun}(\mathcal{C}, \mathcal{D})$ *and* $\mathsf{Fun}^{+}(L_{\text{ét}}(\mathcal{X}), \mathcal{D}) \longrightarrow \mathsf{Fun}(\mathcal{C}^{op}, \mathcal{D})$, *the projections induced by universal property of pullbacks are (-1)-truncated maps of $\infty$-categories.*

*Proof.* Take Proposition 3.1.3 and Proposition 3.1.6. □

### Towards the Motivic Dream

This subsection, although the final act of this paper, is presented with an eye toward a fuller investigation in the future. In particular, in relation to the overarching motivic program. Meanwhile, the main goal is to provide a proof of Theorem C. Furthermore, to present additional 6-functor formalisms that can be factored in the shape Theorem C suggests.

We commence by presenting results concerning suprematic spaces and pregeometries whose underlying $\infty$-categories are animated $S$-stacks. These results observe that suprematic spaces naturally produce such pregeometries. And in making use of the general properties of structured spaces, in our case obtained as geometric envelopes of the said pregeometries, we arrive at Theorem C. Therein, we are finally vindicated in the use of 6-functor formalisms taking values in $\infty$-topoi; which to begin with, does not stray too far from the norm. For instance, when looking at $\mathsf{Pr}^L$, a careful choice of a subcategory is sometimes again a presentable $\infty$-category. Hence, its $\infty$-category of large anima is an $\infty$-topos that approximates the $\infty$-topoi we have so far used; indeed, if we only consider the stable presentable $\infty$-categories among such a subcategory, we will have exactly one of the desired $\infty$-topoi.

We shall, by abusing language, refer to a coverage in an $\infty$-category by a representative covering.

**Proposition 3.2.1.** *Suppose that $\pi^{op} : \mathcal{T} \longrightarrow \mathcal{X}^L$ is a suprematic space. Let $\{U_i \longrightarrow X\}_{i \in I}$ be an admissible covering sieve of the Grothendieck topology on $\mathcal{T}$. Then $\{\underline{U_i}^{\pi^{op}} \longrightarrow \underline{X}^{\pi^{op}}\}_{i \in I}$ is a coverage on $\mathsf{Stk}_{S|_\pi}$ where $(\_)^{\pi^{op}} = \mu_I^{op} \circ \pi^{op}$.*

*Proof.* Taking Proposition 2.2.1 into consideration, it remains to show two things. First, showing that the $\overline{\Sigma}_I$-structure due on the essential image of the subcategory spanned by covering admissible maps, produces a subset of the collection $\overline{\Sigma}_I^{\triangleleft}$-structure that meets all the conditions of $P$ from Definition 2.1.1. And second, showing that each admissible covering map $U \longrightarrow X$ lands in $\xi_{\overline{P}}$ where $\overline{P}$ is $\overline{\Sigma}_I^{\triangleleft}$. Keeping in mind that $\mu_\pi^{op} : \pi^{op}(\mathcal{T}^{ad}) \longrightarrow \mathsf{Stk}_{S|_\pi}$ is left exact and taking into account Proposition 2.1.16, the first statement follows. Unraveling the definition of $\xi_{\overline{P}}$, we observe using property (2) of $\Sigma_I$-structures that it suffices to show property (3) of $\xi_{\overline{P}}$ in order to show the second statement. But this follows from the pasting law for pullbacks. $\square$

**Corollary 3.2.2.** *The subcategory of $\mathsf{Stk}_{S|_\pi}$ spanned by the image, under $(\_)^{\pi^{op}}$, of the admissible covering sieves of $\mathcal{T}$ is an admissibility structure on $\mathsf{Stk}_{S|_\pi}$.*

*Proof.* We know that since $\pi^{op} : \mathcal{T} \longrightarrow \mathcal{X}^L$ is a suprematic space, the subcategory spanned by the admissible covering sieves maps to a $\overline{\Sigma}_I$-structure. Therefore, the image of the morphisms therein is closed under retracts and contains all identities and equivalences; where the latter is guaranteed by Proposition 2.1.4 and recalling the facts established regarding $\xi_{\overline{P}}$ in the proof of the previous statement. The 2-out-of-3 property required of admissible morphisms follows from Proposition 2.1.7. Taking into account Proposition 3.2.1 gives the conclusion. $\square$

We lift the following definition from Definition 3.2.1 of [iv].

**Definition 3.2.3.** Let $\mathcal{T}$ and $\mathcal{T}'$ be pregeometries. A *transformation of pregeometries* from $\mathcal{T}$ to $\mathcal{T}'$ is a functor $F : \mathcal{T} \to \mathcal{T}'$ satisfying the following conditions:

1. $F$ preserves finite products.

2. $F$ carries admissible morphisms in $\mathcal{T}$ to admissible morphisms in $\mathcal{T}'$.

3. Let $\{u_\alpha : U_\alpha \to X\}$ be a collection of admissible morphisms in $\mathcal{T}$ which generates a covering sieve on $X$. Then the morphisms $\{F(u_\alpha) : F(U_\alpha) \to F(X)\}$ generate a covering sieve on $F(X) \in \mathcal{T}'$.

4. Suppose we are given a pullback diagram in $\mathcal{T}$:

$$\begin{array}{ccc} U' & \longrightarrow & U \\ \downarrow & & \downarrow f \\ X' & \longrightarrow & X \end{array}$$

where $f$ is admissible. Then the induced diagram

$$\begin{array}{ccc} F(U') & \longrightarrow & F(U) \\ \downarrow & & \downarrow{\scriptstyle F(f)} \\ F(X') & \longrightarrow & F(X) \end{array}$$

is a pullback square in $\mathcal{T}'$.

Notice that by the axioms of pregeometries, given a pregeometry $\mathcal{T}$, $\mathcal{T}^{\mathrm{ad}}$ is also a pregeometry and the inclusion $\mathcal{T}^{\mathrm{ad}} \subseteq \mathcal{T}$ is a transformation of pregeometries. We shall refer to $\mathcal{T}^{\mathrm{ad}}$, viewed in this way, as the *fine pregeometry*.

**Proposition 3.2.4.** *The subcategory of $\mathcal{T}^{ad}$ spanned by the collection of admissible morphisms generating the covering sieves form an admissibility structure on $\mathcal{T}^{\mathrm{ad}}$ whenever there exists a suprematic space $\phi^{op}: \mathcal{T} \longrightarrow \mathcal{X}^L$.*

*Proof.* This follows from property (1) and property (7) in Definition 2.1.10 once one remembers that for a suprematic space, the image of the relevant subcategory above forms a $\overline{\Sigma}_I$-structure. □

We will call $\mathcal{T}^{\mathrm{ad}}$, viewed with this admissibility structure, as the *coarse pregeometry*.

**Proposition 3.2.5.** *Suppose that $\pi^{op}: \mathcal{T} \longrightarrow \mathcal{X}^L$ is a suprematic space and that $\mathcal{T}^{ad}$ is the fine pregeometry. Then there exists an admissibility structure on $\mathsf{Stk}_{S|_\pi}$ such that $(\_)^{\pi^{\mathrm{op}}}: \mathcal{T}^{\mathrm{ad}} \longrightarrow \mathsf{Stk}_{S|_\pi}$ is a transformation of pregeometries.*

*Proof.* Proposition 3.2.1 enables us to define a topology on $\mathsf{Stk}_{S|_\pi}$ and tells us that condition (3) of Definition 3.2.3 is met. We may now define an admissiblity structure on $\mathsf{Stk}_{S|_\pi}$ as the coarsest admissibility structure on $\mathsf{Stk}_{S|_\pi}$ containing $(\underline{\mathcal{T}^{\mathrm{ad}}})^{\pi^{\mathrm{op}}}$ (cf. [iv] 1.2.8). Now, unraveling the definition of $(\_)^{\pi^{\mathrm{op}}}$, we see that it preserves products and pullbacks. □

Henceforth, we shall refer to the pregeometry described above as the *fine pregeometry induced by the suprematic space $\pi^{op}: \mathcal{T} \longrightarrow \mathcal{X}^L$*. In the case where we are dealing with objects of some $\mathsf{Sup}_{\mathcal{T}}^{\otimes}(\mathcal{X}^L)$, we shall speak of the *fine pregeometry induced by $\mathsf{Sup}_{\mathcal{T}}^{\otimes}(\mathcal{X}^L)$*.

**Proposition 3.2.6.** *Suppose that $\pi^{op}: \mathcal{T} \longrightarrow \mathcal{X}^L$ is a suprematic space and that $\mathcal{T}^{ad}$ is the coarse pregeometry. Then there exists an admissibility structure on $\mathsf{Stk}_{S|_\pi}$ such that $(\_)^{\pi^{\mathrm{op}}}: \mathcal{T}^{\mathrm{ad}} \longrightarrow \mathsf{Stk}_{S|_\pi}$ is a transformation of pregeometries.*

*Proof.* We imbue $\mathsf{Stk}_{S|_\pi}$ with the admissibility structure mentioned in Corollary 3.2.2. Now, Proposition 3.2.1 becomes exactly what we wish to show. □

In the manner done previously, we refer to the pregeometry described above as the *coarse pregeometry induced by the suprematic space $\pi^{op}: \mathcal{T} \longrightarrow \mathcal{X}^L$*. In the case where we are dealing with objects of some $\mathsf{Sup}_{\mathcal{T}}^{\otimes}(\mathcal{X}^L)$, we shall speak of the *coarse pregeometry induced by $\mathsf{Sup}_{\mathcal{T}}^{\otimes}(\mathcal{X}^L)$*.

**Proposition 3.2.7.** *Suppose that both $\mathcal{T}^{\mathrm{ad}}$ and $\mathsf{Stk}_{S|_\pi}$ are the fine pregeometries and that the map $(\_)^{\pi^{\mathrm{op}}}: \mathcal{T}^{\mathrm{ad}} \longrightarrow \mathsf{Stk}_{S|_\pi}$ fully faithful. Then the restriction functor $\mathrm{Str}_{\mathsf{Stk}_{S|_\pi}}(\mathcal{X}^L) \longrightarrow \mathrm{Str}_{\mathcal{T}^{\mathrm{ad}}}(\mathcal{X}^L)$ induced by $(\_)^{\pi^{\mathrm{op}}}$ is left adjoint.*

*Proof.* We observe that since $()^\pi : (\mathcal{T}^{\text{ad}})^{\text{op}} \longrightarrow \mathsf{Stk}_{S|_\pi}^{\text{op}}$ is fully faithful and $\mathsf{CAlg}(\mathcal{X}^L)$ contains all small colimits, the restriction functor $\mathsf{Fun}(\mathsf{Stk}_{S|_\pi}^{\text{op}}, \mathsf{CAlg}(\mathcal{X}^L)) \longrightarrow \mathsf{Fun}((\mathcal{T}^{\text{ad}})^{\text{op}}, \mathsf{CAlg}(\mathcal{X}^L))$ admits a left adjoint seeing that left Kan extensions along $(\_)^\pi$ are to be found. (cf. [i] 4.3.3.6). Furthermore, in general, we see that since the admissible morphisms generating the Grothendieck topology of $\mathsf{Stk}_{S|_\pi}$ are exactly the image of those generating that of $\mathcal{T}^{\text{ad}}$ (under $()^\pi$), right Kan extensions, as soon as they exist and preserve finite limits, produce $\mathsf{Stk}_{S|_\pi}$-structures. In our case, the right Kan extension is taken along a fully faithful map and hence preserves limits. Subsequently, the dual versions of the adjoint functors described at the beginning restrict to the relevant full subcategories. □

**Proposition 3.2.8.** *Suppose that both $\mathcal{T}^{\text{ad}}$ and $\mathsf{Stk}_{S|_\pi}$ are the coarse pregeometries and that all right Kan extensions along $(\_)^{\pi^{\text{op}}} : \mathcal{T}^{\text{ad}} \longrightarrow \mathsf{Stk}_{S|_\pi}$ exist. Then $\mathfrak{Str}_{\mathsf{Stk}_{S|_\pi}}(\mathcal{X}^L) \longrightarrow \mathfrak{Str}_{\mathcal{T}^{\text{ad}}}(\mathcal{X}^L)$, the restriction functor induced by $(\_)^{\pi^{\text{op}}}$, is left adjoint.*

*Proof.* The only benefit of a fully faithful $(\_)^{\pi^{\text{op}}} : \mathcal{T}^{\text{ad}} \longrightarrow \mathsf{Stk}_{S|_\pi}$ in the proof of the previous statement was the existence of right Kan extensions that preserve limits. In our case, since the entire admissibility structure on $\mathsf{Stk}_{S|_\pi}$ lies in the image of $(\_)^{\pi^{\text{op}}}$, any right Kan extension will preserve finite products. Furthermore, it will preserve pullbacks of admissible morphisms since the edges of the admissibility structure admit a $\overline{\Sigma}_I$-structure on $\mathcal{X}^L$. □

**Theorem 3.2.9.** *Suppose that both $\mathcal{T}^{\text{ad}}$ and $\mathsf{Stk}_{S|_\pi}$ are the coarse pregeometries. Then there exists, for any $\mathsf{Sup}_\mathcal{T}^\otimes(\mathcal{X}^L)$ having $\pi^{\text{op}}$ as an object, a monomorphism $\mathsf{Sup}_\mathcal{T}^\otimes(\mathcal{X}^L) \longrightarrow \mathfrak{Str}_{\mathsf{Stk}_{S|_\pi}}(\mathcal{X}^L)$ which maps suprematic spaces to suprematic spaces.*

*Proof.* In general, we can determine a map with codomain $\mathfrak{Str}_{\mathsf{Stk}_{S|_\pi}}(\mathcal{X}^L)$ from the subcategory of $\mathfrak{Str}_{\mathcal{T}^{\text{ad}}}(\mathcal{X}^L)$ spanned by $\mathcal{T}^{\text{ad}}$-structures admitting right Kan extensions along $(\_)^{\pi^{\text{op}}} : \mathcal{T}^{\text{ad}} \longrightarrow \mathsf{Stk}_{S|_\pi}$. This is a specialization of the argument used in the proof of Proposition 3.2.8. Now, to see that the image of a suprematic space is a suprematic space under this map, observe that Kan extensions recover the map being extended; the right Kan extension of any $\pi^{\text{op}} \in \mathsf{Sup}_\mathcal{T}^\otimes(\mathcal{X}^L)$, therefore, shares in all properties of $\pi^{\text{op}}$ when restricted to the admissibility structure on $\mathsf{Stk}_{S|_\pi}$. Proposition 2.2.6 and recalling that maps that preserve pullbacks preserve $n$-truncated maps give the remainder (cf. [xxvi] 9.4.3.21). □

**Remark 3.2.10.** We have talked before about the "efficiency" of factorizing suprematic spaces with the same *geometric content*. These are exemplified by objects of $\mathsf{Sup}_\mathcal{T}^\otimes(\mathcal{X}^L)$. What Theorem 3.2.9 tells us is that the map extending of any one of these objects along $(\_) : \mathcal{T}^{\text{ad}} \longrightarrow \mathsf{Stk}_{S|}$ is an object of $\mathfrak{Str}_{\mathsf{Stk}_{S|}}(\mathcal{X}^L)$. Applying Lemma 3.4.3 in combination with Proposition 1.5.1 of [iv], we deduce that the former has an initial object $\wp$. Therefore, given any $\pi^{op} \in \mathsf{Sup}_\mathcal{T}^\otimes(\mathcal{X}^L)$ there exists a map $h : \mathsf{Stk}_{S|} \times \Delta^1 \longrightarrow \mathcal{X}^L$ such that $h|\mathsf{Stk}_{S|} \times \{0\} = \wp \circ (\_)$ and $h|\mathsf{Stk}_{S|} \times \{1\} = \pi_0^{\text{op}} \circ (\_)$. Therefore, up to natural transformation, $\wp \circ (\_)$ is the most cost-effective factorization of all suprematic spaces belonging to $\mathsf{Sup}_\mathcal{T}^\otimes(\mathcal{X}^L)$.

**Proposition 3.2.11.** *Suppose that $\mathcal{T}$ is any pregeometry. Then there exists an $\infty$-category $\mathcal{M}$*

*that admits finite products and a fully faithful* $p : \mathcal{T} \longrightarrow \mathcal{M}$ *that admits a self-dual essential image and such that if* $\mathcal{O} : \mathcal{T} \longrightarrow \mathcal{X}^L$ *is an object of* $\mathrm{Str}_\mathcal{T}(\mathcal{X}^L)$, *the following are true.*

1. *There exists, up to contractible choice, a unique map* $\widetilde{\mathcal{O}} : \mathcal{M} \longrightarrow \mathcal{X}^L$ *which preserves finite products and such that* $\mathcal{O} \simeq \widetilde{\mathcal{O}} \circ p$.

2. *Given* $\mathcal{O}^{op} : \mathcal{T}^{op} \longrightarrow \mathcal{X}^L$, $\mathcal{O}^{op} \simeq \widetilde{\mathcal{O}} \circ p^{op}$.

3. *There exists* $\mathsf{Fun}^+(\mathcal{M}, \mathcal{X}^L) \subseteq \mathsf{Fun}(\mathcal{M}, \mathcal{X}^L)$ *such that the pullback functors induced by $p$ and $p^{op}$ restrict to monomorphisms of $\infty$-categories as follows:* $\mathsf{Fun}^+(\mathcal{M}, \mathcal{X}^L) \xrightarrow{p^*} \mathsf{Fun}(\mathcal{T}, \mathcal{X}^L)$ *and* $\mathsf{Fun}^+(\mathcal{M}, \mathcal{X}^L) \xrightarrow[p^{op*}]{} \mathsf{Fun}(\mathcal{T}^{op}, \mathcal{X}^L)$.

*Proof.* It suffices, by applying Theorem 3.1.5, to show an embedding $\mathcal{T} \hookrightarrow \mathcal{Y}$ where $\mathcal{Y}$ is a small $\infty$-topos, and the embedding has the following universal property: the pullback functor it induces, $\mathsf{Fun}^*(\mathcal{Y}, \mathcal{X}^L) \longrightarrow \mathrm{Str}_\mathcal{T}(\mathcal{X}^L)$, is a categorical equivalence. Here, $\mathsf{Fun}^*(\mathcal{Y}, \mathcal{X}^L)$ is spanned by edges in $\mathcal{L}\mathcal{T}\mathrm{op}$. But this is exactly what a combination of Lemma 3.4.3 and Proposition 1.4.2 of [iv] guarantees. $\square$

**Remark 3.2.12.** The proof above, using the references therein, implicitly puts to use the fact that any geometry $\mathcal{G}$ admits a *universal* $\mathcal{G}$-structure and that in the first place one can pass from pregeometries to geometries by taking *geometric envelopes*. In particular, the universal geometry in question is nothing more than the following familiar composition: $\mathcal{G} \xrightarrow{\ \ \ \ } \mathcal{P}(\mathcal{G}) \xrightarrow{L} \mathrm{Shv}(\mathcal{G})$ where $L$ is the left exact localization giving rise to sheaves of anima on the site $\mathcal{G}$. It follows from the definition that universal geometries are unique up to categorical equivalence.

**Proposition 3.2.13.** *Let* $p : \mathcal{T} \longrightarrow \mathcal{M}$ *be as discussed previously. Then the following is true.*

1. *The codomain of the dual map* $p^{op} : \mathcal{T}^{op} \longrightarrow \mathcal{M}$ *is a subcategory of* $\mathsf{CAlg}(\mathcal{M}^\otimes)$ *where we have considered the symmetric monoidal structure on* $\otimes : \mathcal{M} \longrightarrow \mathcal{F}\mathrm{in}_*$ *induced by the finite product.*

2. *For* $f \in \mathcal{T}$, $p(f)$ *admits a right adjoint* $p(f)^*$ *such that we have a square* $\Delta^1 \times \Delta^1 \longrightarrow \mathcal{M}$ *visualized as follows.*

$$\begin{array}{ccc} A \times B & \xrightarrow{p(f) \times 1} & B \times B \\ {\scriptstyle 1 \otimes p(f)^*} \downarrow & & \downarrow {\scriptstyle \otimes} \\ A & \xrightarrow{p(f)} & B \end{array}$$

3. *For* $f \in \mathcal{T}$, $p(f)$ *is a conservative map.*

*Proof.* Notice, by applying Theorem 3.1.5 that the image of $p^{op}$ lands in a full subcategory of $\mathcal{M}$ that is categorically equivalent to $(\mathcal{L}\mathcal{T}\mathrm{op}_{\mathrm{\acute{e}t}})_{\mathcal{Y}/}$ where $\mathcal{Y}$ is used as in the situation of the proof of the previous statement. Furthermore, since the objects of $(\mathcal{L}\mathcal{T}\mathrm{op}_{\mathrm{\acute{e}t}})_{\mathcal{Y}/}$ are endowed with a (closed) symmetric monoidal structure induced by finite products, and the maps between them preserve all limits, $p^{op}$ lands in $\mathsf{CAlg}(\mathcal{M}^\otimes)$ since its symmetric monoidal structure is also induced by finite products. It remains to observe Proposition 6.3.5.11 of [i]. $\square$

**Lemma 3.2.14.** *Suppose that $p : \mathcal{T} \longrightarrow \mathcal{M}$ is as discussed previously. Then there exists a geometric setup $(\mathcal{T}, E_{\text{ét}})$ such that the following holds true.*

1. *Suppose given any 2-simplex $\Delta^2 \longrightarrow \mathcal{T}$ depicted as*

$$\begin{array}{ccc} & B & \\ {}^{f}\nearrow & & \searrow^{g} \\ A & \xrightarrow{h} & C \end{array}$$

   *where $g \in E_{\text{ét}}$. Then $f \in E_{\text{ét}}$ if and only if $h \in E_{\text{ét}}$.*

2. *Suppose given a pullback square $\sigma : \Delta^1 \times \Delta^1 \longrightarrow \mathcal{T}$ resulting from taking the pullback of $f \in E_{\text{ét}}$ along some edge of $\mathcal{T}$. Then $p \circ \sigma$ is horizontally right adjointable. Equivalently, $(p \circ \sigma)^{op}$ is vertically left adjointable.*

*Proof.* Consider the essential image $\mathcal{B}$ of $p : \mathcal{T} \longrightarrow \mathcal{M}$ and all the squares therein. From among these squares, consider those that are horizontally right adjointable. Furthermore, take only those from these new squares that satisfy all properties of $P$ in Definition 2.1.1 in light of the fully faithful map $p^{-1} : \mathcal{B} \longrightarrow \mathcal{T}$. Now, set $E_{\text{ét}}$ to be the collection of all edges equivalent to some $f \in \xi_P(\mathcal{T}, p^{-1})$. By Proposition 2.1.3 and Proposition 2.1.5 we are done. $\square$

**Remark 3.2.15.** We should note that $\xi_P(\mathcal{T}, p^{-1})$ contains at least all equivalences. Furthermore, any fully faithful morphism that satisfies condition (3) is also an element of $\xi_P(\mathcal{T}, p^{-1})$. This is because, out of all the conditions of $P$, only the second half of condition (1) is not automatically met by horizontally right adjointable squares (cf. [xxviii] B.2). But this obstruction disappears once all the edges involved in setting up condition (1) are fully faithful morphisms.

**Proposition 3.2.16.** *The map $p^{op} : \mathcal{T}^{op} \longrightarrow \mathcal{M}$ upgrades to a map $\overline{p} \in \mathsf{Fun}^{\otimes,\text{lax}}(\text{Corr}(\mathcal{T}, E_{\text{ét}}), \mathcal{M})$.*

*Proof.* If we combine Lemma 3.2.14 and Proposition 3.2.13, we see that the geometric setup $(\mathcal{T}, E_{\text{ét}})$ meets the conditions of Proposition A.5.10 of [ii]. $\square$

Moving forward, we refer to any map $\mu : \mathcal{C}^{op} \longrightarrow \mathsf{CAlg}(\mathcal{D}^{\otimes})$ of $\infty$-categories meeting the conditions set out in Proposition A.5.10 of [ii] as a *Nagata setup*, as is done in [xxv]. Implicit in this formulation is a geometric setup $(\mathcal{C}, E_{\mu})$. We will also replace $\mathcal{M}$, as discussed previously, with $L_{\text{ét}}(\mathcal{T})$. The notation here is doubly suggestive. It references the fact that we are using the $\infty$-category $\mathcal{L}\mathsf{Top}_{\text{ét}}$ and the 'L' is suggestive of J. Lurie, whose work is extensively throughout the article. If need be, in words, we will say *L-étale extension* of $\mathcal{T}$ to describe $L_{\text{ét}}(\mathcal{T})$.

**Theorem 3.2.17.** *Suppose that $\mathcal{O}^{op} : \mathcal{T}^{op} \longrightarrow \mathcal{X}^L$ is a Nagata setup as well as the dual map of $\mathcal{O} \in \mathsf{Str}_{\mathcal{T}}(\mathcal{X}^L)$. Then if $p : \mathcal{T} \longrightarrow L_{\text{ét}}(\mathcal{T})$ is the map discussed before, there exists a geometric setup $(\mathcal{T}, E)$, a map $\overline{p} : \text{Corr}(\mathcal{T}, E) \longrightarrow L_{\text{ét}}(\mathcal{T})$, and a map $\mathcal{O}^{\square} : L_{\text{ét}}(\mathcal{T}) \longrightarrow \mathcal{X}^L$ such that if $\overline{\mathcal{O}} : \text{Corr}(\mathcal{T}, E) \longrightarrow \mathcal{X}^L$ is the map induced by the Nagata setup, $\overline{\mathcal{O}} \simeq \mathcal{O}^{\square} \circ \overline{p}$.*

*Proof.* We begin by setting $E = E_{\mathcal{O}^{op}} \cap E_{\text{ét}}$. Lemma 2.4.2 of [xxv] allows one to reduce the showing of the intended result to showing that both $\mathcal{O}^{op} : \mathcal{T}^{op} \longrightarrow \mathcal{X}^L$ and its dual $\mathcal{O} : \mathcal{T} \longrightarrow \mathcal{X}^L$ factor through $\mathcal{O}^\square : L_{\text{ét}}(\mathcal{T}) \longrightarrow \mathcal{X}^L$. But this is precisely the content of Proposition 3.2.11. □

**Theorem 3.2.18.** *There exists a geometric setup* $(\mathsf{Stk}_{S|}, \overline{E})$ *and a lax symmetric monoidal map* $\chi : \mathrm{Corr}(\mathsf{Stk}_{S|}, \overline{E}) \longrightarrow L_{\text{ét}}(\mathsf{Stk}_{S|})$, *whose image is a subcategory of* $\mathsf{Pr}^L$, *such that for any* $\mathcal{D} \in \mathsf{Fun}^{\otimes,\mathrm{lax}}(\mathrm{Corr}(\mathsf{Stk}_{S|}, \overline{E}), \mathcal{X}^L)$ *in the image of* $f : \mathsf{Sup}_{\mathcal{T}}^{\otimes}(\overline{E}, \mathcal{X}^L)^{op} \longrightarrow \mathsf{Fun}^{\otimes,\mathrm{lax}}(\mathrm{Corr}(\mathsf{Stk}_{S|}, \overline{E}), \mathcal{X}^L)$, *there exists a map* $\tilde{\mathcal{D}} : L_{\text{ét}}(\mathsf{Stk}_{S|}) \longrightarrow \mathcal{X}^L$ *so that* $\mathcal{D} \simeq \tilde{\mathcal{D}} \circ \chi$.

*Proof.* Theorem 3.2.9 tells us that when we consider the coarse pregeometry $\mathsf{Stk}_{S|}$ induced by $\mathsf{Sup}_{\mathcal{T}}^{\otimes}(\overline{E}, \mathcal{X}^L)$, then for every $\mathcal{D} \in \mathsf{Fun}^{\otimes,\mathrm{lax}}(\mathrm{Corr}(\mathsf{Stk}_{S|}, \overline{E}), \mathcal{X}^L)$ in the image of $f : \mathsf{Sup}_{\mathcal{T}}^{\otimes}(E, \mathcal{X}^L)^{op} \longrightarrow \mathsf{Fun}^{\otimes,\mathrm{lax}}(\mathrm{Corr}(\mathsf{Stk}_{S|}, \overline{E}), \mathcal{X}^L)$, $\mathcal{D}|\mathrm{Corr}(\mathsf{Stk}_{S|}, \mathrm{isom})$ is a $\mathsf{Stk}_{S|}$-structure. Therefore, meeting the requisites of 3.2.17. It remains to intersect the admissible morphisms of $\mathsf{Stk}_{S|}$ with $E_{\text{ét}}$ to obtain the result. Notice that for $E$ obtained in Proposition 2.2.1, $E \supseteq \overline{E}$. But Remark 2.2.9 ensures that $f$ exists in the first place. □

**Remark 3.2.19.** It is tempting to think that because we have a Nagata setup $\chi|\mathrm{Corr}(\mathsf{Stk}_{S|}, \mathrm{isom})$, any map $\mathcal{D}^\square \in \mathsf{Fun}^+(L_{\text{ét}}(\mathsf{Stk}_{S|}), \mathcal{X}^L)$ is automatically extended to a Nagata setup $\mathsf{Stk}_{S|}^{op} \longrightarrow \mathcal{X}^L$. It is, however, not the case that every vertex in $\mathsf{Fun}^+(L_{\text{ét}}(\mathsf{Stk}_{S|}), \mathcal{X}^L)$ preserves the data associated with adjoint functors.

**Remark 3.2.20.** It is stipulated in Theorem 5.2 of [xxix] that for a noetherian scheme $S$ and the $\infty$-category $\mathsf{N}(\mathrm{Sm}_S^{ft})$ of smooth finite type $S$-schemes, the $\infty$-category $\mathcal{SH}(S)$ is characterized by the universal property that it admits a fully faithful functor $\mathsf{Fun}^L(\mathcal{SH}(S), \mathcal{D}) \longrightarrow \mathsf{Fun}(\mathsf{N}(\mathrm{Sm}_S^{ft}), \mathcal{D})$ for any $\infty$-category $\mathcal{D}$ that admits small colimits. Here $\mathsf{Fun}^L(\mathcal{SH}(S), \mathcal{D})$ is the full subcategory of $\mathsf{Fun}(\mathcal{SH}(S), \mathcal{D})$ spanned by cocontinuous vertices. Furthermore, the essential image of the above map is the full subcategory of its target spanned by maps sending Nisnevich covers to effective epimorphisms and that also satisfy $\mathbb{A}^1$-invariance. Roughly speaking, as discussed in Remark 3.2.12, the construction of $L_{\text{ét}}(\mathsf{Stk}_{S|})$ comes down to the existence of a similar universal property; at least without taking $\mathbb{A}^1$-invariance into account. This universal property is the statement of Proposition 6.2.3.20 of [i]. At the same time, one should observe that $L_{\text{ét}}(\mathsf{Stk}_{S|})$ contains a full subcategory spanned by étale geometric morphism of the kind $\mathcal{S}hv(\mathsf{Stk}_{S|})/U \longrightarrow \mathcal{S}hv(\mathsf{Stk}_{S|})/V$ where both $U$ and $V$ are objects of $\mathsf{Stk}_{S|}$. Setting $S = \mathbb{Z}$, it will be nice to know how for any object $U \in \mathsf{Stk}_|$, the Grothendieck topology on $\mathsf{Stk}_|$ transfers to $\mathsf{Stk}_{S|}/U$ and, therefore, how $\mathcal{S}hv(\mathsf{Stk}_|)/U$ and $\mathcal{S}hv(\mathsf{Stk}_|/U)$ relate. In the case that these two $\infty$-categories are categorically equivalent, the universal property established in Proposition 6.2.3.20 allows us to see $L_{\text{ét}}(\mathsf{Stk}_{S|})$ as a kind of $\infty$-category of $\infty$-categories which could potentially be coincident with an $\infty$-category of objects of the kind $\mathcal{SH}(S)$.

In closing, it is worth mentioning that although one needs the service of Theorem 3.2.9 to arrive at Theorem 3.2.18, it is nevertheless noticeable that it is an easy consequence of Theorem 3.2.17. However, it stands as a bridge between the ideas and technologies employed throughout this article. Suprematic spaces, yet again justifying the architecture of their formulation, are seen weaving out of the three sections a coherent whole.

# Appendix

We prove the main (reconstruction) theorem of tensor triangulated geometry using topos theoretic methods. We will work exclusively with 1-categories. The main aim of this section is to provide hands-on experience of the groundwork needed to motivate the article. We hope to keep this section brief and will only present results to the extent that they are needed for this article.

We will work only with tensor triangulated categories. However, it should be noted that similar results for other relevant tensor exact categories hold as special variants of those obtained in the setting of the former. For example, those involving reconstruction of qcqs a scheme $X$ from $QC(X)$ or a noetherian scheme $X$ from $Coh(X)$. When necessary, we will utilize the results we obtained in the article regarding packeted prosites. And when confronting size issues, we will call a class that is a set *small* and *large* otherwise.

### Tensor Triangulated Categories

We understand by a tensor triangulated category, an essentially small triangulated category endowed with a symmetric monoidal structure that is compatible with the triangulated structure of the category. That is, both left and right tensoring preserve triangles. This implies that the tensor product commutes with direct sums. We will label such a category as $(\mathcal{K}, \otimes, 1)$ or simply as $\mathcal{K}$ where there is no danger of ambiguity.

We will call a functor $F : (\mathcal{K}, \otimes, 1) \longrightarrow (\mathcal{H}, \otimes, 1)$ *tensor triangulated* ($\otimes$-triangulated) if it is triangulated, additive and symmetric monoidal. We will often write "$\otimes$-triangulated" to mean the latter. Subsequently, we have a category $\text{TriCat}^{\otimes}$.

**Definition A.1.1.** A subcategory $A \subseteq \mathcal{K}$ is called triangulated if it is triangulated as a category. $A$ is called a *thick tensor ideal* if it has the following properties.

1. It is triangulated.
2. It is full, and the inclusion $A \hookrightarrow \mathcal{K}$ is an isofibration.
3. It admits finite direct sums and the inclusion into $\mathcal{K}$ preserves them.
4. It is idempotent complete.
5. Given $x \in A$ and $y \in \mathcal{K}$, $x \otimes y \in \mathcal{K}$.

If there is no danger of imprecision, we will drop "thick". Notice that given any $a \in \mathcal{K}$, the intersection of all thick tensor ideals is again another thick tensor ideal. We will label this thick tensor ideal as $[a]$ and say that it is generated by $a$. If $a = 1$, we observe that $[1] = \mathcal{K}$. That is, the tensor unit generates $\mathcal{K}$. Notice that we can define $[M]$ for any collection $M$ of objects of $\mathcal{K}$.

We say a tensor ideal $A$ is prime if and only if whenever $x \otimes y \in \mathcal{K}$, then $x \in A$ or $y \in A$.

**Definition A.1.2.** We will call a subcategory of $(\mathcal{K}, \otimes, 1)$ a $\Sigma$-structure if it has the following qualities.

1. It is a full triangulated tensor subcategory of $\mathcal{K}$. This means that if it contains two objects in a triangle of $\mathcal{K}$, then it contains all the objects of the triangle. Furthermore, that inherits the tensor structure from $\mathcal{K}$.
2. It admits finite direct sums and the inclusion into $\mathcal{K}$ preserves them.

3. It is idempotent complete.

**Definition A.1.3.** For objects $a$ and $b$ of $\mathcal{K}$, we say $a \leq b$ if whenever $a$ is contained in a $\Sigma$-structure, then so is $b$. It follows that the relation $\leq$ on the objects of $\mathcal{K}$ is a proset. We will label it $\mathcal{K}^{\leq}$.

**Proposition A.1.4.** $\mathcal{K}^{\leq}$ *is finitely complete.*

*Proof.* We notice that since each $\Sigma$-structure is idempotent complete, $a \oplus b$ is a direct sum of objects in $\mathcal{K}$, $a \oplus b \leq a$ and $a \oplus b \leq b$. In fact, if $c \leq a$ and $c \leq b$, then $c \leq a \oplus b$ since each $\Sigma$-structure admits direct sums which are coincident with direct sums in $\mathcal{K}$. Moreover, zero object is an object of each $\Sigma$-structure $\Rightarrow a \leq 0$ for all $a \in \mathcal{K}^{\leq}$. □

**Definition A.1.5.** Assume that $a$ is an object for which the thick tensor subcategory it generates is not $\mathcal{K}$. Let $\Sigma^{a,0}$ be a $\Sigma$-structure containing $a$ and for which the following properties are true:

1. It is isomorphic in $\text{TriCat}^{\otimes}$ to a thick tensor ideal of $A \subseteq \mathcal{K}$.

2. If $a \leq b$ and $b = x \otimes y$, then it contains either $x$ or $y$.

We will call such a $\Sigma$-structure *almost-prime* of $a$. We observe that $\Sigma^{a,0}$ is not uniquely defined for $a$ and that their collection is large. However, up to equivalence in $\text{TriCat}^{\otimes}$, this collection is small since $\mathcal{K}$ is essentially small. We will work up to equivalences of almost-primes in $\text{TriCat}^{\otimes}$. Furthermore, we will collect all almost-primes into a category whose morphisms are their inclusions into each other. We will label this category $\Sigma^{*,0}$.

Recall that the Balmer spectrum of a tensor triangulated category is the set of its prime tensor ideals (up to equivalence in $\text{TriCat}^{\otimes}$) given the topology generated by the sub-basis $\{U_a\}_{a \in \mathcal{K}}$ where $U_a = \{\mathfrak{p} : a \in \mathfrak{p}, \mathfrak{p} \text{ is a prime tensor ideal}\}$. For $\mathcal{K} \in \text{TriCat}^{\otimes}$, we will label its Balmer spectrum as $\text{Spec}(\mathcal{K})$.

**Definition A.1.7.** Let $\mathcal{C} = (\mathcal{K}, \leq)$ be a proset. A $J$-prime filter on $\mathcal{C}$ is a subset $\mathfrak{F} \subseteq ob(\mathcal{C})$ such that:

1. $\mathfrak{F}$ is non-empty.

2. $a \in \mathfrak{F}$ implies $b \in \mathfrak{F}$ whenever $a \leq b$.

3. For any $a, b \in \mathfrak{F}$ there exists $c \in \mathfrak{F}$ such that $c \leq a$ and $c \leq b$.

4. For any $J$-covering sieve $\{a_i \longrightarrow a\}_{i \in I}$ if $a \in \mathfrak{F}$ then there exists $i \in I$ such that $a_i \in \mathfrak{F}$.

**Proposition A.1.8.** *Suppose that $\mathcal{K}_{\alpha}^{\leq}$ is a sub-proset of $\mathcal{K}^{\leq}$ whose objects do not generate $\mathcal{K}$ and that it admits a packeting $\Gamma : (\mathcal{K}_{\alpha}^{\leq})^{op} \longrightarrow \Sigma^{*,0}$ such that for every $k \in \mathcal{K}_{\alpha}^{\leq}$, $\Gamma(k) \subseteq \mathcal{K}_{\alpha}^{\leq}$. Then the space of points of $\boldsymbol{Sh}(\mathcal{K}_{\alpha}^{\leq}, J^{\Gamma})$ embeds into $\text{Spec}(\mathcal{K})$.*

*Proof.* Recall Theorem 1.1.1 and observe that it suffices to show that the full subcategory of $\mathcal{K}$ containing the objects of any $J^{\Gamma}$-prime filter is prime. But since we are working up to equivalence in $\text{TriCat}^{\otimes}$, it suffices to show that the $J^{\Gamma}$-prime filters contain objects of prime tensor ideals. Let $\mathfrak{F}$ be such a prime filter. Recalling the definition of prime filters above, we see that $\mathfrak{F}$ contains direct summands by condition (2). Additionally, condition (3) tells us that $\mathfrak{F}$ is closed under direct sums. (this, in combination with what follows, becomes useful when showing that $\mathfrak{F}$ is closed under taking cones). Moreover, if $a \in \mathfrak{F}$ then all the objects of the thick tensor ideal generated by $[\Gamma(a)]$ are in $\mathfrak{F}$. To see this, take condition (4) and that any covering sieve $\{f : a_i \longrightarrow a\}_{i \in I}$ implies $a_i \leq a'$ for some $a'$ in $\Gamma^\lambda(a)$ and all $i \in I$. Applying condition (2) $\Rightarrow a' \in \mathfrak{F}$. The combination of closure under direct sums and the latter implies that the smallest full subcategory of $\mathcal{K}$ containing objects of $\mathfrak{F}$ is closed under taking cones.

Altogether, the above conclusions, with the non-emptiness criterion included, tell us that $\mathfrak{F}$, viewed as a set, is a set of objects of a thick tensor ideal. That this tensor ideal is prime, follows when we combine condition (4) of prime filters and property (2) of almost-primes. $\square$

**Proposition A.1.9.** *Suppose that $\mathfrak{p}_0 \subseteq \mathfrak{p}$ is a nesting of prime tensor ideals. Then there exists $\Gamma : (\mathcal{K}_\alpha^\leq)^{op} \longrightarrow \Sigma^{*,0}$ such that both $\mathfrak{p}$ and $\mathfrak{p}_0$ are $J^\Gamma$-prime filters.*

*Proof.* Consider $\mathfrak{p}^\leq$, the sub-proset of $\mathcal{K}^\leq$ determined the collection of all objects of $\mathfrak{p}$. And let $\Gamma : (\mathfrak{p}^\leq)^{op} \longrightarrow \Sigma^{*,0}$ acting as $a \longmapsto \mathfrak{p}$ for all $a \in \mathfrak{p} \backslash \mathfrak{p}_0$ and acting as $a \longmapsto \mathfrak{p}_0$ for all $a \in \mathfrak{p}_0$. Observe that this is well-defined since if $a \leq b$, then $a \in \mathfrak{p}_0$ only if $b \in \mathfrak{p}_0$. It is straightforward to see that $\mathfrak{p}$ is a $J^\Gamma$-prime filter. Checking Definition A.1.7 then gives that $\mathfrak{p}_0$ is a $J^\Gamma$-prime filter. $\square$

**Remark A.1.10.** The conclusion above tells us that in order to get a hold of all prime tensor ideals of $\mathcal{K}$, we need only look at packetings on sub-prosets of the kind $\mathfrak{p}^\leq$ where $\mathfrak{p}$ is a maximal prime ideal. We will utilize this fact when looking at $D_{\text{perf}}(X)$ where $X$ is a qcqs scheme with $j : X \backslash \overline{\{x\}}$ *quasi-perfect* for all closed points $x \in X$.

Let $\{\Gamma^\alpha : (\mathcal{K}_\alpha^\leq)^{op} \longrightarrow \Sigma^{*,0}\}_{\alpha \in \Omega}$ be the collection of all packetings, up to equivalence of prosets and natural isomorphisms of functors, of the kind described in Proposition A.1.7. Moreover, let the spaces of points resulting from the topoi of sheaves of sets on the sites induced by the packetings be labeled as such: $\{X_\alpha\}_{\alpha \in \Omega}$.

Next, set $X_{\alpha,\alpha'} = X_\alpha \cap X_{\alpha'}$ and fix $j_\alpha : X_{\alpha,\alpha'} \hookrightarrow X_\alpha$ as the obvious inclusion. Let $X_0$ be the coequalizer of the following diagram

$$\prod_{\alpha,\alpha' \in \Omega} X_{\alpha,\alpha'} \underset{j_{\alpha'}}{\overset{j_\alpha}{\rightrightarrows}} \prod_{\alpha \in \Omega} X_\alpha$$

**Theorem A.1.11.** *There exists an isomorphism $X_0 \longrightarrow \text{Spec}(\mathcal{K})$ of topological spaces.*

*Proof.* We observe that by the universal property of coequalizers, there exists a natural map $X_0 \longrightarrow \text{Spec}(\mathcal{K})$ since we are made aware by Proposition A.1.8 of an embedding $X_\alpha \hookrightarrow \text{Spec}(\mathcal{K})$ for each $\alpha \in \Omega$. Now the forgetful functor $U : \text{Top} \longrightarrow \text{Set}$ is left adjoint and hence coequalizers are computed in Set (cf. [x] 5.29.1). It follows that the natural map $X_0 \longrightarrow \text{Spec}(\mathcal{K})$ is an embedding. Proposition A.1.9 gives the desired result. $\square$

**Functoriality**

We examine the behavior of this construction under functors of TriCat$^\otimes$. In what is to follow, we will let $F : (\mathcal{K}, \otimes, 1) \longrightarrow (\mathcal{H}, \otimes, 1)$ be a $\otimes$-triangulated functor.

**Proposition A.2.1.** *$F$ updates to a map of prosets $F : \mathcal{K}^\leq \longrightarrow \mathcal{H}^\leq$*

*Proof.* This reduces to showing that $\Sigma$-structures are stable under pullbacks along $F$. And since $F$ preserves right exact and symmetric monoidal the result follows. $\square$

**Proposition A.2.2.** *$F$ preserves almost-primes under pre-images.*

*Proof.* It suffices to show $F$ preserves thick tensor ideals as well as prime thick tensor ideals. The latter is Lemma 7.2 of [vii] and since $F$ commutes with finite directs and is symmetric monoidal, preimage of a thick tensor ideal is also thick tensor ideal. Repleteness holds for any preimage of a replete subcategory for general categories and general functors. $\square$

**Corollary A.2.3.** *$F$ has target $(\mathcal{H}, \otimes, 1)$ $\otimes$-triangulated and $\Gamma : (\mathcal{H}_\alpha^\leq)^{op} \longrightarrow \Sigma^{*,0}$ is a packeting. Then $F^{-1}(\Gamma) : F^{-1}(\mathcal{H}_\alpha^\leq)^{op} \longrightarrow \Sigma^{*,0}$ is a packeting.*

*Proof.* Apply Proposition A.2.2 pointwise. $\square$

**Proposition A.2.4.** *Let $X_0(\mathcal{K})$ be the space we have obtained above for $\mathcal{K}$. Then $F$ induces a map $F^{op} : X_0(\mathcal{H}) \longrightarrow X_0(\mathcal{K})$*

*Proof.* We observe that for each packeting $F^{-1}(\Gamma) : F^{-1}(\mathcal{H}_\alpha^\leq)^{op} \longrightarrow \Sigma^{*,0}$, we obtain a map of prosites $(F^{-1}(\mathcal{H}_\alpha^\leq), J^{F^{-1}(\Gamma)}) \longrightarrow (\mathcal{H}_\alpha^\leq, J^\Gamma)$. This is because $F$ is left exact as a map of prosets and by the construction of the packets, we have $F(F^{-1}(\Gamma)(a)) \subseteq \Gamma(F(a))$. This means each $X_\alpha(\mathcal{H})$ admits a map into some $X_\beta(\mathcal{K})$. By the universal property of coequalizers we have desired map. $\square$

**Definition A.2.5.** This construction is due to Section 7 of [vii]. Let $\mathcal{X}_0$ be the collection of the points of $X_0$ seen as prime filters. Define $\chi^c : obj(\mathcal{K}) \longrightarrow 2^{\mathcal{X}_0}$ as $a \longmapsto \{\mathfrak{F} \in \mathcal{X}_0 : a \notin \mathfrak{F}\}$.

**Proposition A.2.6.** *For each open subset $U \subseteq X_0$, let $\mathcal{K}_U := \{a \in \mathcal{K} : \chi^c(a) \cap U = \varnothing\}$. Then $\mathcal{K}_U$ is thick tensor ideal.*

*Proof.* The condition met by $a \in \mathcal{K}_U$ is that for all $\mathfrak{F} \in U$, $a \in \mathfrak{F}$. But we have seen that $\mathfrak{F}$ are thick tensor ideals. But this implies that $\mathcal{K}_U$ is the intersection of all $\mathfrak{F} \in U$. $\square$

In Section 7 of [vii], one is able to use *only* the above fact to construct a sheaf of rings on $X_0$. Specifically, the Verdier quotient $\mathcal{K} \longrightarrow \mathcal{K}/\mathcal{K}_U$ is considered, and then the association $U \longmapsto \mathrm{End}_{\mathcal{K}/\mathcal{K}_U}(1)$ is made. This promotes $X_0$ to a ringed space. Furthermore, Lemma 7.2 of [vii] establishes that $\otimes$-triangulated maps are dual to maps of ringed spaces.

**Remark A.2.7.** We will find it worthwhile to observe that Proposition A.2.6 only uses the fact that all $\mathfrak{F} \in U$ are tensor ideals. Hence, it is enough to obtain a space made out of tensor ideals in order to get a ringed space out of a triangulated category. Moreover, the proof of Lemma 7.2 only uses general properties of functor and the fact that $\chi^c(F(a)) = F^{-1}(\chi^c(a))$ to show the functoriality of the ringed space construction. If we are able to construct spaces that meet these two conditions, then we are able to construct accompanying functors into RS, the category of ringed spaces and ringed maps.

Let $\Sigma(\mathcal{K})$ be the 2-category sub-prosets of $\mathcal{K}^{\leq}$ that correspond to radical tensor ideals of $\mathcal{K}$ and their inclusions as maps. In what follows, by abuse of notation, we will write $\mathcal{K}^{\leq}$ to mean all the objects of $\mathcal{K}^{\leq}$ that are contained in some prime tensor ideal.

**Proposition A.2.8.** *Let* $\Gamma : (\mathcal{K}^{\leq})^{op} \longrightarrow \Sigma(\mathcal{K})$ *be the map* $k \longmapsto \{a \in \cap_{\mathfrak{F} \in \mathcal{X}_0} \mathfrak{F} : k \in \mathfrak{F}\}$. *Then*

1. $\Gamma$ *is a packeting.*

2. *Given $X_1$, the space corresponding to the locale determined by the site $(\mathcal{K}^{\leq}, J^{\Gamma})$, there exists a topological embedding* $\mathrm{Spec}(\mathcal{K}) \hookrightarrow X_1$.

3. *There exists a ringed space $(X_1, \mathcal{O}_1)$ such that the embedding $\mathrm{Spec}(\mathcal{K}) \hookrightarrow X_1$ promotes to an embedding of ringed spaces.*

*Proof.* Notice that for each $a \in \mathcal{K}$, $\Gamma(a)$ is the intersection of all primes that contain $a$. Therefore, if $d \in \Gamma(a)$, $\Gamma(d) \subseteq \Gamma(a)$. Closer under pullbacks follows since it is an intersection of sub-prosets with pullbacks. That $\mathcal{K}^{\leq}$ is finitely complete is Proposition A.1.4. Next, we determine that each point of $X_1$ is semi-prime. This follows the same argument as used in Proposition A.1.9 except we note that we replace the prime condition with the semi-prime condition. Set theoretically, we obtain an inclusion $\mathrm{Spec}(\mathcal{K}) \hookrightarrow X_1$. Considering the sub-basis as determined in Theorem 1.1.1 and noticing that each semi-prime ideal that is not prime is contained in a prime (cf. Lemma 1.5 of [vii]), we obtain the topological embedding.

Following Remark A.2.7, we construct a sheaf of rings on $X_1$. However, we notice that if $\varnothing \neq U \subseteq V$ are open sub-basis subsets of $X_1$ and $U = \mathrm{Spec}(\mathcal{K}) \cap V$, then $\mathcal{K}_V = \mathcal{K}_U$ by the observation of Lemma 1.5 of [vii]. The only other case to check is where $U = \varnothing$. But this sends $V$ to the trivial ring once it is pulled back to $\mathrm{Spec}(\mathcal{K})$. Therefore, pointwise, we either have an isomorphism or a surjection of rings. Subsequently, we obtain an embedding ringed spaces. □

**Remark A.2.9.** The spaces obtained above are spectral; as we have shown regarding packeted prosites and as is known about the Balmer spectrum. Therefore, every characterization made regarding the spaces of points of the locales in question transfers to said locales.

### Scheme Theoretic Applications

In this section, we focus primarily on the tensor triangulated category, $(D_{\mathrm{perf}}(X), \otimes, 1)$ where we have placed an imposition on a qcqs scheme $X$. The main aim is to reinterpret the procedures of topos theory, performed in the construction of the Balmer spectrum above, as procedures initiated by certain kinds of functors in $\mathrm{TriCat}^{\otimes}$.

Grothendieck duality establishes that when $f : X \longrightarrow Y$ is a map of qcqs schemes, then $f_* : D_{\text{qc}}(X) \longrightarrow D_{\text{qc}}(Y)$ admits a right adjoint $f^! : D_{\text{qc}}(Y) \longrightarrow D_{\text{qc}}(X)$. We are interested in the situation where this occurence is detected at the level of the compact objects of the categories above. This is the case when $f^!$ commutes with direct sums.

**Definition A.3.1.** A map of qcqs schemes is called *quasi-perfect* if and only if $f^!$ commutes with direct sums.

**Example A.3.2.** We lift Theorem 1.2 from [xxxii]. For a map $f : X \to Y$ of qcqs schmes, the following are equivalent:

1. $f$ is quasi-perfect (resp. perfect).
2. $f$ is quasi-proper (resp. pseudo-coherent) and has finite Tor-dimension.
3. $f$ is quasi-proper (resp. pseudo-coherent) and $f^*$ is bounded.

**Remark A.3.3.** We are interested in this characterization because we would like to situate the construction above as an example of a general technique to obtain maps $\text{TriCat}^\otimes \longrightarrow \text{Stk}_\mathbb{Z} \longrightarrow \text{Top}$. So far, we have established that we need to specify for objects of $\text{TriCat}^\otimes$ prosets which are compatible with maps therein. In order to upgrade these prosets into prosites, we need to specify packetings. These packetings, are locally for objects of each prosets, intersections of thick tensor ideals and a prime ideal. And we have already determined a characterization of "primeness" in an ambient category that contains $\text{TriCat}^\otimes$ (cf. Definition 0.1.1). Now, tightening the definition therein by using strict pullbacks, we can specify a prime to be a fully faithful colocalization, $q : K \longrightarrow (D_{\text{perf}}(X), \otimes, 1)$ that obeys the projection formula $(q(x \otimes q^*y) \simeq qx \otimes y)$, and whose pullback along $\otimes : D_{\text{perf}}(X) \times D_{\text{perf}}(X) \longrightarrow D_{\text{perf}}(X)$ is as in Definition 0.1.1. Let $\mathcal{M}$ be the collection of all primes in $D_{\text{perf}}(X)$. Henceforth, we will assume that the proset on $D_{\text{perf}}(X)$ is as described for the general case of triangulated categories.

In what is to follow, we will assume that $X$ is a qcqs scheme where for each closed point $x \in X$, the inclusion $j : X/\overline{\{x\}} \hookrightarrow X$ is quasi-perfect.

**Theorem A.3.4.** *Let $\Gamma : (D_{\text{perf}}(X)^{\leq})^{\text{op}} \longrightarrow \Sigma(D_{\text{perf}}(X))$ be the map $k \longmapsto \{a \in \cap_{\mathfrak{p} \in \mathcal{M}} \mathfrak{p} : k \in \mathfrak{p}\}$. Then*

1. *$\Gamma$ is a packeting.*
2. *Given $X_1$, the space corresponding to the locale determined by the site $(D_{\text{perf}}(X)^{\leq}, J^\Gamma)$, there exists a topological embedding $\text{Spec}(D_{\text{perf}}(X)) \hookrightarrow X_1$.*
3. *There exists a ringed space $(X_1, \mathcal{O}_1)$ such that the embedding $\text{Spec}(D_{\text{perf}}(X)) \hookrightarrow X_1$ promotes to an embedding of ringed spaces.*

*Proof.* This follows from combining Proposition A.2.8 with the fact that maximal prime ideals correspond to kernels of localizations of $D_{\text{perf}}(X)$ induced by the closed immersions $\overline{\{x\}} \hookrightarrow D_{\text{perf}}(X)$ where $x$ is closed; and in our setup, such kernels are recovered as primes owing to the fact that $j : X/\overline{\{x\}} \hookrightarrow X$ is quasi-perfect. □

## References

[i]. Lurie, Jacob. *Higher Topos Theory*. Annals of Mathematics Studies 170. Princeton, NJ: Princeton University Press, 2009.

[ii]. Mann, Lucas. *A p-Adic 6-Functor Formalism in Rigid-Analytic Geometry*. Preprint, June 2022. arXiv:2206.02022.

[iii]. Scholze, Peter. *Six-Functor Formalisms*. Lecture notes, Winter 2022/23. Accessed March 30, 2025. Six Functors.pdf.

[iv]. Lurie, Jacob. *Derived Algebraic Geometry V: Structured Spaces*. May 4, 2009. arXiv:0905.045.

[v]. Balmer, Paul. *Presheaves of Triangulated Categories and Reconstruction of Schemes*. May 28, 2002. arXiv:math/0111049v2.

[vi]. Balmer, Paul. *The Spectrum of Prime Ideals in Tensor Triangulated Categories*. September 22, 2004. arXiv:math/0409360v2.

[vii]. Buan, Aslak Bakke, Henning Krause, and Øyvind Solberg. *Support Varieties—An Ideal Approach*. August 20, 2005. arXiv:math/0508379v1.

[viii]. Garkusha, Grigory. *Classifying Finite Localizations of Quasicoherent Sheaves*. St. Petersburg Mathematical Journal 21, no. 3 (2010): 433–458.

[ix]. Balmer, Paul, Ivo Dell'Ambrogio, and Beren Sanders. *Grothendieck-Neeman Duality and the Wirthmüller Isomorphism*. December 12, 2015. arXiv:1501.01999.

[x]. The Stacks Project Authors. *The Stacks Project*. Stacks Project. Accessed April 14, 2025.

[xi]. Caramello, Olivia. *A Topos-Theoretic Approach to Stone-Type Dualities*. March 17, 2011. arXiv:1103.3493v1.

[xii]. Pepin Lehalleur, Simon. *Six Functor Formalism; Why Solid Abelian Groups?* November 11, 2022. Accessed April 14, 2025. Notes: Six Functors.

[xiii]. Lurie, Jacob. *Tannaka Duality for Geometric Stacks*. March 23, 2005. arXiv:math/0412266.

[xiv]. Fukuyama, Hiroshi, and Isamu Iwanari. *Monoidal Infinity Category of Complexes from Tannakian Viewpoint*. September 27, 2012. arXiv:1004.3087.

[xv]. Balmer, Paul. *A Guide to Tensor-Triangular Classification*. December 19, 2019. Available at arXiv:1912.08963.

[xvi]. Lurie, Jacob. *Higher Algebra*. August 3, 2012. Available at lurie/papers/HA2012.

[xvii]. nLab Authors. *Structured $(\infty,1)$-Topos*. Accessed May 1, 2025. Structured $(\infty,1)$-topos

[xviii]. Aoki, Ko. *The Sheaves-Spectrum Adjunction*. February 8, 2023. arXiv:2302.04069v1.

[xix]. Levine, Marc. *Six Lectures on Motives*. Lecture notes from the ICTP Workshop on K-theory and Motives, May 14–25, 2007. Accessed May 14, 2025. LevineMotiveLecture.

[xx]. Drew, Brad, and Martin Gallauer. 2022. *The Universal Six-Functor Formalism*. Annals of K-Theory 7: 599–649. The Universal Six-Functor Formalism.

[xxi]. Johnstone, Peter T. *Stone Spaces*. Cambridge Studies in Advanced Mathematics, no. 3. Cambridge: Cambridge University Press, 1982.

[xxii]. Hochster, Melvin. *Prime Ideal Structure in Commutative Rings.* Transactions of the American Mathematical Society 142: 43–60. 1969. doi.org/10.1090/S0002-9947-1969-0251026-X.

[xxiii]. Aoki, Ko. *Tensor Triangular Geometry of Filtered Objects and Sheaves.* Mathematische Zeitschrift 303: 62. 2023. doi.org/10.1007/s00209-023-03210-z.

[xxiv]. Porta, Mauro, and Yue Yu, Tony. *Derived Non-Archimedean Analytic Spaces.* Selecta Mathematica, New Series 24 (2):609–665. 2018. doi.org/10.1007/s00029-017-0310-1.

[xxv]. Kuijper, Josefien. *An Axiomatization of Six-Functor Formalisms.* September 7, 2024. Available at arXiv:2309.11449v3.

[xxvi]. Lurie, Jacob. *Kerodon. Online Textbook: Homotopy-Coherent mathematics.* Accessed June 18, 2025. Kerodon.

[xxvii]. Tarrío, L. Alonso, López A. Jeremías, and Sancho de Salas Fernando. *Relative Perfect Complexes.* Mathematische Zeitschrift 304, no.3, 42. June 2023. doi.org/10.1007/s00209-023-03294-7.

[xxviii]. Porta, Mauro. 2015. *Derived Complex Analytic Geometry II: Square-Zero Extensions.* July 23, 2015. arXiv:1507.06602 [math.AG].

[xxix]. Robalo, Marco. *Noncommutative Motives I: A Universal Characterization of the Motivic Stable Homotopy Theory of Schemes.* June 17, 2013. arXiv:1206.3645v3 [math.AG].

[xxx]. Porta, Mauro, Emily Riehl, and Dominic Verity. *Fibrations and Yoneda's Lemma in an $\infty$-Cosmos.* Journal of Pure and Applied Algebra 221, no.3 (March): 499–564. 2017. doi.org/10.1016/j.jpaa.2016.07.003.

[xxxi]. Lurie, Jacob. *($\infty$,2)-Categories and the Goodwillie Calculus I.* 2009. arXiv:0905.0462.

[xxxii]. Lipman, Joseph, and Amnon Neeman. *Quasi-Perfect Scheme-Maps and Boundedness of the Twisted Inverse Image Functor.* Illinois Journal of Mathematics 51, no. 1: 209–236. 2007. https://doi.org/10.1215/.